\documentclass[12pt,leqno]{amsbook}
\usepackage{amssymb,amsmath,amsfonts,amstext,amscd,latexsym}
\def\ds{\displaystyle}
\def\beginproof{\noindent{\it Proof. }}

\newtheorem{theorem}{Theorem}[section]
\newtheorem{proposition}[theorem]{Proposition}
\newtheorem{example}[theorem]{Example}
\newtheorem{lemma}[theorem]{Lemma}
\newtheorem{corollary}[theorem]{Corollary}
\newcommand{\CC}{{\Bbb C}}

\newcommand{\R}{{\mathbb{R}}}
\newcommand{\N}{{\mathbb{N}}}
\newcommand{\D}{{\mathbb{D}}}
\newcommand{\G}{{\mathbb{G}}}
\newcommand{\Om}{{\Omega}}
\renewcommand{\O}{{\mathcal{O}}}
\newcommand{\CY}{{\mathcal C}}
\newcommand{\HH}{{\mathcal H}}
\newcommand{\eps}{\varepsilon}

\newcommand{\Ima}{\operatorname{Im}}

\newcommand{\ad}{\operatorname{ad}}

\newcommand{\rank}{\operatorname{rank}}
\newcommand{\tr}{\operatorname{tr}}
\newcommand{\spe}{\operatorname{sp}}
\newcommand{\Aut}{\operatorname{Aut}}
\newcommand{\Jac}{\operatorname{Jac}}

\def\PSH{\operatorname{PSH}}
\def\SH{\operatorname{PSH}}
\def\diag{\operatorname{diag}}

\def\C{\Bbb C}
\def\G{\Bbb G}
\def\D{\Bbb D}
\def\M{\mathcal M}
\def\N{\Bbb N}

\def\phi{\varphi}
\def\wdtl{\widetilde}
\def\wdht{\widehat}
\def\CD{\mathcal D}
\def\O{\mathcal O}
\def\Om{\Omega}
\def\eps{\varepsilon}
\def\ds{\displaystyle}
\def\diam{\operatorname{diam}}

\def\supp{\operatorname{supp}}

\def\span{\operatorname{span}}
\def\id{\operatorname{id}}
\def\bs{\boldsymbol}
\def\Sp{|\bs{p}|}

\def\phi{\varphi}
\def\dist{\operatorname{dist}}
\def\Ree{\operatorname{Re}}

\def\Ree{\operatorname{Re}}
\def\Re{\operatorname{Re}}

\def\Im{\operatorname{Im}}
\def\CO{\mathcal O}
\def\OO{\mathcal O}
\def\CA{\mathcal A}
\def\C{\Bbb C}
\def\ord{\operatorname{ord}}
\def\vol{\operatorname{vol}}

\def\NN{\mathbb N}

\def\RR{\mathbb R}
\def\CC{\mathbb C}
\def\DD{\mathbb D}

\begin{document}
\thispagestyle{empty}

\centerline{\large\bf NIKOLAI NIKOLOV}
\bigskip\bigskip
\centerline{\large\bf INVARIANT FUNCTIONS AND METRICS}
\bigskip
\centerline{\large\bf IN COMPLEX ANALYSIS} \vfill

Institute of Mathematics and Informatics,

Bulgarian Academy of Sciences,

Acad.~G.~Bonchev 8,

1113 Sofia, Bulgaria
\smallskip

E-mail address: nik@math.bas.bg

\tableofcontents
\chapter*{Introduction}\label{intr}
\setcounter{page}{1} \hfill{\it To Veneta} \vspace{2mm}

One of the most beautiful and important results in the classical
complex analysis is the Riemann Mapping Theorem stating that any
nonempty simply connected open subset of the complex number plane,
other than the plane itself, is biholomorphic to the open unit
disk $\D\subset\C$. On the other hand, H. Poincar\'e (1907) has
proven that the groups of (holomorphic) automorphisms of the open
polydisc and of the open ball in $\CC^2$ are not isomorphic; hence
these two topologically equivalent domains are not
biholomorphically equivalent. Therefore it is important that any
domain $D$ in $\CC^n$ can be associated with some
biholomorphically equivalent object. Generalizing the Schwarz-Pick
Lemma, C. Carath\'eo\-dory (1926) provided the first example of
such an object, different from the automorphism group; this object
was later called the Carath\'eodory pseudodistance. That is the
largest Poincar\'e distance between the images of two points from
$D$ under the all holomorphic mappings from $D$ into $\D.$
Somewhat later (1933) S. Bergman started to consider the
generating kernel of the Hilbert space of square-integrable
holomorphic functions on $D$ with the natural Hermitian metric and
distance (later his name is given to these three invariants). In
1967 S. Kobayashi introduced a pseudodistance, dual in some sense
to the Carath\'eodory pseudodistance. More precisely, it is the
greatest pseudo\-distance not exceeding the so called Lempert
function -- the infimum of the Poincar\'e distances between
preimages of pairs of points from $D$ under an arbitrary
holomorphic mapping from $\D$ to $D.$

In Chapter \ref{chap.lem} we discuss the basic properties of
various invariant functions and their infinitesimal forms called
(pseudo)metrics.

The estimates and the limit behavior of invariant
(pseudo)distances (or more generally, of functions) and
(pseudo)metrics, as well as Bergman kernel, play an important role
in numerous problems from complex analysis like asymptotic
estimates of holomorphic functions (of various classes),
continuation of holomorphic mappings, biholomorphic
(non)\-equivalence of domains, description of domains with
noncompact groups of automorphisms etc. (see e.g.
\cite{Jar-Pfl1,Isa-Kra,Kra2,Pol-Sha}). We will only mention that
one of the basic points in the classification theorem of bounded
convex domains of finite type in $\C^n$ with noncompact groups of
automorphisms (see \cite{BP}) is an estimate for the Kobayashi and
Carath\'eodory pseudodistances (see also Proposition
\ref{bound0}). In Chapter \ref{chap.bound} we obtain estimates of
these metrics, as well as the Bergman kernel and Bergman metric of
the so-called $\C$-convex domains.

Let us note that the exact calculation of some of the invariants
or finding estimates thereof leads e.g. to criteria for
solvability of correspon\-ding interpolation problems or to
restrictions on the solvability. The research in Chapter
\ref{chap.spec} is partially motivated by two examples of such
types of problems.

What follows is a brief review of the contents of the chapters.
\vspace{2mm}

{\bf Chapter \ref{chap.lem}.} Its aim is the introduction of the
basic invariant functions, distances and metrics together with
their basic properties.

The Lempert function $l_M$ of a given complex manifold $Ì$ is the
greatest holomorphically contractible function (i.e. decreasing
under holomorphic mappings), coinciding with the M\"obius distance
$m_\D$ on the unit disc $\D$. The Kobayashi function $k^\ast_M$ is
the greatest pseudodistance not exceeding $l_M$. The
Carath\'eodory function $c^\ast_M$ is the least holomor\-phically
contractible function. Then $k_M=\tanh^{-1}k^\ast_M$ and
$c_M=\tanh^{-1}c^\ast_M$ are the Kobayashi and Carath\'eodory
(pseudo)distances, respectively (coinciding on $\D$ with the
Poincar\'e distance $p_\D$). The (pseudo)metrics of Kobayashi
$\kappa_M,$ of Kobayashi--Buseman $\hat\kappa_M$ and of
Carath\'eodory $\gamma_M$ are the infinitesimal forms of $l_M,$
$k_m$ and $c_M,$ respectively. In a natural way are defined the
Lempert functions $k^{(m)}_M$ of higher order that are between
$\tanh^{-1}l_D(=k^{(1)}_M)$ and $k_M=(=k^{(\infty)}_M),$ as well
as their infinitesimal forms -- the Kobayashi metrics
$\kappa_M^{(m)}$ of higher order.

In Section \ref{lem-der} we note that the objects under
consideration are upper semicontinuous (see Proposition
\ref{der.pr0} and the comment preceding it). The main result in
this section, namely Theorem \ref{der.th1}, states that if $z\in
M$ and the function $\kappa_M$ is continuous and positive in
$(z;X)$ for each nonzero vector $X$, then the "derivative"\ of
$k^{(m)}_M$ at $z$ in the direction of $X$ coincides with
$\kappa^{(m)}_M(z;X).$ An essential step in the proof is
Proposition \ref{der.pr2} claiming that the "upper derivative"\ of
$k^{(m)}_M$ does not exceed $\kappa^{(m)}_M$ in the general case.
Theorem \ref{der.th1} generalizes some results of M.-Y. Pang
\cite{Pang} and M. Kobayashi \cite{KobM} concerning taut manifolds
(domains). We provide examples to demonstrate that the assumptions
in the theorem are essential.

In Section \ref{lem-bal} we find some relationships between the
Minkowski functions of a balanced domain or of its
convex/holomorphic hull and some of the previously defined
biholomorphic invariants of that domain whenever one of their
arguments is the origin. Some of these relationships are used in
the subsequent chapter.

In Section \ref{kob-buz} we prove that the Kobayashi--Buseman
metric $\hat\kappa_M$ equals the Kobayashi metric
$\kappa_M^{(2n-1)}$ of order $2n-1$ and this number is the least
possible in the general case. A similar result for $2n$ instead of
$2n-1$ can be found in the paper \cite{KobS1} of S. Kobayashi,
where $\hat\kappa_M$ is introduced.

In Section \ref{arak} we prove a general statement, Theorem
\ref{mer.th4}, for approximation and interpolation over the
so-called Arakelian sets. To this aim we use a well-known
interpolation-approximation result of P. M. Gauthier and W.
Hengartner \cite{Gau-Hen}, and A. Nersesyan \cite{Ner}.

This theorem is in the base of the proof of Theorem \ref{decr.th1}
from the subsequent section stating that the so-called generalized
Lempert function (of a given domain) does not decrease under
addition of poles. The last assertion is proven by Wikstr\"om
\cite{Wik1} for convex domains; he had left the general case as an
open question in \cite{Wik2}.

In Section \ref{desc} we discuss the product property of the
generalized Lempert function in order to reject a hypothesis of D.
Coman \cite{Com2} for equality between this function and the
generalized pluricomplex Green function. In Proposition
\ref{prod1} we find a necessary and sufficient condition for the
Lempert function of the bidisc with (fixed argument and) poles in
the cartesian product of two two-point subsets of
$\D_\ast=\D\setminus\{0\}$ to equal each of the two corresponding
functions of $\D.$
\vspace{2mm}

{\bf Chapter \ref{chap.spec}.} To understand better the geometry
of the so-called symmetrized polydisc we need some notions for
complex convexity of domains and their interrelations; this is the
aim of the first part of Section \ref{con}. In \cite{APS,Hor2} one
can find a detailed discussion of their role and application. We
only note that the $\C$-convexity is closely related to some
important properties of the Fantappi\`e transformation, and, as a
deep conclusion, to the question of solvability of linear PDEs and
in the class of holomorphic functions. A domain $D\subset\C^n$ is
called $\Bbb C$-convex if its nonempty intersections with complex
lines is connected and simply connected. Some other notions for
complex convexity of $D$ are the linear one (each point in the
complement of $D$ is contained in a complex hyperplane, disjoint
from $D$), the weak linear one (the same, but for the points in
$\partial D$) and the weak locally linear one. $\C$-convexity
implies weak convexity, and all the four notions coincide for
bounded domains with $\mathcal C^1$-smooth boundaries. In the
general case their place is between convexity and pseudoconvexity.

Indeed, in Corollary \ref{c.pr00} we give an affirmative answer to
a question of D. Jacquet \cite[p. 58]{Jac2} whether each weakly
locally linearly convex domain is pseudoconvex. On the other hand,
in Proposition \ref{c.pr1} we show that each weakly linearly
convex balanced domain is convex, which strengthens the same
observation for complete Reinhardt domains in \cite[Example
2.2.4]{APS}.

Theorem \ref{c.pr4} implies that a $\C$-convex domain is either a
cartesian product of $\C$ and another $\C$-convex domain, or is
biholomorphic to a bounded domain; in the latter case it is
$c$-finitely compact (i.e. the balls with respect to the
Carath\'eodory distance are relatively compact). This generalizes
the result of T. J. Barth in \cite{Bar} for convex domains.

Let $\mathcal L$ denote the class of domains $D$ such that the
least and the largest invariant functions of $D$ (from
complex-analytic viewpoint), namely Carath\'eodory function
$c^\ast_D$ and the Lempert function $l_D$, coincide (and in
particular they coincide with the Kobayashi function $k^\ast_D$).

Recently D. Jacquet \cite{Jac1} proved that each bounded $\Bbb
C$-convex domain with $\mathcal C^2$-smooth boundary can be
exhausted by $\Bbb C$-convex domains with ($\mathcal
C^\infty$-)smooth boundaries. Then the fundamental Lempert theorem
\cite{Lem1,Lem2} can be formulated like this:
\smallskip

{\it Each bounded $\Bbb C$-convex domain with $\mathcal
C^2$-smooth boundary belongs to the class $\mathcal L.$}
\smallskip

This property is carried over convex domains, as they can be
exhausted with smooth (even strictly) smooth domains. It was an
open question whether a bounded pseudoconvex domain from $\mathcal
L$ must be biholomorphic to a convex domain \cite{Zna,Jar-Pfl3}. A
recently found counterexample is the so-called symmetrized bidisc
$\Bbb G_2\subset\Bbb C^2$ -- the image of the bidisc $\Bbb
D^2\subset\Bbb C^2$ under the mapping with coordinate components
the two elementary symmetric functions of two complex variables.
This domain appears in the spectral Nevanlinna--Pick problem,
related to questions from control theory and applications in
engineering mathematics (see e.g. \cite{Agl-You1,Agl-You4,HMY} and
the references therein). J. Agler and  N. Young \cite{Agl-You3}
showed that $\G_2\in\mathcal L,$ by calculating $l_{\Bbb G_2}.$ On
the other hand, C. Costara \cite{Cos2} proved that $\G_2$ is not
biholomorphic to a convex domain. In addition, let $\mathcal E$
denote the class of domains that can be exhausted with domains,
biholomorphic to convex domains. The Lempert theorem implies that
$\mathcal E\subset\mathcal L.$ A. Edigarian \cite{Edi} showed that
even $\Bbb G_2\not\in\mathcal E$ (see Proposition \ref{b.pr4}). In
connection with this and the mentioned result of D. Jacquet let us
note the following
\smallskip

\noindent{\bf\cite[\bf Problem 2]{Zna}, a hypothesis of L. A.
Aizenberg \cite{Aiz}} {\it Can each $\Bbb C$-convex domain be
exhausted by smooth $\Bbb C$-convex domains?}
\smallskip

Theorem \ref{c.th5} (i) says that $\Bbb G_2$ is $\Bbb C$-convex
domain. Moreover, P. Pflug and W. Zwonek \cite{Pfl-Zwo} have
recently shown that $G_2$ can be exhausted by $\Bbb C$-convex
domains with real-analytic boundaries (see \ref{c.th51}). This
gives an alternative proof of the fact that $\G_2\in\mathcal L.$
We may also formulate the following weaker version of the above
hypothesis,
\smallskip

\noindent{\bf\cite[\bf Problem 4']{Zna}} {\it Does each bounded
$\Bbb C$-convex domain belong to $\mathcal L$?}
\smallskip

An affirmative answer of Problem 4' would follow from an
affirmative answer of \smallskip

\noindent{\bf\cite[\bf Problem 4]{Zna}} {\it Is each bounded $\Bbb
C$-convex domain biholomorphic to a convex domain?}
\smallskip

Theorem \ref{c.th5} (i) together with the result of A. Edigarian
gives a negative answer to the last question.

In a similar way for $\G_2$ one can define the symmetrized
polydisc $\Bbb G_n\subset\Bbb C^n.$ It is natural to ask whether
$\Bbb G_n$ for $n\ge 3$ has the same properties as $\Bbb G_2.$ For
example M. Jarnicki and P. Pflug pose the following question:
\smallskip

\noindent{\bf\cite[\bf Problem 1.2]{Jar-Pfl3}} {\it Does $\Bbb
G_n\in\mathcal L$ or even $\Bbb G_n\in\mathcal E$?}
\smallskip

Clearly if $\Bbb G_n\not\in\mathcal L$ then $\Bbb G_n\not\in\mathcal E.$

Chronologically, first the author proved that $\Bbb
G_n\not\in\mathcal E$ for $n\ge 3$ as Theorem \ref{b.th6}. In its
proof the approach from \cite{Cos2,Edi} is applied on the
so-called generalized balanced domains. In Theorem \ref{b.th3} we
prove that if such a domain in $\Bbb C^n$ belongs to $\mathcal E,$
then its intersection with a special linear subspace of $\Bbb C^n$
is necessarily convex. This is in accordance with the fact that a
(usual) balanced domain is in the class $\mathcal E$ exactly when
it is convex (Corollary \ref{cor}). Theorem \ref{b.th6} follows
from Theorem \ref{b.th3} by showing that the corresponding
intersections for $\Bbb G_n$ are not convex.

Let us note that $\G_n$ for $n\ge 3$ is linearly convex domain,
but it is not $\Bbb C$-convex by Theorem \ref{c.th5} (ii).

Theorem \ref{b.th6} is also a direct corollary from $\Bbb
G_n\not\in\mathcal L,$ $n\ge 3.$ This question is discussed in
Section \ref{ine}. To this aim one uses the infinitesimal forms of
$c_{\Bbb G_n},$ $l_{\Bbb G_n}$ and $k_{\Bbb G_n},$ namely the
Carath\'eodory, Koba\-yashi and Kobayashi--Buseman metrics:
$\gamma_{\Bbb G_n},$ $\kappa_{\Bbb G_n}$ and $\hat\kappa_{\Bbb
G_n}.$ We also introduce a naturally emerging distance $m_{\Bbb
G_n}$ on $\Bbb G_n$ (an analogue to the M\"obius distance $m_{\Bbb
D}$) and its infinitesimal form at the origin, $\rho_n.$

In \cite{Agl-You3} J. Agler and N. Young have shown that
$$l_{\Bbb G_2}=k^\ast_{\Bbb G_2}=c^\ast_{\Bbb G_2}=m_{\Bbb G_2},$$
and $m_{\Bbb G_2}$ is (almost) explicitly calculated. The proof is
based on the method of complex geodesics; their complete
description for $\G_2$ can be found in the work of P. Pflug and W.
Zwonek \cite{Pfl-Zwo}.

In \cite{Cos3} this identity is also obtained for some special
pairs of points from $\Bbb G_n,$ $n\ge 3.$ However in this case it
turns out that $$l_{\Bbb G_n}(0,\cdot)\gneq k^\ast_{\Bbb
G_n}(0,\cdot)\ge c^\ast_{\Bbb G_n}(0,\cdot)\gneq m_{\Bbb
G_n}(0,\cdot)\quad\hbox{(Corollary \ref{i.cor2})}.$$ These
inequalities are directly obtained from the corresponding
inequalities between the infinitesimal forms, which are basically
considered on the coordinate directions (Theorem \ref{i.th1}).

In Proposition \ref{e.pr2} we get an estimate for the difference
between $\gamma_{\Bbb G_{2n+1}}(0;\cdot)$ and $\rho_{2n+1}$ in the
first direction where they do not coincide (Proposition
\ref{e.pr2}). This estimate is based on a "polynomial"\,
description of $\gamma_{\Bbb G_n}.$ Using this description and
computer calculations it is shown that
$$\hat\kappa_{\Bbb G_3}(0;\cdot)\neq\gamma_{\Bbb G_3}(0;\cdot)\quad\hbox{(Theorem \ref{e.th3})}$$
so the Carath\'eodory and Kobayashi metrics do not coincide on
$\Bbb G_3.$ It can be expected that the approach in the proof is
applicable in the higher dimensions too.

The fact that $\G_n$ for $n\ge 3$ has quite different properties
from $\G_2$ is confirmed by Theorem \ref{l.th2}, which gives an
affirmative answer to the following question of M. Jarnicki and P.
Pflug.
\smallskip

\noindent{\bf\cite[\bf Problem 3.2]{Jar-Pfl3}} {\it
Does (unlike $\G_2$) the Bergman kernel of $\G_n$ have zeroes?}
\smallskip

The proof is based on an explicit formula obtained by A. Edigarian
and W. Zwonek in \cite{Edi-Zwo}.

As we noted, the symmetrized polydisc appears in connection with
the spectral Nevanlinna--Pick problem, i.e. an interpolation
problem from the unit disc $\Bbb D$ in the spectral ball
$\Omega_n$ -- the set of complex $n\times n$ matrices of spectral
radius less than 1 (i.e. with eigenvalues in $\D$). The
infinitesimal form of this problem is the spectral
Carath\'eodory-Fej\'er problem. The easiest forms of these
problems are reduced to finding $l_{\Om_n}$ and $\kappa_{\Om_n},$
while the continuous dependance on the given data -- to the
continuity of these two functions. In the case of cyclic matrices
(i.e. ones with a cyclic vector) they coincide with the
corresponding functions on the taut domain $\Bbb G_n$, so they are
continuous.

In Section \ref{cyc} we provide some equivalent conditions for a
matrix to be cyclic (part of these are used in the last sections
of the chapter).

In Section \ref{np} we gather the basic properties of the above
problems and their reduction to similar problems on the
symmetrized polydisc in the case of cyclic matrices. (As this is a
taut domain, the problems there "depend"\ on the data in a
continuous manner). This also determines the corresponding
relationships with the Lempert function and the Kobayashi metric
on the symmetrized polydisc.

Section \ref{cont} is dedicated to the continuity of $l_{\Om_n}$
(in the general case). The main result there (Theorem \ref{septh})
states that $l_{\Om_n}(A;\cdot)$ is a continuous function exactly
when $A$ is a scalar matrix or $n=2$ and $A$ has (two) equal
eigenvalues. This result is based on Proposition \ref{np.pr4},
which is obtained from the basic Proposition \ref{i.pr4} (iii).

In Section \ref{zeko} we discuss the (dis)continuity of
$\kappa_{\Om_n}$ by studying its zeroes. In particular we have
found all matrices $A\in\Om_3$ such that $\kappa_{\Om_3}(A;B)>0$
for $B\neq 0$ (a relatively easy question for $n=2$).

As an application, in the last section \ref{k-l} we have shown
that the Kobayashi metric of a pseudoconvex domain is not equal to
the weak "derivative"\ of the Lempert function in the general case
(this gives a partially affirmative answer to a question from
Section \ref{lem-der}). The counterexample is $\Omega_3$ (or, of a
lower dimension, the domain of zero-trace matrices in $\Omega_3$).

Finally, we point out that a detail study of another biholomorphic
invariant -- the pluricomplex Green function -- on the the
spectral ball and the symmetrized polydisc can be found in
\cite{TTZ}.
\vspace{2mm}

{\bf Chapter \ref{chap.bound}.} The main purpose of this chapter
is obtaining estimates (in a geometric way) for the
Carath\'eodory, Kobayashi and Bergman metrics, as well as for the
Bergman kernel (on the diagonal), of an arbitrary $\C$-convex
domain $D\subset\C^n$ not containing complex lines, in terms of
the distance $d_D(z;X)$ from the point $z\in D$ to the boundary
$\partial D$ in the direction $X\in(\C^n)_\ast.$ These estimates
show that on such a domain these three metrics coincide up to a
constant, depending only on $n$ (Corollary \ref{ber2}). Similar
results in the special case of a $\mathcal C^\infty$-smooth
bounded $\C$-convex domain of a finite type, with quite hard
proofs, are the main results of the dissertations of S. Blumberg
\cite{Blu} and M. Lieder \cite{Lie}. In addition, the constants
there depend on the domain. Earlier similar results for convex
domains can be found in the Ph.~D.~thesis \cite{Chen} of J.-H.
Chen and in the works \cite{McN2,McN3} of J. D. McNeal; however
their proofs are experiencing essential weaknesses.

Using the $1/4$-Theorem of K\"obe, it is easily shown in Proposition
\ref{bound1} that \vspace{-1.5pt}
$$1/4\le\gamma_D(z;X)d_D(z;X)\le\kappa_D(z;X)d_D(z;X)\le 1.$$
The two (absolute) constants are exact, and $1/4$ can be replaced by $1/2$
for convex domains (see \cite{BP} or Proposition
\ref{bound0}).

As an application of these estimates, in Section \ref{type} we get
that the standard and linear multitypes of D'Angelo and Catlin
coincide for a smooth boundary point of a $\C$-convex domain. This
generalizes a result of M. Conrad \cite{Con} and J. Yu \cite{Yu1}
(see also the works of J. D. McNeal \cite{McN1}, and H. P. Boas
and E. Straube \cite{Boas-Str}), while our proof is essentially
different and much shorter. It is based upon a easy result from
\cite{Yu3} and moving the main result in \cite{Lee} from the
convex to the $\C$-convex case (the considerations here are easier
than in \cite{Lee}).

The main result of Chapter \ref{chap.bound} Theorem \ref{ber1}
states that there is an inequality for the Bergman metric, similar
to the above one; the corresponding constants depend only on the
dimension $n$ of the domain. To this aim in Theorem \ref{bound3}
we get some estimates for the Bergman kernel, which are also of
independent interest. The constants there depend only on $n$ and
are exact for the class of convex domains. These estimates are
connected with the so called minimal basis (for a point in a given
domain), introduced by T. Hefer \cite{Hef1} for the smooth case
(of finite type) and somewhat later, but independent, by the
author and P. Pflug \cite{Nik-Pfl3} in the general case. It is
used in the proof of Theorem \ref{ber1} and almost all arguments
are geometrical. One can define a minimal basis for a point in a
given open set (not containing complex lines) by induction: the
first vector of the base is directed towards the closest point
from the set boundary, and the next ones are from the basis of the
intersection of the set with the complex hyperplane through the
point, orthogonal to that vector. The main (and trivial) property
of that basis that is used for weakly linearly convex domains is
the orthogonality of the intersections of complex "support"\
hyperplanes through the emerging boundary points and the
corresponding vectors form the basis. The geometrical arguments
are completed by the stability of $\C$-convexity under
projections.

In the previously mentioned works \cite{Chen,McN2,McN3}, apart
from $\overline\partial$-technique, the authors use a similar
(however notably more complicated) method in the terms of a basis
that we will call maximal. In Section \ref{count} we provide a
natural counterexample for the main "property"\ of this basis (the
same as for the minimal one) that is used in these and other works
for various problems (e.g. for the linear and D'Angelo types in
the already quoted paper \cite{McN1}). Anyway in Section \ref{max}
we show how the estimates obtained in the minimal basis imply
those for the maximal one (using some combinatorial arguments).

Another aim of this chapter is to state the local character of the
obtained results by showing that the estimates near a given
boundary point $a$ of a domain remain true if the domain is weakly
locally linearly convex near $a$ and the boundary does not contain
analytic discs through $a.$ Such a domain with a $\mathcal
C^2$-smooth boundary near $a$ turns out to be locally $\C$-convex
(Proposition \ref{loc-con}). Then the local character of the
estimates for the Kobayashi metric (if the domain is bounded)
follows from the general localization Proposition \ref{loc-kob}.
Its proof permits to obtain immediately the exact boundary
behavior of this metric near an isolated point of a planar domain,
having at least one more point in its boundary. This essentially
strengthens a main result from \cite{KLZ}.

The local character of the estimates for the Bergman kernel and
Bergman metric is determined in Section \ref{loc-ber}, where the
domain is assumed to be pseudoconvex (but not necessarily bounded)
and locally convex around a boundary point not contained in
analytic (or, equiva\-lently, linear) discs from the boundary. The
proof is based on the existence of a locally holomorphic peak
function at this point (Proposition \ref{con-peak}) and the
localization theorem \ref{peak} for the Bergman kernel and Bergman
metric (if such a function exists). In the case of a bounded
pseudoconvex domain this theorem is contained in the fundamental
work \cite{Hor1} of L. H\"ormander as an application of the
obtained $L^2$-estimates for the $\overline\partial$-problem. Our
proof is a variation of this technique. As a corollary we get a
stronger variant of the main result of G. Herbort in \cite{Her1}
without a usage of $\overline\partial$-technique of
Ohsawa--Takegoshi (see the remark at the end of Section
\ref{loc-ber}). The proof also implies weak localization of the
Bergman kernel and Bergman metric for a planar domain with a
non-polar complement (Corollary \ref{antipeak}).

The last section we get the exact boundary behavior of the
invariant metrics under consideration near a $\mathcal C^1$-smooth
boundary point of an arbitrary planar domain, once again using a
geometric argument (the Pinchuk scaling method).
\vspace{2mm}

This work is the author's D.~Sc.~dissertation originally written
in Bulgarian and defended in October, 2010.

The results are published as follows:

Chapter \ref{chap.lem} -- in \cite{Nik-Pfl5,Nik-Pfl6,Nik-Pfl7,Nik-Pfl8,Nik-Pfl9,Nik-Zwo1};

Chapter \ref{chap.spec} -- in
\cite{Nik3,Nik-Pfl10,NPTZ1,NPTZ2,NPZ1,NPZ2,Nik-Tho1,Nik-Tho2,NTZ,Nik-Zwo2};

Chapter \ref{chap.bound} -- in
\cite{Jar-Nik,Nik1,Nik2,Nik-Pfl2,NPT1,Nik-Pfl3,NPZ3}.

Some of the mentioned results are from
\cite{Nik0,Nik-Pfl4,NPT2,Nik-Sar}.
\vspace{2mm}

\noindent{\bf Acknowledgments.} The author like to thank Peter
Pflug, Pascal J.~Thomas, W\l odzimierz Zwonek and Marek Jarnicki
for the fruitful collaboration. He is also grateful to P.~Pflug
and M.~Jarnicki for their remarks which improved the manuscript.
The author is indebted to the Carl von Ossietzky University
(Oldenburg), the Paul Sabatier University (Toulouse) and the
Jagiellonian University (Krak\'ow) for the hospitality during his
stays there, as well as to DFG, DAAD and CNRS for their supports.

\chapter{Lempert functions and Kobayashi metrics}\label{chap.lem}
\setcounter{equation}{0}
\section{Lempert functions and their "derivatives"}\label{lem-der}

In this section we introduce the Lempert functions of higher order
and their infinitesimal forms, the Kobayashi metrics of higher
order, for an arbitrary complex manifold. (see also
\cite{Jar-Pfl1,KobS2}).

Our main aim is to prove that if the Kobayashi metric of a complex
manifold is continuous and positive at a given point for each
nonzero tangent vector, then the "derivatives"\ of the Lempert
functions exist and are equal to the corresponding Kobayashi
metrics at this point. This generalizes some results of M.-Y. Pang
\cite{Pang} and M. Kobayashi \cite{KobM} for taut
domains/manifolds.

As usual $\Bbb D\subset\Bbb C$ denotes the unit disc. Let $M$ be
an $n$-dimensional complex manifold. Let us recall the definitions
of the Lempert function $l_M$ and the Kobayashi--Royden (in short,
Kobayashi) (pseudo)metric $\kappa_M$ of $M$:
$$\aligned
l_M(z,w)&=\inf\{|\alpha|:\exists f\in\mathcal
O(\Bbb D,M):f(0)=z,f(\alpha)=w\},\\
\kappa_M(z;X)&=\inf\{|\alpha|:\exists f\in\mathcal O(\Bbb D,M):
f(0)=z,\alpha f_{\ast,0}(d/d\zeta)=X\},
\endaligned$$
where $X$ is a complex tangent vector to $M$ at $z$. Such $f$
always exist (see e.g. \cite{Win}; also according to \cite[p.
49]{Din} this has been known even earlier by J. Globevnik).

Note that if $F:M\to N$ is a holomorphic mapping between two
manifolds, then $$l_M(z,w)\ge l_N(F(z),F(w)).$$ In particular, if
$F$ is a biholomorphism, then we get an equality, i.e. the Lempert
function is invariant under biholomorphisms. The above inequality
also shows that this function is the largest holomorphically
contractible function that coincides on $\D$ with the M\"obius
distance $m_\D.$ On the other hand, the smallest such function is
the Carath\'eodory function:
$$c^\ast_M(z,w)=\sup\{m_\D(f(z),f(w)):f\in\mathcal O(M,\D)\}.$$
If in this definition we replace $m_\D$ by the Poincar\'e distance $p_\D,$
we get the Carath\'eodory (pseudo)distance
$$c_M=\tanh^{-1}c^\ast_M.$$
As $$\kappa_M(z;X)\ge\kappa_N(F(z),F_{\ast,z}(X)),$$ the Kobayashi
metric is the largest holomorphically contractible pseudometric
such that $\kappa_\D(0;X)=|X|.$ The smallest such pseudometric is
the Carath\'eodory--Reiffen (ni short, Carath\'eodory) metric
$$\gamma_M(z;X)=\sup\{|f_{\ast,z}(X)|:f\in\mathcal O(M,\D)\}$$
(we can assume $f(z)=0$).

As in the case of domains, the Kobayashi distance $k_M$ can be
defined as the largest pseudodistance not exceeding the Lempert
function of first order $$k_M^{(1)}=\tanh^{-1}l_M$$ (for
convenience we distinguish this function from the Lempert function
$l_M$). By the Kobayashi function we mean $$k^\ast_M=\tanh k_M.$$

Let us note that if $k_M^{(m)}$ denotes the Lempert function of
order $m$ ($m\in\Bbb N$), i.e.
$$k_M^{(m)}(z,w)=\inf\{\sum_{j=1}^m k_M^{(1)}(z_{j-1},z_{j}):
z_0,\dots,z_m\in M,z_0=z,z_m=w\},$$ then
$$k_M(z,w)=k_M^{(\infty)}:=\inf_m k_M^{(m)}(z,w).$$
Now let us recall that a manifold $M$ is called a taut manifold if
the family $\O(\D,M)$ is normal. Every taut domain in $\C^n$ is
pseudoconvex. Conversely, every bounded domain with a $\mathcal
C^1$-smooth boundary is hyperconvex (i.e. has an exhausting
negative plurisubharmonic function), so it is a taut domain.

According to a result of M.-Y. Pang \cite{Pang}, the Kobayashi
metric is the "derivative"\ of the Lempert function, if the
corresponding domain is a taut domain: $$\kappa_D(z;X)=\lim_{t\to
0}\frac{l_D(z,z+tX)}{t}$$ (in this limit, as well as in some
similar ones below, we can replace $l_D$ by $k^{(1)}_D$ and, in
general, an invariant function with values in $[0,1)$, by
$\tanh^{-1}$ of it, or vice versa).

In the general case the Kobayashi metric at a given point of a
domain is not a pseudonorm (vectorwise), i.e. its indicatrices are
not convex domains. To avoid this defect, S. Kobayashi
\cite{KobS1} introduces a new invariant metric, later called
Kobayashi--Buseman metric. As in the case with the Kobayashi
distance, this metric $\hat\kappa_M$ can be defined by putting
$\hat\kappa_M(z;\cdot)$ to be the largest pseudonorm not exceeding
$\kappa_M(z;\cdot)$. Clearly
$$
\hat\kappa_M(z;X)=\inf\{\sum_{j=1}^m\kappa_M(z;X_j):m\in\Bbb N,\
\sum_{j=1}^mX_j=X\}.
$$ Hence it is natural to consider the functions $\kappa_M^{(m)}$,
$m\in\Bbb N$ defined as follows:
$$
\kappa_M^{(m)}(z;X)=\inf\{\sum_{j=1}^m\kappa_M(z;X_j):
\sum_{j=1}^mX_j=X\}.
$$
We call the function $\kappa_M^{(m)}$ a Kobayashi metric of order
$m.$ Clearly $\kappa_M^{(m)}\ge\kappa_M^{(m+1)}.$ Also one can
easily observe that if
$\kappa_M^{(m)}(z;\cdot)=\kappa_M^{(m+1)}(z;\cdot)$ for some $m$,
then $\kappa_M^{(m)}(z;\cdot)=\kappa_D^{(j)}(z;\cdot)$ for each
$j>m$. Furthermore, as we will see in the next section,
$\kappa_M^{(2n-1)}=\kappa_M^{(\infty)}:=\hat\kappa_M,$ with $2n-1$
being the least possible number in the general case.

Let us note that all introduced objects are upper semicontinuous;
for $\kappa_M$ (and hence for $\kappa_M^{(m)}$ and $\hat\kappa_M$) see also
\cite{KobS2}. To convince ourselves for $k_M^{(m)},$ it suffices to check it
for $l_M$.

\begin{proposition}\label{der.pr0} For each complex manifold $M$, the function
$l_M$ is upper semicontinuous.
\end{proposition}

\beginproof We use a standard procedure (see \cite{Roy}). Let $r\in(0,1)$ and $z,w\in M$.
Let $f\in\mathcal O(\Bbb D,M)$, $f(0)=z$ and
$f(\alpha)=w$. Then $\tilde f=(f,\hbox{id}):\Delta\to\tilde M=M\times\Delta$ is
an immersion.
Put $\tilde f_r(\zeta)=\tilde
f(r\zeta)$; now \cite[Lemma 3]{Roy} implies that there is a Stein neighborhood
$S\subset\tilde M$ of $\tilde f_r(\Bbb D)$. As well known,
$S$ can be immersed as a closed complex manifold in $\Bbb C^{2n+1}.$ Let
$\psi$ be the corresponding immersion.
Then there is an (open) neighborhood $V\subset\Bbb C^{2n+1}N$ of $\psi(S)$ and a
holomorphic retraction $\theta:V\to\psi(S).$ For $z'$
near $z$ and $w'$ near $w$ we can find (in a standard way) $g\in\mathcal O(\Bbb D,V)$,
so that $g(0)=\psi(z',0)$ and
$g(\alpha/r)=\psi(w',\alpha)$. Denote by $\pi$ the natural projection of $\tilde M$
onto $M.$ Then
$h=\pi\circ\psi^{-1}\circ\theta\circ g\in\mathcal O(\Bbb D,M)$,
$h(0)=z'$ and $h(\alpha/r)=w'$. Consequently
$rl_M(z',w')\le\alpha.$ This shows that $\limsup_{z'\to z,w'\to
w}l_M(z',w')\le l_M(z,w).$\qed
\smallskip

To extend the previously mentioned result of Pang, we define the
"derivatives"\ of $k^{(m)}_M$, $m\in\Bbb N^\ast=\Bbb
N\cup\{\infty\}.$ Let $(U,\varphi)$ be a holomorphic chart near
$z.$ We put
$$\mathcal D k^{(m)}_M(z;X)=\limsup_{t\rightarrow 0,w\to
z,Y\to\varphi_\ast X}
\frac{k^{(m)}_M(w,\varphi^{-1}(\varphi(w)+tY))}{|t|}.$$ This
definition does not depend of the chart; also,
$$\mathcal D k^{(m)}_M(z;\lambda X)=|\lambda|\mathcal D
k^{(m)}_M(z;X),\quad \lambda\in\Bbb C.$$ Replacing $\limsup$ by
$\liminf$, we can define $\underline{\mathcal D} k_M^{(m)}$.

A result of M. Kobayashi \cite{KobM} shows that if $M$ is a
complex taut manifold, then
$$
\hat\kappa_M(z;X)=\mathcal D k_M(z;X)=\underline{\mathcal D}
k_M(z;X),$$ i.e. the Kobayashi--Buseman metric is the
"derivative"\ of the Kobayashi distance. The proof of this result
allows us to learn something more: $$\kappa_M^{(m)}(z;X)=\mathcal
D k^{(m)}_M(z;X)=\underline{\mathcal D} k_M^{(m)}(z;X),\quad
m\in\Bbb N^.$$

Let us not that for the Carath\'eodory metric of an arbitrary complex manifold $Ì$
one has (see \cite{Jar-Pfl1} for domains in $\C^n$)
\begin{equation}\label{car}
\gamma_M=\mathcal D c_M=\underline{\mathcal D} c_M
\end{equation}
(the definitions of the last two invariants are obvious).

To formulate in full generality the main result of this section we need the following
notion. A complex manifold is called
hyperbolic at the point $z\in M,$ if $k_M(z,w)>0$ for each $w\neq z.$ (Recall that
$M$ is hyperbolic if it is hyperbolic at each of its points,
i.e. $k_M$ is a distance). Then the following assertions are equivalent:

(i) $M$ is hyperbolic at $z;$

(ii) $\liminf_{z'\to z,w\in M\setminus U}l_M(z',w)>0$ for each neighborhood $U$ of $z;$

(iii) $\underline{\kappa}_M(z;X):=\liminf_{z'\to z,X'\to X}\kappa_M(z';X')>0$ for each
$X\neq 0;$

The implications (i)$\Rightarrow$(ii)$\Rightarrow$(iii) are (almost) trivial, while
(iii)$\Rightarrow$(i) follows from the fact that
$k_M$ is the integrated form of $\kappa_M.$

In particular, if $M$ is hyperbolic at $z$, then it is hyperbolic at each point $z'$
near $z$.

If $M$ is a taut manifold, then it is hyperbolic and $\kappa_M$ is
a continuous function. This shows that the theorem below
generalizes the previously mentioned result of M. Kobayashi.

\begin{theorem}\label{der.th1} Let $M$ be a complex manifold and $z\in M$.

(i) If $M$ is hyperbolic at $z$ and $\kappa_M$ is continuous at $(z,X)$, then
$$\kappa_M(z;X)=\mathcal Dl_M(z;X)=\underline{\mathcal D} l_M(z;X).$$

(ii) If $\kappa_M$ is continuous and positive at $(z,X)$ for each $X\neq 0,$ then
$$\kappa_M^{(m)}(z;\cdot)=\mathcal D k_M^{(m)}(z;\cdot)=\underline{\mathcal D}
k_M^{(m)}(z;\cdot),\quad m\in\Bbb N^\ast.$$
\end{theorem}

The first step of the proof of Theorem \ref{der.th1} is the following

\begin{proposition}\label{der.pr2} For each complex manifold $M$ one has
$$\kappa_M^{(m)}\ge\mathcal D k_M^{(m)}, \quad m\in\Bbb N^\ast.$$
\end{proposition}

Note that if $M$ is a domain, a weaker variant of Proposition
\ref{der.pr2} can be found in \cite{Jar-Pfl1}, namely
$\hat\kappa_M\ge\mathcal D k_M$ (the proof is based on the fact
that $\mathcal D k_M(z;\cdot)$ is a pseudonorm).
\smallskip

\noindent{\it Proof of Proposition \ref{der.pr2}.} Let us first
consider the case $m=1.$ A main role will be played by the
following

\begin{theorem}\label{der.th3}\cite{Roy}*\footnote{*Instead of Theorem
\ref{der.th3},
one can use the approach from the proof of the semicontinuity of
$l_M$.} Let $M$ be a complex manifold and the mapping
$f\in\mathcal O(\Bbb D,M)$ be regular at $0$. Let $r\in(0,1)$ and
$D_r=r\Bbb D\times\Bbb D^{n-1}$. Then there is a mapping
$F\in\mathcal O(D_r,M)$ that is singular at $0$ and $F|_{r\Bbb
D\times\{0\}}=f$.
\end{theorem}

Since $\kappa_M(z;0)=\mathcal Dl_M(z;0)=0$, one can assume that
$X\neq 0$. Let $\alpha>0$ and $f\in\mathcal O(\Bbb D,M)$ be such
that $f(0)=z$ and $\alpha f_{\ast,0}(d/d\zeta)=X$. Let $r\in(0,1)$
and $F$ be as in Theorem \ref{der.th3}. Since $F$ is regular at
$0,$ there are neighborhoods $U=U(z)\subset M$ and $V=V(0)\subset
D_r$, such that $F|_V:V\to U$ is a biholomorphism. Therefore
$(U,\varphi),$ where $\varphi=(F|_V)^{-1}$ is a map near $z$. Let
us note that $\varphi_{\ast,z}(X)=\alpha e_1$, where
$e_1=(1,0,\dots,0).$

If $w$ and $Y$ are close enough to $z$ and $\alpha e_1$, then
$g(\zeta)=F(\varphi(w)+\zeta Y/\alpha)$ belongs to $\mathcal
O(r^2\Bbb D,M),$ $g(0)=w$ and
$g(t\alpha)=\varphi^{-1}(\varphi(w)+tY),$ $t<r^2/\alpha.$
Consequently $r^2l_M(w,\varphi^{-1}(\varphi(w)+tY))\leq t\alpha$.
Thus $r^2\l_M(z;X)\le\alpha$. For $r\to 1$ and
$\alpha\to\kappa_M(z;X)$ we get $\mathcal
Dl_M(z;X)\le\kappa_M(z;X).$

Now let $m\in\Bbb N$. Recall that $\kappa_M^{(m)}(z;\cdot)$ is the
largest function with the following property:

For each $X=\sum_{j=1}^mX_j$ it follows that
$\kappa_M^{(m)}(z;X)\le\sum_{j=1}^m\kappa_M(z;X_j)$.

To prove that $\kappa_M^{(m)}\ge\mathcal D k_M^{(m)},$ it is
sufficient to check that $\mathcal D k^{(m)}_M(z;\cdot)$ has this
property Using the above notation and choosing
$Y_j\to\varphi_{\ast,z}X_j$ so that $\sum_{j=1}^m Y_j=Y$, we put
$w_0=w$ and $w_j=\varphi^{-1}(\varphi(w)+t\sum_{k=1}^{j}Y_j)$.
Since $$k^{(m)}_M(w,w_q)\le\sum_{j=1}^mk^{(1)}_M(w_{j-1},w_j),$$
from the case $m=1$ it follows that $$\mathcal D
k^{(m)}_M(z;X)\le\sum_{j=1}^m\mathcal D
k_M(z;X_j)\le\sum_{j=1}^m\kappa_M(z;X_j).$$

Finally, let $m=\infty$ and $n=\dim M$. Since
$\hat\kappa_M=\kappa_M^{(2n-1)}$ and $k_M\le k_M^{(2n-1)}$, the
case $m=2n-1$ shows that $\mathcal D k_M\le\hat\kappa_M.$\qed
\smallskip

\noindent{\it Proof of Theorem \ref{der.th1}.} We can assume that
$X\neq 0$. Having in mind Proposition \ref{der.pr2}, we just have
to prove that $$\kappa_M^{(m)}(z;X)\le\underline{\mathcal
D}k^{(m)}_M(z;X)$$ under the corresponding assumptions. For
simplicity we assume that $M$ is a domain in $\Bbb C^n$ (the
changes in the general case of a manifold are obvious).

(i) Fix a neighborhood $U=U(z)\Subset M.$ By hyperbolicity of $M$
at $z$, there exist a neighborhood $V=V(z)\subset U$ and a number
$\delta\in (0,1)$ such that if $h\in\mathcal O(\Bbb D,M)$ and
$h(0)\in V$, then $h(\delta\Bbb D)\subset U$. From here by the
Cauchy inequalities it follows that $||h^{(k)}(0)||\leq
c/\delta^k$, $k\in\Bbb N$ ($||\cdot||$ is the Euclidean norm).

Now choose sequences $w_j\to z$, $t_j\to 0$ and $Y_j\to X,$ such that
$$
\frac{l_M(w_j,w_j+t_jY_j)}{|t_j|}\to\underline{\mathcal
D}l_M(z;X).
$$
Let the holomorphic discs $g_j\in\mathcal O(\Bbb D,M)$ and the
numbers $\beta_j\in (0,1)$ be such that $g_j(0)=w_j$,
$g_j(\beta_j)=w_j+t_jY_j$ and $\beta_j\le
l_M(w_j,w_j+t_jY_j)+|t_j|/j.$ Note that $l_M(w_j,w_j+t_jY_j)\le
c_1||t_jY_j||\le c_2 |t_j|.$

Let
$$
w_j+t_jY_j=g_j(\beta_j)=w_j+g_j'(0)\beta_j+ h_j(\beta_j)
$$
Then $$||h_j(\beta_j)||\leq c\sum_{k=2}^\infty
\left(\tfrac{\beta_j}{\delta}\right)^k\leq c_3|\beta_j|^2\le
c_4|t_j|^2,\quad j\ge j_0.$$

We put $\hat Y_j=Y_j-h_j(\beta_j)/t_j.$ Then
$g_j(0)=w_j$ and $\beta_jg_j'(0)/t_j=\hat Y_j\to X$.
Consequently,
$$\kappa_M(w_j;\hat
Y_j)\leq\frac{\beta_j}{|t_j|}\leq\frac{l_M(z_j,w_j+t_jY_j)}{|t_j|}+\frac{1}{j}.
$$
For $j\to\infty$ we get
$\kappa_M(z;X)=\underline{\kappa}_M(z;X)\leq\underline{\mathcal
D}l_M(z;X)$.

(ii) The proof of the case $m\in\Bbb N$ is similar to the one below and we omit it.
Now let $m=\infty$.

Our assumptions show that $M$ is hyperbolic at $z.$ Also it easily follows
(say by contradiction) that

\begin{equation}\label{der1}
\forall\varepsilon>0\
\exists\delta>0:||w-z||<\delta,||Y-X||<\delta||X|| \Rightarrow
\end{equation}
$$|\kappa_M(w;Y)-\kappa_M(z;X)|<\varepsilon\kappa_M(z;X).$$
Also, the proof of (i) shows that
\begin{equation}\label{der2}
k^{(1)}_M(a,b)\ge\kappa_M(a;b-a+o(a,b)),\hbox{ where
}\lim_{a,b\to z}\frac{o(a,b)}{||a-b||}=0.
\end{equation}

Now choose sequences $w_j\to z$, $t_j\to 0$ and $Y_j\to X$ such that
$$\frac{k_M(w_j,w_j+t_jY_j)}{|t_j|}\to\underline{\mathcal
D}k_M(z;X).$$ Let $w_{j,0}=w_j,\dots,w_{j,m_j}=w_j+t_jX_j$ be
points from $M,$ such that
\begin{equation}\label{der3}
\sum_{k=1}^{m_j}k^{(1)}_M(w_{j,k-1},w_{j,k})\le
k_M(w_j,w_j+t_jY_j)+\frac{1}{j}.
\end{equation}
Put $w_{j,k}=w_j$ for $k>m_j$. Since
$$k_M(w_j,w_{j,l})\le\sum_{j=1}^lk^{(1)}_M(w_{j,k-1},w_{j,k})\le
k_M(w_j,w_j+t_jY_j)+\frac{1}{j}\le c_2|t_j|+\frac{1}{j},$$
$k_M(w_j,w_{j,l})\to 0$ uniformly on $l.$ The hyperbolicity of $M$
at $z$ implies that $w_{j,l}\to z$ uniformly on $l.$ Indeed,
assuming the contrary and choosing a subsequence, we can consider
that $w_{j,l_j}\not\in U$ for some $U=U(z).$ Then
$$0=\lim_{j\to\infty}k_M(w_j,w_{j,l})\ge\liminf_{z'\to z,w\in
M\setminus U}l_M(z',w)>0,$$ which is a contradiction.

Finally let us fix $R>1.$ Then (\ref{der1}) shows that
$$\kappa_M(z;w_{j,k}-w_{j,k-1})\le R\kappa_M(w_{j,k};w_{j,k}-w_{j,k-1}
+o(w_{j,k},w_{j,k-1})),\ j\ge j(R).$$ From this inequality, (\ref{der2})
and (\ref{der3}) it follows that
$$\sum_{k=1}^{m_j}\kappa_M(z;w_{j,k}-w_{j,k-1})
\le Rk_M(w_j,w_j+t_jYj)+\frac{R}{j}.$$ Since $\hat\kappa_M(z;t_jY_j)$
is bounded by the above sum, we get that
$$\hat\kappa_M(z;Y_j)\le R\frac{k_M(w_j,w_j+t_jYj)+1/j}{|t_j|}.$$
It remains to use that $\hat\kappa_M(z;\cdot)$ is a continuous function.
Then for $j\to\infty$ and $R\to 1$ it follows that
$\hat\kappa_M(z;X)\le\underline{\mathcal D}k_M(z;X)$.\qed
\smallskip

\noindent{\bf Remark.} From the above proofs, by a standard diagonal process,
it follows that if $M$ is hyperbolic at $z,$ then
$\underline{\kappa}_M(z;\cdot)=\underline{\mathcal D}l(z;\cdot).$
\smallskip

The subsequent examples show that the assumptions of continuity in Theorem
\ref{der.th1} are essential.

\smallskip

$\bullet$ Let $A$ be a countable dense subset of $\Bbb
C_\ast(=\C\setminus\{0\}).$ In \cite{Die-Sib} (see also
\cite{Jar-Pfl1}) there is an example of a pseudoconvex domain
$D\subset\Bbb C^2$ such that:

(i) $(\Bbb C\times\{0\})\cup(A\times\Bbb C)\subset D;$

(ii) if $z_0=(0,t)\in D,$ $t\neq 0,$ then $\kappa_D(z_0;\cdot)\ge
C||\cdot||$ for some $C>0$. (It can be even shown that $\mathcal
Dl_D(z_0;\cdot)\ge C||\cdot||.$)

\noindent Then it is easily concluded that $\underline{
\kappa}_D(\cdot;e_2)=\mathcal D k^{(3)}_D(\cdot;e_2)=k^{(5)}_D=0$
è $\hat\kappa_D(z_0;\cdot)\ge c||\cdot||,$ where $e_2=(0,1)$ and $c>0.$
Therefore
$$\hat\kappa_D(z_0;X)>0=\underline{ \kappa}_D(z_0;e_2)=\mathcal D
k^{(3)}_D(z_0;e_2)=\mathcal D k^{(5)}_D(z_0;X),\quad X\in(\C^2)_\ast.$$

This phenomenon clearly appears also in $\Bbb C^n$, $n>2$ (say for $D\times\Bbb D^{n-2}$).
Thus the inequalities in
Proposition \ref{der.pr2} are strict in the general case.

\smallskip
$\bullet$ There exists a bounded pseudoconvex domain
$D\subset\C^2$ containing the origin such that (see e.g.
\cite[Example 4.2.10]{Zwo1})
$$\kappa_D(0;e_1)=\mathcal D k_D(0;e_1)=\limsup_{t\to
0}\frac{l_D(0,te_1)}{|t|}$$
$$>\liminf_{t\to 0}\frac{l_D(0,te_1)}{|t|}\ge
\underline{\mathcal D}k_D(0;e_1).$$

We conclude this section by the following
\smallskip

\noindent{\bf Question.} Is $\kappa_D\neq\mathcal Dl_D$ in the general case?
Is $\mathcal D k_D$ a holomorphically contractible invariant?
(For this question see also \cite{Jar-Pfl3}.)
\smallskip

A partially positive answer will be given in Section \ref{k-l} by
showing that there is a pseudoconvex domain $D\subset\C^8$ and a
point $(z,X)\in D\times\C^n$ such that
$$\kappa_D(z;X)>0=\limsup_{t\to 0}\frac{l_D(z,tX)}{|t|}.$$

\setcounter{equation}{0}
\section{Balanced domains}\label{lem-bal}

The biholomorphic invariants can be explicitly calculated for a few classes of
domains, usually contained in the class of Reinhardt domains.
Each complete Reinhardt domain is balanced. In this section we determine some
relationships between the Minkowski functions of a balanced domain
or of its convex/holomorphically convex hull and some biholomorphic invariants
of the domain when one of their arguments is the origin.

Recall that a domain $D\subset\Bbb C^n$ is called balanced if
$\lambda z\in D$ for each $(\lambda,z)\in\overline{\Bbb D}\times
D$ (for this definition and a part of the facts below see e.g.
\cite{Jar-Pfl1}). We naturally associate to such a domain its
Minkowski function
$$h_D(z)=\inf\{t>0:z/t\in D\},\ z\in\Bbb C^n.$$ The function $h_D\ge0$ is
upper semicontinuous and
$$h_D(\lambda z)=|\lambda|h_D(z),\ \lambda\in\Bbb C,z\in\Bbb C^n,$$
$$D=\{z\in\Bbb C^n:h_D(z)<1\}.$$ Let us note that $D$ is pseudoconvex
exactly when $\log h\in\PSH(\Bbb C^n),$ which in this case is
equivalent to $h\in\PSH(\Bbb C^n).$ Also recall that $D$ is a taut
domain exactly when it is bounded and $h_D$ is a continuous
plurisubharmonic function. This shows that, for a balanced domain,
being hyperconvex or taut is the same. Let us note that the
hyperbolicity of $D$ is equivalent to its boundedness. More
general results concerning the so-called Hartogs domains can be
found in the paper \cite{Nik-Pfl4} of the author and P. Pflug.

Clearly, the convex hull $\hat D$ of a balanced domain $D$ is balanced.
Let us recall the well-known relationships between
$h_D,$ $\hat h_D=h_{\hat D}$ and some invariant functions and metrics.

\begin{proposition}\label{bala.pr1}
Let $D\subset\Bbb C^n$ be a balanced domain and $a\in\ D.$ Then:

(i) $\gamma_D(0;\cdot)=\hat\kappa_D(0;\cdot)=\hat h_D.$

(ii) $\hat h_D\le c_D^\ast(0,\cdot)\le k^\ast_D(0,\cdot)\le
l_D(0,\cdot)\le h_D$ and $\hat h_D\le\kappa(0;\cdot)\le h_D;$

(iii) $c^\ast(0,a)=h_D(a)\Leftrightarrow k^\ast_D(0,a)=h_D(a)
\Leftrightarrow h_D(a)=\hat h_D(a).$

\noindent If in addition $D$ is pseudoconvex, then

(iv) $l_D(0,\cdot)=h_D$ and $\kappa(0;\cdot)=h_D.$
\end{proposition}

The Lempert theorem (mentioned in the introduction) implies that
$c^\ast_D=k^\ast_D=l_D$ for each convex domain $D.$ Then by the
above proposition we get

\begin{corollary}\label{bala.cor0} For a pseudoconvex balanced domain
$D\subset\C^n$ the following are equivalent:

(i) $D$ is convex (i.e. $h_D=\hat h_D$);

(ii) $c^\ast_D=l_D;$

(iii) $c^\ast_D(0,\cdot)=l_D(0,\cdot);$

(iv) $k^\ast_D=l_D;$

(v) $k^\ast_D(0,\cdot)=l_D(0,\cdot);$
\end{corollary}

Put $(k^{(m)}_D)^\ast=\tanh k^{(m)}_D.$ Proposition
\ref{bala.pr1} (iii) shows that at the point $a\in D$ the value of
$k_D(0,\cdot)$ is maximal exactly when $D$ is "convex"\
in the direction of $a,$ i.e. $h_D(a)=\hat h_D(a)$. The next result shows
that something more is true.

\begin{proposition}\label{bala.pr2} Let $D\subset\Bbb C^n$ be a balanced domain
and $a\in D.$ The following are equivalent:

(i) $\hat h_D(a)=h_D(a);$

(ii) $(k^{(3)}_D(0,a))^\ast=h_D(a);$

(iii) $\kappa^{(2)}_D(0;a)=h_D(a).$
\end{proposition}

Since $k^{(m)}_D(0,a)\le k^{(3)}_D(0,a)\le h_D$ for $3\ge m\le\infty$
($k_D=k^{(\infty)}_D$) and $\kappa^{(l)}_D(0;a)\le
\kappa^{(2)}_D(0;a)$ for $2\le l\le\infty$
($\hat\kappa_D=\kappa^{(\infty)}_D$), for these $m$ and $l$ it follows that
$\hat h_D(a)=h_D(a)\Leftrightarrow
(k^{(m)}_D(0,a))^\ast=h_D(a)\Leftrightarrow \kappa^{(2)}_D(0;a).$
\smallskip

\noindent{\bf Remark.} We do not know whether the number 3 can be
replaced by 2 (it cannot be replaced by 1 according to Proposition
\ref{bala.pr1}(iv)).
\smallskip

\beginproof The implication (i)$\Rightarrow$(ii) follows from Proposition
\ref{bala.pr1}.

Assume (iii). If $a_1+a_2=a,$ then by
$\kappa^{(2)}_D(0;a)\le\kappa_D(0;a_1)+\kappa_D(0;a_2)$ and
$\kappa_D(0;\cdot)\le h_D$ it follows that $h_D(a)\le
h_D(a_1)+h_D(a_2),$ so (i) holds.

It remains to prove (ii)$\Rightarrow$(iii). We first prove that (ii) implies
\begin{equation}\label{bala1}(k^{(2)}_D(0,\lambda a))^\ast=|\lambda|h_D(a),\
\lambda\in\Bbb D.
\end{equation}
We can assume that $h_D(a)\neq 0.$ Considering the analytic disc
$\phi(\zeta)= a\zeta/h_D(a)$ as a competitor *\footnote{*This
means that $\phi$ belongs to the set over which we take the
infimum in the definition of $l_D$} for $l_D(\lambda a,a),$ we get
$$l_D(\lambda a,a)\le m(h_D(\lambda a),h_D(a)).$$  Hence by the inequality
$$p(0,h_D(a))=k^{(3)}_D(0,a)\le k^{(2)}_D(0,\lambda a)+k^{(1)}_D(\lambda a,a)$$
we get
$$p(0,|\lambda|h_D(a))=p(0,h_D(a))-p(|\lambda|h_D(a),h_D(a))\le
k^{(2)}_D(0,\lambda a).$$ Then $$(k^{(2)}_D(0,\lambda
a))^\ast\ge|\lambda|h_D(a).$$ It remains to note that the opposite
inequality is always true.

Now (\ref{bala1}) shows that $$\lim_{\lambda\to
0}\frac{k^{(2)}_D(0,\lambda a)}{|\lambda|}=h_D(a).$$ On the other
hand, by Proposition \ref{der.pr2} this limit does not exceed
$\kappa^{(2)}_D(0;a)\le h_D(a)$, so (iii) is proved.\qed
\smallskip

Having in mind Propositions \ref{bala.pr1} and \ref{bala.pr2}, it
is natural to ask whether the minimality (rather than maximality)
of some $k^{(m)}_D(0,a),$ i.e. $l_D(0,a)=h_G(a)$ for a domain
$G\supset D$ implies some "convex"\ property. We have the
following

\begin{proposition}\label{bala.pr8} Let $D\subset\Bbb C^n$ be a bounded
balanced domain and $G\subset\Bbb C^n$ be a pseudoconvex balanced
domain containing $D$. Suppose that $h_D$ is continuous at some
$a\in D,$ $h_G(a)\neq 0$ and $\overline G$ does not contain
(nontrivial) analytic discs through $a/h_G(a).$ Then the following
are equivalent:

(i)  $h_D(a)=h_G(a);$

(ii) $l_D(0,a)=h_G(a);$

(iii) $\kappa_D(0;a)=h_G(a).$
\end{proposition}

\beginproof
It suffices to prove that $$l_D(0,a)=h_G(a)\Rightarrow h_D(a)\le
h_G(a), \kappa_D(0;a)=h_G(a)\Rightarrow h_D(a)\le h_G(a).$$ Let
$(\varphi_j)\subset\mathcal O(\Bbb D,D)$ and $\alpha_j\to h_G(a)$
so that $\varphi_j(0)=0$ and $\varphi_j(\alpha_j)=a$
(correspondingly, $\alpha_j\varphi'_j(0)=a$). Expressing
$\varphi_j$ in the form
$\varphi_j(\lambda)=\lambda\psi_j(\lambda),$ by maximality
principle $h_G\circ\psi_j\le 1$ so $\psi_j\in\mathcal O(\Bbb
D,\overline G).$ As $D$ is bounded, then by going to a subsequence
we can assume that $\varphi_j\to\varphi\in\mathcal O(\Bbb
D,\overline D)$, hence $\psi_j\to\psi\in\mathcal O(\Bbb
D,\overline G).$ In particular,
$$\psi(h_G(a))=\lim_{j\to\infty}\psi_j(\alpha_j)=\lim_{j\to\infty}
\frac{a}{\alpha_j}=\frac{a}{h_G(a)}=:b\ \mbox{(and }\psi(0)=b),$$
respectively. On the other hand, as $\overline G$ does not contain
analytic discs through $b,$ we get $\psi(\Bbb D)=b.$ The
continuity of $h_D$ at $b$ implies that
$$1>h_D(\varphi_j(\lambda))\to|\lambda|h_D(b),\ \lambda\in\Bbb
D.$$ When $\lambda\to 1$ we get $h_D(b)\le 1,$ ò.å $h_D(a)\le
h_G(a).$\qed
\smallskip

\noindent{\bf Remarks.} a) Since the holomorphic hull $\mathcal E(D)$
of a balanced domain $D$
is a balanced domain (see e.g.
\cite[Remark 3.1.2(b)]{Jar-Pfl2}), the above result can be also
applied for $G=\mathcal E(D)$. Of course,
it can be also applied for $G=\hat D.$

b) If $h_G$ is continuous near $a$ and $\partial G$ does not
contain analytic discs through $a/h_G(a),$ then by maximality
principle it follows that $\overline G$ does not contain analytic
discs through $a/h_G(a)$ either.

â) In connection with Proposition \ref{bala.pr8} it is natural to
ask whether if $h_D=l_D(0,\cdot)$ for a balanced domain $D,$ then
it needs to be pseudoconvex. The answer of this question is
unknown to us.
\smallskip

The next example shows that in Proposition \ref{bala.pr8} the continuity
assumption for $h_D$ is essential.

\begin{example}\label{bala.ex9} If $D=\Bbb D^2\setminus\{(t,t):|t|\ge 1/2\},$
$d=(t,t),$ $|t|<1/2,$ then
$$h_D(d)=2|t|,\hbox{\ \ but\ \ }l_D^\ast(0,d)=|t|=h_{\Bbb D^2}(d).$$
On the other hand, $\mathcal E(D)=\D^2$ and $\overline{\Bbb D^2}$
does not contain analytic discs through any point from
$\partial\Bbb D\times\partial\Bbb D.$
\end{example}

\beginproof We need to prove only that $$l_D(0,d)\le|t|.$$
For each $r\in(|t|,1)$ we can choose $\alpha\in\Bbb D$ so that
$t=\varphi(t/r),$ where
$\varphi(\lambda)=\lambda\frac{\lambda-\alpha}{1-\overline{\alpha}\lambda}.$
Then the disc $\psi(\zeta)=(r\zeta,\varphi(\zeta))$ is a
competitor for $l_D(0,d),$ whence it follows that
$l_D(0,d)\le|t|/r.$ It remains to leave $r\to 1.$\qed
\smallskip

\noindent{\it Appendix.} Note that even $$l_D(0,\cdot)=l_{\Bbb
D^2}(0,\cdot).$$ To this aim it suffices to prove that
$l_D(0,a)\le|a_1|$ for $a=(a_1,a_2)\in D,$ $a_1\neq a_2,$
$|a_1|\ge|a_2|.$ We get this easily by considering
$\psi(\lambda)=(\lambda,\lambda a_2/a_1)$ as a competitor for
$l_D(0,a).$

On the other hand, if $a_1=(0,b)$ and $a_2=(b,0),$ $b\in\Bbb D,$ then
$$l_D(a_1,a_2)=l_{\Bbb D^2}(a_1,a_2)\Leftrightarrow |b|\le 4/5.$$

Indeed, using the M\"obius transformation
$\psi_b(\lambda)=\frac{\lambda-b}{1-\overline{b}\lambda},$
we get $l_D(a_1,a_2)=l_{D_b}(0,a),$ where $a=(b,-b)$ and
$D_b=\Bbb D^2\setminus\{(\psi_b(\lambda),\lambda):1/2\le |\lambda|<1\}.$

For $|b|<4/5$ we easily check that
$\varphi=(\mbox{id},-\mbox{id})\in\mathcal O(\Bbb D, D_b).$ Then
$l_{D_b}(0,a)\le|b|$ so $l_D(a_1,a_2)=l_{\Bbb D^2}(a_1,a_2).$

To get this for $|b|=4/5,$ it suffices to consider $r\varphi$ for
$r\in(0,1)$ as a competitor for $l_{D_b}(0,a)$, then to leave
$r\to 1.$

Now assume that $l_D(a_1,a_2)=l_{\Bbb D^2}(a_1,a_2)$ for
$|b|>4/5.$  Then we can find discs $\varphi_j\in\mathcal O(\Bbb
D,D_b),$ such that $\varphi_j(0)=0$ and $\varphi_j(\alpha_j)=a,$
where $\alpha_j\to b.$ The Schwarz--Pick Lemma implies that
$\varphi_j\to\varphi.$ On the other hand, $\varphi(\Bbb
D)\cap\{(\psi_b(\lambda),\lambda):1/2<|\lambda|<1\}$ is a
singleton, contradicting the Hurwitz Theorem.
\smallskip

\noindent{\bf Remark.} By \cite[Theorem 3.4.2]{Jar-Pfl1}
(see also \cite{Pol-Sha}) it follows that if
$D_n=\Bbb D^n\setminus\{(t,\dots,t):|t|\ge 1/2\},$ $n\ge 3,$
then $l_{D_n}=l_{\Bbb D^n}.$
\smallskip

The next example shows that in Proposition \ref{bala.pr8}
the assumption on discs is essential.

\begin{example}\label{bala.ex10} Let $0<a<1$ and
$$D=\{z\in\Bbb D^2:|z_2|^2-a^2<2(1-a^2)|z_1|\}.$$
Then $D$ is a balanced Reinhardt domain, $h_D$ is continuous
function and $\mathcal E(D)=\Bbb D^2$
(see e.g. \cite{Jar-Pfl2}).
On the other hand, if $c=(0,d)$, $|d|<a,$ then
$$h_D(c)=|d|/a>l_D(0,c)=|d|=h_{\Bbb D^2}(c).$$
\end{example}

\beginproof We have to show just that $l_D(0,c)\le d$
for $d\in(0,a).$
To this aim it suffices to show that
$\varphi=(\psi,\mbox{id})\in\mathcal O(\Bbb D,D),$ where
$\psi(\lambda)=\lambda\frac{\lambda-d}{1-d\lambda}.$ It is easily
seen that $|\psi(\lambda)|\ge x\frac{x-d}{1-xd}$ for $x=|\lambda|$
and is suffices to check that
$$x^2-a^2<2(1-a^2)\frac{x(x-d)}{1-dx},\hbox{ i.e.}$$
$$dx^3+(1-2a^2)x^2-d(2-a^2)x+a^2>0.$$ This is clear for $x=0$.
Since $x\in (0,1)$ and $d\in(0,a)$,
we need to prove that $$ax^3+(1-2a^2)x^2-a(2-a^2)x+a^2 \ge 0,$$
which is equivalent to the obvious inequality $(x-a)^2(ax+1)\ge 0.$\qed
\smallskip

\noindent{\bf Remark.} Some propositions and examples in the
spirit of the above for $k^{(m)}$ can be found in the paper
\cite{Nik-Pfl9} of the author and P. Pflug.

\setcounter{equation}{0}
\section{Kobayashi--Buseman metric}\label{kob-buz}

The main aim of this section is to prove that the
Kobayashi--Buseman metric for an arbitrary domain equals the
Kobayashi metric of order $2n-1$ and this number is the least
possible. A similar result for $2n$ instead of $2n-1$ is contained
in the work \cite{KobS1} of S. Kobayashi, where this metric is
introduced.

\begin{theorem}\label{kb.th} For each domain $D\subset\Bbb C^n$ one has
\begin{equation}\label{kb1}
\kappa_D^{(2n-1)}=\hat{\kappa}_D.
\end{equation}

On the other hand, if $n\ge 2$ and $$D_n=\{z\in\Bbb
C^n:\sum_{j=2}^n(2|z_1^3-z_j^3|+|z_1^3+z_j^3|)<2(n-1)\},$$ then
\begin{equation}\label{kb2}
\kappa_{D_n}^{(2n-2)}(0;\cdot)\neq\hat\kappa_{D_n}(0;\cdot).
\end{equation}
\end{theorem}

The proof below shows that the identity (\ref{kb1}) remains true
for an arbitrary $n$-dimensional complex manifold.

Theorems \ref{kb.th} and \ref{der.th1} lead to the following

\begin{corollary}\label{kb.cor2} For every taut domain $D\subset\Bbb C^n$
one has
$$
\lim_{w\to z}\frac{k_D^{(2n-1)}(z,w)}{k_D(z,w)}=1
$$
locally uniformly on $z.$ The number $2n-1$ is the least possible
in the general case.
\end{corollary}

\noindent{\bf Remarks.} a) Corollary \ref{kb.cor2} remains true
for an arbitrary $n$-dimensional complex taut manifold.

b) Corollary \ref{kb.cor2} can be viewed as an affirmative answer
of the infinitesimal version of a question of S. Krantz
\cite{Kra1}: For an arbitrary strictly pseudoconvex domain
$D\in\Bbb C^n$, is there some $m=m(D)\in\Bbb N,$ such that
$k_D=k_D^{(m)}$? Unlike the infinitesimal case,  $m$ cannot depend
only on $n,$ as shown in \cite[p. 109]{Jar-Pfl1}.)
\smallskip

For $z\in D\subset C^n$, denote by $I_{D,z}$ the indicatrix of
$\kappa_D(z;\cdot)$, i.e. $I_{D,z}=\{X\in\Bbb
C^n:\kappa_D(z;X)<1\}.$ Note that $I_{D,z}$ is a balanced domain.
In particular, it is starlike with respect to the origin. Then the
identity $\kappa_D^{(2n)}=\hat\kappa_D$ is obtained from the
following application of a lemma of C. Carath\'eodory (see e.g.
\cite{KobM}):
\begin{equation}\label{kb3}
\aligned\hat h_S=\inf\{\sum_{j=1}^mh_S(X_j):m\le 2n,\
\sum_{j=1}^mX_j=X,\\ X_1,\dots, X_m\mbox{ are }\Bbb
R\mbox{-linearly independent}\},\endaligned
\end{equation}
where $h_S$ and $\hat h_S$ are the functions of Minkowski of an
arbitrary domain $S\subset\C^n$ that is starlike with respect to
the origin (i.e. $ta\in S$ for $a\in S$ and $t\in[0,1]$) and of
its convex hull $\hat S$, respectively (it is easily seen that in
this case the number $2n$ is the least possible).

In order to replace the number $2n$ by $2n-1,$ we will use the
fact that $I_{D,z}$ is balanced rather than starlike. For $m\in\N$
we put
$$h_S^{(m)}(X)=\inf\{\sum_{j=1}^mh_S(X_j):\sum_{j=1}^mX_j=X\}.$$

\noindent{\it Proof of (\ref{kb1}).} Follows directly from the following

\begin{proposition}\label{kb.prop3} If $B\subset\Bbb C^n$ is a balanced
domain, then
\begin{equation}\label{kb4}
\hat h_B=h_B^{(2n-1)}.
\end{equation}
\end{proposition}

To prove Proposition \ref{kb.prop3} we need the following

\begin{lemma}\label{kb.lm4} Every balanced domain can be exhausted
by bounded balanced domains with continuous Minkowski functions.
\end{lemma}

\beginproof Let $B\subset\Bbb C^n$ be a balanced domain. For $z\in\Bbb C^n$
and $j\in\Bbb N$ we put $F_{n,j,z}=\overline{\Bbb
B_n(z,||z||^2/j)}$ ($\Bbb B_n\subset(a,r)\C^n$ is the ball of
center $a$ and radius $r$). We can assume that $\Bbb
B_n(0,1)\Subset B.$ Let
$$ B_j=\{z\in \Bbb B_n(0,j): F_{n,j,z}\subset B\},\quad j\in\Bbb N.
$$
Then $(B_j)$ is an exhaustion of $B$ by nonempty bounded open sets.
We will show that
$B_j$ is a balanced domain with continuous
Minkowski function $h_{B_j}.$

To this aim let us note that if $z\in B_j$ and
$\lambda\in(\overline{\Bbb D})_\ast,$ then $F_{n,j,\lambda
z}\subset\lambda F_{n,j,z}\subset B.$ Now it easily follows that
$B_j$ is a balanced domain.

As $h_{B_j}$ is upper semicontinuous, it remains to prove that it
is also lower semicontinuous. Assuming the contrary, we can find a
sequence of points $z_k$ tending to some $z,$ and a number $c>0$
such that $h_{B_j}(z_k)<1/c<h_{B_j}(z)$ for each $k.$ Note that
$F_{n,j,cz_k}\subset B$, so $\Bbb B_n(cz,c^2\|z\|^2/j)\subset B$.
On the other hand, let us choose $t\in (0,1)$ such that
$h_{B_j}(tcz)>1$. Then $F_{n,j,tcz} \subset\Bbb
B_n(cz,c^2\|z\|^2/j)\subset B$, so $h(tcz)<1$ -- a
contradiction.\qed
\smallskip

\noindent{\it Proof of Proposition \ref{kb.prop3}.} We will first
prove (\ref{kb4}) in the case when $B\subset\Bbb C^n$ is a bounded
balanced domain with a continuous Minkowski function. Let us fix a
vector $X\in(\Bbb C^n)_\ast.$ Then $\hat h_B(X)\neq 0$ and we can
assume that $\hat h_B(X)=1.$ As $h_B$ is continuous, by
(\ref{kb3}) there exist $\Bbb R$-linearly independent vectors
$X_1,\dots,X_m$ ($m\le 2n$) such that $\sum_{j=1}^m X_j=X$ and
$\sum_{j=1}^mh_B(X_j)=1.$ As $\hat h_B$ is a norm and $\hat h_B\le
h_B,$ by triangle inequality $h_B(X_j)=\hat h_B(X_j),$\;
$j=1,\dots,m.$ To prove (\ref{kb4}), it suffices to show that
$m\neq 2n.$ Let $H$ be a support hyperplane for $\hat B$ at
$X\in\partial\hat B.$ We can assume that $H=\{z\in\Bbb
C^n:\Re\langle z-X,X_0\rangle=0\},$ where $X_0\in\Bbb C^n$
($\langle\cdot,cdot\rangle$ is the Hermitian scalar product).
Suppose that $m=2n.$ Then $H=\{\sum_{j=1}^m\alpha_j
X_j/h_B(X_j):\sum_{j=1}^m\alpha_j =1,\
\alpha_1,\dots,\alpha_m\in\Bbb R\}.$ In particular, $\partial\hat
B$ contains a relatively open set in $H.$ As $\hat B$ is a
balanced domain, its intersection with the complex line through X,
directed at $\overline X_0,$ is a disc containing a line segment
within its boundary. This contradiction proves (\ref{kb4}) for a
bounded balanced domain with a continuous Minkowski function.

Now let $B\subset\Bbb C^n$ be an arbitrary balanced domain. If
$(B_j)$ is an exhaustion of $B$ as in Lemma \ref{kb.lm4}, then
$h_{B_j}\searrow h_B$ pointwise. Then (\ref{kb3}) shows that $\hat
h_{B_j}\searrow\hat h_B.$ Now (\ref{kb3}) follows from the
inequalities $\hat h_B\le h_B^{(2n-1)}\le h_{B_j}^{(2n-1)}$ and
the equality $\hat h_{B_j}=h_{B_j}^{(2n-1)}$ from above.\qed
\smallskip

\noindent{\it Proof of (\ref{kb2}).} Observe that the domain $D_n$
from Theorem \ref{kb.th} is pseudoconvex and balanced. Then
$\kappa_{D_n}(0;\cdot)=h_{D_n}$ (see Proposition \ref{bala.pr1}
(iv)) so $\kappa_{D_n}^{(m)}(0;\cdot)=h_{D_n}^{(m)}.$ Thus
(\ref{kb2}) is equivalent to
\begin{equation}
\hat h_{D_n}\neq h_{D_n}^{(2n-2)}.
\end{equation}

To prove this inequality, let $L_n=\{z\in \Bbb C^n:z_1=1\}.$ By
triangle inequality $D_n\subset\Bbb D\times\Bbb C^{n-1}$ and
$$
F_n:=\partial D_n\cap L_n=\{z\in\Bbb C^n:z_1=1,\ z_j^3=1\},\ 2\le j\le n\}.
$$
Hence $\partial\hat D_n\cap L_n=\hat
F_n=\{1\}\times\hat\Delta{n-1},$ where $\Delta$ is the triangle of
vertices $1,$ $e^{2\pi i/3},$ $e^{4\pi i/3}$ together with its
interior. Note that $\partial\hat D_n\cap L_n$ is a
$(2n-2)$-dimensional convex set. Put $\tilde F_n=\{Y\in\hat
F_n:h_{D_n}^{(2n-2)}(Y)=1\}.$ If $X\in\tilde F_n,$ then there
exist vectors $X_1,\dots,X_m\in(\Bbb C^n)_\ast,$ $m\le 2n-2$ such
that $\sum_{j=1}^m X_j=X$ and $\sum_{j=1}^m h_{D_n}(X_j)=1$ (as
$D_n$ is a taut domain). Then
$X_1/h_{D_n}(X_1),\dots,X_m/h_{D_n}(X_m)\\ \in F_n$ and the convex
hull of these vectors contains $X.$ As $F_n$ is a finite set, it
is contained $\tilde F_n$ in a finite union of not more than
$(2n-3)$-dimensional convex sets. So $\hat F_n\neq\tilde F_n,$
which shows that $\hat h_{D_n}\neq h_{D_n}^{(2n-2)}.$\qed
\smallskip

Thus Theorem \ref{kb.th} is proved.

\setcounter{equation}{0}
\section{Interpolation in the Arakelian theorem}\label{arak}

The aim of this section is to prove a general statement for
approximation and interpolation over the so-called Arakelian sets.
This statement will be used in the proof of Theorem \ref{decr.th1}
from the next section.

Let us first recall the well-known theorem of Mergelian that
generalizes the theorems of Weierstrass and Runge.

\begin{theorem}\label{mer.th1} The complement of a compact $K\subset C$
is a connected set if and only if
each continuous function on $K$ that is holomorphic in the interior of $K$
can be uniformly approximated on $K$ by polynomials.
\end{theorem}

The most popular generalization of Theorem \ref{mer.th1} belongs
to N. Arakelian.*\footnote{*After the proof by N. Arakelian of
Theorem \ref{mer.th2}, J.-P. Rosay and W. Rudin \cite{Ros-Rud}
showed how this theorem follows from the Mergelian theorem
itself.}

A relatively closed subset $E$ of a domain $D\subset\Bbb C$ is
called an Arakelian set, if $D^\ast\setminus E$ is connected and
locally connected, where $D^\ast$ is the one-point
compactification of $D$.

Denote by $\CA(E)$ the set of continuous functions on $E$ that are
holomorphic in the interior $E^0$ of $E.$

\begin{theorem}\label{mer.th2}\cite{Ara} A relatively closed subset $E$
of a domain $D\subset\Bbb C$ is an Arakelian set if and only if
each function from $\CA(E)$ can be uniformly approximated on $E$
by holomorphic functions from $D.$
\end{theorem}

The next result, proven independently by P. M. Gauthier and W.
Hengartner and A.~Nersesyan, provides an opportunity for
approximation in Theorem \ref{mer.th2}.

\begin{theorem}\label{mer.th3}\cite{Gau-Hen,Ner}
Let $D\subset\Bbb C$ be a domain, let  $E\subset D$ be an
Arakelian set, and let $\Lambda$ be a sequence of points in
$E\setminus E^0$ without an accumulation point in $D$. For every
$\lambda\in\Lambda$, a finite sequence
$(\beta_{\lambda}^{\nu})_{\nu=1}^{\nu(\lambda)}$ of complex
numbers is given. Then for each $f\in\CA(E)$ and each
$\varepsilon>0$, there exists a $g\in\CO(D)$ such that
$|g(z)-f(z)|<\varepsilon$ for $z\in E,$ $g(\lambda)=f(\lambda)$
and $g^{(\nu)}(\lambda)=\beta_{\lambda}^{\nu}$ for
$\lambda\in\Lambda$ and $\nu=1,\dots,\nu(\lambda).$
\end{theorem}

Now let us formulate an extension of Theorem \ref{mer.th3}.

\begin{theorem}\label{mer.th4} Let  $D$, $E$, $\Lambda$ (possibly
$\Lambda=\varnothing$), $\beta_\lambda^\nu$ be as in Theorem
\ref{mer.th3} and let  $b_1,\dots b_k\in E^0$. Then for each
$f\in\CA(E)$, $\varepsilon>0$ and $m\in\NN^\ast$ there exists a
$g\in\CO(D)$ with the properties of Theorem \ref{mer.th3} and such
that $g^{(\nu)}(b_j)=f^{(\nu)}(b_j)$ for $j=1,\dots,k$ and
$\nu=0,\dots,m$.
\end{theorem}

\beginproof We can clearly assume that $E\neq D$.

The proof will be divided into four steps.

{\it Step 1.} For each $j=1,\dots,k$, there is a function
$s_j\in\CO(D)$, bounded on $E$ and such that $s_j'(b_j)\neq 0$,
$s_j(b_j)=0$ and $s_j(b_q)\neq 0$ for an arbitrary $q\neq j$.
*\footnote{*This is clear if $D$ is biholomorphic to a bounded
domain; in particular, if $\overline D\neq \C$.}

Indeed, choose a point $c\in D\setminus E.$ As $E\cup\{c\}\subset
D$ is an Arakelian set, Theorem \ref{mer.th2} implies the
existence of a $\tilde s\in\CO(D)$ such that $|\tilde s|<1$ on $E$
and $|\tilde s(c)-2|<1$. Put $\hat s_j=\tilde s-\tilde s(b_j)$. As
$\hat s_j(c)\neq 0,$ we get $\hat s_j\not\equiv 0$. Now as $|\hat
s_j|<2$ on $E$, the function
$$s_j(z)=\frac{(z-b_j)\hat s_j(z)}{\prod_{q=1}^k(z-b_q)^{\ord_{b_q}\hat{s}_j}},\quad z\in D,
$$ has the required properties.

{\it Step 2.} There exists a function $p\in\CO(D),$ bounded on $E$
and such that $p(b_j)\neq 0$ for $j=1,\dots,k$ and $\ord_\lambda
p\ge\nu(\lambda)+1$ for an arbitrary $\lambda\in\Lambda$.

Indeed, if $q=0$ on $E$ and $q(c)=1,$ where $c\in D\setminus E,$
we can apply Theorem \ref{mer.th3} for $E\cup\{c\}$, $q$,
$\varepsilon=1$ and $\beta_{\lambda}^{\nu}=0$,
$\nu=1,\dots,\nu(\lambda)+1$, $\lambda\in\Lambda$. Thus we get a
nonconstant function $\tilde p\in\CO(D)$ such that $|\tilde p|<1$
on $E$ and $\ord_\lambda\tilde p\geq\nu(\lambda)+1$,
$\lambda\in\Lambda$. It remains to put
$$p(z)=\frac{\tilde p(z)}{\prod_{j=1}^k(z-b_j)^{\ord_{b_j}\tilde p}},\quad z\in D.$$

{\it Step 3.} Let $s_j$ be the function from Step 1, $j=1,\dots
k,$ and let $p$ be the function from Step 2. For each
$\nu\in\NN^\ast$ we put
$$\tilde h_j^\nu=\frac{p}{s_j}\prod_{q=1}^ks_q^{\nu+1}, $$
Then
$$h_j^\nu=\frac{\tilde h_j^\nu}{(\tilde h_j^\nu)^{(\nu)}(b_j)} $$
is well defined on $D.$ The function $$M_\nu=\sup_E\sum_{j=1}^{k}|h_j^\nu|$$
will be also needed in the last step.

{\it Step 4.} We are ready to prove the theorem by induction on $m.$

Let $m=0$ and $g$ be the function from Theorem \ref{mer.th3} for
$\Lambda$, $(\beta_\lambda^\nu)_{\nu=1}^{\nu(\lambda)}$ and
$\tfrac{\varepsilon}{M_0+1}$. It is easily checked that the function
$$g_0=g+\sum_{j=1}^k(f(b_j)-g(b_j))h_j^0$$ has the required properties.

Put $d=\min_{1\le j\le k}\hbox{dist}(b_j,\Bbb C\setminus E^0)$.
Assume that Theorem \ref{mer.th4} is true for some $m\ge 0$ and let $g_m$ be
the corresponding function for $\varepsilon(1+M_{m+1}(m+1)!d^{-m-1})^{-1}$.
By the Cauchy Inequality, the function
$$g_{m+1}=g_m+\sum_{j=1}^k(f^{(m+1)}(b_j)-g_m^{(m+1)}(b_j))h_j^{m+1}$$
has the required properties for $m+1.$

This finishes the induction step. The theorem is proven.\qed

\setcounter{equation}{0}
\section{Generalized Lempert function}\label{lem-decr}

In this section we define the generalized Lempert function
(introduced by D. Coman \cite{Com2}) and prove that it decreases
under adding poles. This function is introduced as an easier and
more flexible (in some sense) version of the so-called generalized
(pluricomplex) Green function (see e.g. \cite{Jar-Pfl3}).

Let $D\subset\CC^n$ be a domain and $p\gneq 0$ be a function on $D.$ Put
$$\Sp=\{a\in D:\bs p(a)>0\}.$$
For $z\in D$ we define
\begin{multline*}
l_D(\bs p,z)=\inf\{\prod|\lambda_{\psi,a}|^{\bs
p(a)}:\exists\psi\in\CO(\DD,D),  \psi(0)=z\\
\psi(\lambda_{\psi,a})=a \text{ for each } a\in\Sp\}
\end{multline*}
(for any $a\in\Sp$ we take one $\lambda_{\psi,a}$).

From the proof of Theorem \ref{decr.th1} below it follows that
such a $\psi$ exists, if $\Sp$ is finite or countable. If $\Sp$ is
uncountable and such a $\psi$ exists, then it is easily seen that
$$0=l_D(\bs p,z)=\inf\{l_D(\bs p_B,z):B \subset \Sp,\
0<\#B<\infty\},$$ where $p_B=p\chi_B.$

If there is no such $\psi,$ we can define
$$l_D(\bs p,z)=\inf\{l_D(\bs p_B,z): B\subset \Sp,\ 0<\#B<\infty\}.$$

The function $l_D(\bs p,\cdot)$ so introduced is called Lempert
function of $D$ with respect to $\bs p$ (a generalized Lempert
function). If $A$ is a nonempty subset of $D$ and $\chi_A$ is its
characteristic function, then we put $l_D(A,z)= l_D(\chi_A,z).$
This function is called the Lempert function with poles in $A.$
Let us note that $l_D(\{a\},z)$ is the usual Lempert function
$l_D(a,z).$

Using the Lempert theorem, F. Wikstr\"om \cite{Wik1} showed that
if $A$ and $B$ are subsets of a convex domain $D\subset\CC^n$ such
that $A\subset B,$ then $l_D(B,\cdot)\le l_D(A,\cdot),$ i.e. the
Lempert function decreases under adding poles.

On the other hand, in \cite{Wik2} there is an example of a complex
space not satisfying the inequality (under the same definition of
a Lempert function) and it is asked whether this inequality is
true for arbitrary domains in $\CC^n$.

The main aim of this section is to give an affirmative answer to
this question. We will use Theorem \ref{mer.th4}, that is, the
possibility for interpolation in the Arakelian approximation
theorem.

\begin{theorem}\label{decr.th1} If $D\subset\CC^n$ is a domain and
$\bs p\gneq 0$ is a function on $D,$ then
$$l_D(\bs p,\cdot)=\inf\{l_D(\bs p_B,\cdot):
B\subset\Sp,\ 0<\#B<\infty\}.$$ In particular, $l_D(\bs
p,\cdot)=\inf\{l_D(\bs p_B,\cdot):\varnothing\neq B\subset\Sp\}.$
\end{theorem}

\begin{corollary}\label{decr.cor1} If $D\subset\CC^n$ is a domain and
$\bs p,\bs q$ are functions on $D$ such that $0\lneq\bs p\leq\bs
q,$ then $l_D(\bs q,\cdot)\leq l_D(\bs p,\cdot)$.
\end{corollary}

\beginproof By the above remark, the theorem follows in the case when
$\Sp$ is uncountable.

Now let $\Sp=(a_j)_{j=1}^l$ ($l\in\Bbb N^\ast$) be a countable or
finite nonempty set. Let $z\in D.$

We first prove the inequality
\begin{equation}\label{gl1}
l_D(\bs p,z)\le\inf\{l_D(\bs p_B,z):B\subset\Sp,\ 0<\#B<\infty\}.
\end{equation}
Let $B\neq\varnothing$ be a finite subset of $\Sp.$ We can assume that
$B=A_m:=(a_j)_{j=1}^m$ for some $m\le l.$

Let us consider an arbitrary $\phi:\Bbb D\to D$ such that
$\phi(\lambda_j)=a_j,$ $0\le j\le m,$ where $\lambda_0=0$ and
$a_0=z.$ Let $t\in[\max_{0\le j\le m}|\lambda_j|,1)$ and
$\lambda_j=1-(1-t)/j$,\; $j\in A(m)$, where $A(m)=\{m+1,\dots,l\}$
for $l<\infty$ and $A(m)=\{j\in\NN:j>m\}$ for $l=\infty$. Consider
a continuous curve $\phi_1:[t,1)\to D$ such that
$\phi_1(t)=\phi(t)$ and $\phi_1(\lambda_j)=a_j$, $j\in A(m)$. Put
$$
f=\begin{cases}
\phi|_{\overline{t\Bbb D}}\\
\phi_1|_{[t,1)}
\end{cases}ò
$$
on $F_t=\overline{t\Bbb D}\cup[t,1)\subset\DD$. Clearly $F_t$ is
an Arakelian set for $\Bbb D,$ $f\in\mathcal A(F_t,D)$ and
$\Lambda=(\lambda_j)_{j=1}^l$ satisfies the conditions in Theorem
\ref{mer.th4}. Let $d(z)=\dist(f(z),\partial D),$ $z\in F_t,$
where the distance is generated by the $L^\infty$-norm. Choose a
continuous real-valued function $\eta$ on $F_t$ such that
$$\eta\le\log d\mbox{ on }[t,1),\qquad \eta=\min_{\overline{t\DD}}
\log d\mbox{ on }\overline{t\DD}.$$ By theorem \ref{mer.th3},
there exists a $\zeta\in\CO(\DD)$ such that $|\zeta-\eta|<1$
on $F_t.$ By Theorem \ref{mer.th4}, applied to the components of
$e^{\zeta-1}f,$ one can find a $q_t\in\mathcal O(\Bbb D)$ such that
$q_t(\lambda)=f(\lambda),$ $\lambda\in\Lambda$ and
$$||q_t-f||<|e^{\zeta(z)-1}|<e^{\eta(z)}\le d(z),\quad z\in F_t.$$
Thus $q_t(F_t)\subset D$ and so there exists a simply connected
domain $E_t$ such that $F_t\subset E_t\subset\DD$ and
$q_t(E_t)\subset D$.

Let $\rho_t:\DD\to E_t$ be the corresponding Riemann (conformal)
mapping, satisfying $\rho_t(0)=0,$ $\rho_t'(0)>0$ and
$\rho_t(\lambda_j^t)=\lambda_j$. Considering the analytic discs
$q_t\circ\rho_t:\DD\to D$ we get
$$
l_D(\bs p,z)\le\prod_{j=1}^l|\lambda_j^t|^{\bs
p(a_j)}\leq\prod_{j=1}^m|\lambda_j^t|^{\bs p(a_j)}.
$$
Note that by the Carath\'eodory Kernel Theorem, $\rho_t$ for $t\to
1$ tends locally uniformly on $\DD$ to $\id.$ Hence the latter
product above tends to $\prod_{j=1}^m|\lambda_j|^{\bs p(a_j)}$.
Since $\phi$ was an arbitrary competitor for $l_D(\bs
p|_{A_m},z)$, we get the inequality (\ref{gl1}).

On the other hand, the existence of analytic discs containing $z$
and $\Sp$ easily implies
$$l_D(\bs p,z)\ge\limsup_{m\to\infty}l_D(\bs p|_{A_m},z),$$
which concludes the proof of the theorem.\qed
\smallskip

\noindent{\bf Remark.} The Lempert function does not decrease
strictly under adding of poles; for example \cite[Theorem
2.1]{Dieu-Trao} shows that
$$l_{\Bbb D^2}(\{a_1,a_2\}\times\{a_1\},0)=|a_1|=l_{\Bbb D^2}(\{a_1\}\times\{a_1\},0).$$

The next example shows that our definition of a generalized
Lempert function, in the case of nonexistence of a corresponding
disc, is more "sensitive"\ than that from \cite{Jar-Pfl3} (where
in this case the function is fixed to $1$).
\smallskip

\noindent{\bf Example.} Let $A\subset\DD$ be an uncountable set.
Then there is no analytic disc $\phi\in\CO(\DD,\DD^2),$ containing
$A\times\{0\}$ and $(0,w)$, $w\in\DD_\ast$.

Let $B$ be an arbitrary finite subset of $A.$ From \cite[Theorem 2.1]{Dieu-Trao},
$$l_D(B\times\{0\},(0,w))=\max\{l_{\DD}(B,0),l_{\DD}(0,w)\}=
\max\{\prod_{b\in B}|b|,|w|\}.$$ So $l_D(A\times\{0\},(0,w))=|w|.$
\smallskip

Finally let us note that the generalized Lempert function is
clearly biholomorphically invariant, but in general not
contractible under holomorphic mappings even when they are proper
coverings.
\smallskip

\noindent{\bf Example.} Let $\pi(z)=z^2$. Clearly $\pi:\D_\ast\to
\D_\ast$ is a proper covering ($D_\ast=\D\setminus\{0\}$). Let
$a_1=-a_2\in\DD_\ast,$ $c=a_1^2$ and $z\in\DD_\ast,$ $z\neq
a_1,a_2.$ By \cite[Theorem 3.3.7]{Jar-Pfl1}
$$l_{\DD_\ast}(c,z^2)=\min\{l_{\DD_\ast}(a_1,z),l_{\DD_\ast}
(a_2,z)\}>l_{\DD_\ast}(a_1,z)l_{\DD_\ast}(a_2,z).$$ On the other
hand, by \cite[Theorem 2.1]{Dieu-Trao} the last product equals
$l_{\DD_\ast}(\{a_1,a_2\},z)$. Therefore $$l_{\DD_\ast}(\bs p
,\pi(z))>l_{\DD_\ast}(\bs p\circ\pi,z)\mbox{ for }\bs
p=\chi_{\{c\}}.$$

We conclude this section with the following comment. The proof of
Theorem \ref{decr.th1} is contained in the paper \cite{Nik-Pfl5}
by the author and P. Pflug. Later, based on the same idea, F.
Forstneric and J. Winkelman \cite{For-Win} proved that, for every
connected complex manifold, the holomorphic discs with dense
images form a dense subset of the set of all discs. To this aim a
nontrivial approximation statement is used and the result is the
following.

Let $M$ be a connected complex manifold, $d$ is a distance
generated by a complete Riemann metric, $A$ is a countable subset
of $M$, $f\in\CO(\Bbb D,X)$ and $r\in(0,1)$. Then there exists a
$g\in\CO(\Bbb D,X)$ such that $A\subset g(\Bbb D)$ and
$d(f(z),g(z))<1-r$ for each $z\in r\Bbb D.$

A modification of the proof of this fact shows that if, apart from
$A, f, r$, we are given a finite subset $\Lambda$ of $\Bbb D,$
then there exists a $g$ as above, as well as a sequence
$(\mu_\lambda)_{\lambda\in\Lambda}\subset\Bbb D$,
$r|\mu_\lambda|<|\lambda|$ such that $f(\lambda)=g(\mu_\lambda)$,
$\lambda\in\Lambda$. Letting $r\to 1$ one can prove that Theorem
\ref{decr.th1} remains true for complex manifolds.

\setcounter{equation}{0}
\section{Product property}\label{desc}

Let $l_D(\bs p,\cdot)$ and $l_G(\bs q,\cdot)$ be generalized
Lempert functions of domains $D\subset\C^n$ and $G\subset\C^m.$
They generate a generalized Lempert function $l_{D\times G}(\bs
r,\cdot)$ of the product product $G\times D,$ where
$$\bs r(\zeta,\eta)=\bs p(\eta)\bs q(\zeta),\quad \zeta\in D,\ \eta\in G.$$

In this section we discuss when does the generalized Lempert
function have a basic property, namely the product property, i.e.
$$l_{D\times G}(\bs r,(z,w))=\max\{l_D(\bs p,z),l_G(\bs q,w)\}.$$

Let us note that the Lempert functions, the Kobayashi functions,
and the Carath\'eodory functions have this property; a similar
property is true for their infinitesimal forms (see e.g.
\cite{Jar-Pfl1,Jar-Pfl3}).

We need the pluricomplex Green function $g_D$ defined as follows:
$$g_D(z,w)=\sup u(w),$$ where the supremum is over all negative
functions $u\in\PSH(D)$ such that
$u(\cdot)\le\log||\cdot-z||+O_u(1).$ Then
$$c^\ast_D\le g_D^\ast:=\exp g_D\le l_D$$
so for the infinitesimal form of $g_D,$ the so-called Azukawa
(pseudo)metric,
$$A_D(z;X)=\limsup_{\lambda\to 0}{g_D^\ast(z,z+\lambda X)\over|\lambda|},$$
we have
$$\gamma_D\le A_D\le\kappa_D.$$

For example, the theorem of Lempert implies that if $D$ is a
convex domain, then in both the chains of inequalities we have in
fact equalities.

Recall that $\Sp=\{a\in D:\bs p(a)>0\}.$ The next proposition
provides a necessary and sufficient condition for the product
property when the support of one of the functions is a singleton.

\begin{proposition}\label{one} If $(z,w)\in D\times G,$ $\Sp=\{a\}\subset D,$
then
$$l_{D\times G}(\bs r,(z,w))=\max\{l_D(\bs p,z),l_G(\bs q,w)\}$$ for each
function $q\gneq 0$ on $G$ if and only if $l_D(a,z)=g_D(a,z).$
\end{proposition}

A special case of Proposition \ref{one} was used is Section
\ref{lem-decr} with a quote of \cite[Theorem 2.1]{Dieu-Trao}. In
fact, this theorem is Proposition \ref{one} in the special case
when $\bs q$ is a characteristic function of a finite set; in the
general case the proof is similar and we omit it. We only note
that it is based on the inequality
$$l_{D\times G}(\bs r,(z,w))\ge\max\{l_D(\bs p,\cdot),
l_G(\bs q,\cdot)\},\quad \#\Sp=1.$$ The proof of this inequality,
given in \cite{Dieu-Trao} for the mentioned special case, contains
an essential flaw, corrected in the paper \cite{Nik-Zwo1} of the
author and W. Zwonek.
\smallskip

Similarly to the generalized Lempert function, one can define a
generalized Green function. This function does not exceed the
corresponding generalized Lempert function. In addition, by a
result of A. Edigarian, it possesses the product property (see
e.g. \cite{Edi0,Jar-Pfl3}).

On the other hand, D. Coman \cite{Com2} showed that the Lempert
and Green functions of a ball that have two poles coincide. He was
asking (see also \cite{Jar-Pfl3}) whether, like the Lemert
theorem, this property remains true for every convex domain for
every finite number of poles.

To give a negative answer of this question is one of the reasons
for our interest in the product property of the generalized
Lempert function (which shows that this function does not have so
typical properties as the generalized Green function). More
precisely, there are two-element subsets $A,B$ of $\D$ and a point
$z\in\D$ such that
\begin{equation}\label{tt}
l_{\D^2}(A\times B,(z,w))>\max\{l_\D(A,z),l_\D(B,0)\}.
\end{equation}
As
\begin{equation}\label{for}
l_\D(C,z)=g_\D(C,z)=\prod_{c\in C}m_\D(c,z),
\end{equation}
we get that
\begin{eqnarray*}l_{\D^2}(A\times
B,(z,0))&>&\max\{l_\D(A,z),l_\D(B,0)\}\\
&=&\max\{g_\D(A,z),g_\D(B,0)\}=g_{\D^2}(A\times
B,(z,0)).\end{eqnarray*}

The inequality (\ref{tt}) is first established in by P. J. Thomas
and N. V. Trao in \cite{Tho-Tra1} and independently, but somewhat
later, by the author and W. Zwonek in \cite{Nik-Zwo1}. In the
latter work the proof is considerably shorter and includes a
complete characterization of the two-element subsets $A$ and $B$
of $\D,$ for which we have the critical double equality
$$l_{\D^2}(A\times B,(z,w))=l_\D(A,z)=l_\D(B,w).$$
(This characterization shows that the product property is not
typical for the generalized Lempert function.) By applying an
automorphism of $\D^2,$ it suffices to consider only the case
$z=w=0.$

Note that, as above, $$l_{\D^2}(\bs r,(z,w))\ge g_{\D^2}(\bs r,(z,w))$$
$$=\max\{g_\D(\bs p,z),l_\D(\bs q,w)\}=\max\{l_\D(\bs p,z),l_\D(\bs q,w)\};$$
in particular, always
\begin{equation}\label{bi}
l_{\D^2}(A\times B,(z,w))\ge\max\{l_\D(A,z),l_\D(B,0)\}.
\end{equation}

\begin{proposition}\label{prod1} If $A=\{a_1,a_2\}\subset\D_\ast$
and $B=\{b_1,b_2\}\subset\D_\ast,$ then
\begin{equation}\label{tri}
l_{\D^2}(A\times B,(0,0))=l_\D(A,0)=l_\D(B,0)
\end{equation}
if and only if there us a rotation that maps $A$ in $B.$

In addition, if $B=e^{i\theta}A,$ $\theta\in\Bbb R,$ then the
extremal discs*\footnote{*As $\D^2$ is a taut domain, the infimum
in the definition of $l_{\D^2}$ is attained and the corresponding
discs are called extremal} for $l_{\Bbb D^2}(A\times B,(0,0))$ are
of the form
$\zeta\to(e^{i\varphi}\zeta,e^{i(\varphi+\theta)}\zeta),$
$\varphi\in\R.$
\end{proposition}

\noindent{\bf Remark.} From the last statement it follows that the
extremal discs for $l_{\Bbb D^2}(A\times B,(0,0))$ pass through
two points from the four-element set $A\times B$, although the
Lempert function decreases under adding of poles, according to
Corollary \ref{decr.cor1}.
\smallskip

\beginproof Let $\psi=(\psi_1,\psi_2)$ be an extremal disc for
$l_{\Bbb D^2}(A\times B,(0,0)).$
Then we can find a set $J\subset \{1,2\}\times\{1,2\}$
and points $z_{k,l}\in\Bbb D,$ $(k,l)\in J$ such that
$$\psi(z_{k,l})=(a_k,b_l)\hbox{ and } \prod_{(k,l)\in J}|z_{k,l}|=
l_{\Bbb D^2}(A\times B,(0,0)).$$

First let (\ref{tri}) be true. If $\#J=1,$ we can assume that
$J=\{(1,1)\}.$ Then
$$|z_{1,1}|=l_{\Bbb D^2}(A\times
B,(0,0))=l_{\D}(A,0)=|a_1a_2|<|a_1|=|\psi_1(z_{1,1})|\le|z_{1,1}|$$
(according to the Schwarz-Pick lemma), a contradiction.

If $\#J=3,$ we can assume that $J=\{(1,1),(1,2),(2,2\}.$ As above,
$$|z_{1,1}z_{1,2}z_{2,2}|=|a_1a_2|.$$ On the other hand, as
$\phi_1\in\O(\D,\D),\ \phi_1(0)=0,\
\phi_1(z_{1,1})=\phi_1(z_{1,2})=a_1,\ \phi_1(z_{2,2})=a_2,$ we get
$$|z_{1,1}z_{1,2}|\le|a_1|,\  |z_{2,2}|\le|a_2|,$$
with equalities attained when $\phi_1$ is a Blaschke product of
order 2 and a rotation, respectively -- a contradiction.

Let $\#J=4.$ We can assume that
$$\psi_1(z)=z\Phi_\alpha(z),\ \psi_2(z)=e^{it}z\Phi_\beta$$ for
some $\alpha,\beta\in\Bbb D,t\in\Bbb R$. Then
$$z_{1,1}\Phi_\alpha(z_{1,1})=z_{1,2}\Phi_\alpha(z_{1,2}),\
z_{2,1}\Phi_\alpha(z_{2,1})= z_{2,2}\Phi_\alpha(z_{2,2}),$$
$$z_{1,1}\Phi_\beta(z_{1,1})= z_{2,1}\Phi_\beta(z_{2,1}),\
z_{1,2}\Phi_\beta(z_{1,2})= z_{2,2}\Phi_\beta(z_{2,2}).$$
Consequently,
$$z_{1,1}=\Phi_\alpha(z_{1,2})=\Phi_\beta(z_{2,1}),\
z_{1,2}=\Phi_\beta(z_{2,2}),\ z_{2,1}=\Phi_\alpha(z_{2,2}).$$
Hence $z_{1,1}=\Phi_\alpha\circ\Phi_\beta(z_{2,2})=
\Phi_\beta\circ\Phi_\alpha(z_{2,2}).$ After some calculations we
get the equality
$$(2-\alpha\bar\beta-\bar\alpha\beta)(z_{2,2}^2(\bar\alpha-\bar\beta)
+z_{2,2}(\alpha\bar\beta-\bar\alpha\beta)+\beta-\alpha)=0.$$ It is
easily seen that if $\alpha\neq\beta,$ then the two roots of the
equation
$$z^2(\bar\alpha-\bar\beta)+z(\alpha\bar\beta-\bar\alpha\beta)=\alpha-\beta$$
lie on the unit circle. Therefore $\alpha=\beta,$
$z_{1,2}=z_{2,1},$ $z_{1,1}=z_{2,2}$, a contradiction.

Let now $\#J=2.$ We can assume that $J=\{(1,1),(2,2\}.$ Then it
easily follows that $\psi_1(z)=e^{i\theta_1}z,$
$\psi_2(z)=e^{i\theta_2}z,$ $\theta_1,\theta_2\in\Bbb R,$ and so
$B=e^{i\theta}A,$ where $\theta=\theta_1-\theta_2.$

Conversely, if $B=e^{i\theta}A,$ the mapping
$(\id,e^{i\theta}\id)\in\OO(\Bbb D,\Bbb D^2)$ is a competitor for
$l_{\Bbb D^2}(A\times B,(0,0))$ and so
$$l_{\Bbb D}(A,0)=l_{\Bbb D}(B,0)\geq\ l_{\Bbb D^2}(A\times B,(0,0)).$$
To get (\ref{tri}), it remains to use (\ref{bi}).\qed

\begin{corollary}\label{prod2} If $A$ and $B$ are two-point subsets
of $\Bbb D$ and $z\in\Bbb D\setminus A,$ then the set of points
$w\in\Bbb D$ such that $$l_{\Bbb D}(A,z)=l_{\Bbb D}(B,w)<l_{\Bbb
D^2}(A\times B,(z,w))$$ has Hausdorff dimension 1.
\end{corollary}

\beginproof  It suffices to note that the set of points $w\in\Bbb D$
such that $l_{\Bbb D}(A,z)=l_{\Bbb D}(B,w)$ has Hausdorff
dimension 1, and there are at most two points $w,$ for which there
is an automorphism of $\Bbb D$ that maps $z$ in $w$ and $A$ in
$B.$\qed
\medskip

We don't know whether Proposition \ref{prod1} remains true for
sets of equal cardinality, greater than 2. Anyway, for a given
point $(z,w)\in\Bbb D^2$ this proposition and (\ref{for}) provides
a large class of counterexamples for the product property of
$l_{\Bbb D^2}(A\times B,(z,w))$, where $A$ and $B$ have an
arbitrary number of elements, greater than 1.

\begin{proposition}\label{prod3} Let $z,w\in\D,$ $A,B\subset\D$
and $q\in(0,1)$ such that
$$\max\{l_\D(A,z),l_\D(B,w)\}=ql_{\D^2}(A\times B,(z,w))>0.$$
Then $$\max\{l_D(A\cup A_1,z),l_G(B\cup B_1,w)\}<l_{D\times
G}((A\cup A_1)\times(B\cup B_1),(z,w)),$$ if $A_1,B_1\subset\D,$
$A\cap A_1=B\cap B_1=\varnothing$ and $l_\D(A_1,z)l_\D(B_1,w)>q.$
\end{proposition}

\beginproof We have
$$l_{D\times G}((A\cup A_1)\times(B\cup B_1),(z,w))$$
$$\ge l_{D\times G}(A\times B,(z,w))l_{D\times G}(A\times B_1,(z,w))
l_{D\times G}(A_1\times (B\cup B_1),(z,w))$$
$$\ge l_{D\times G}(A\times B,(z,w))l_G(B_1,w)l_D(A_1,z)$$
$$>\max\{l_D(A,z),l_G(B,w)\}\ge\max\{l_D(A\cup A_1,z),l_G(B\cup B_1,w)\}$$
(the first inequality is checked immediately; for the second one
see (\ref{bi}); for the fourth one -- Corollary
\ref{decr.cor1}).\qed

\chapter{The symmetrized polydisc and the spectral ball}\label{chap.spec}
\setcounter{page}{38}
\section{Preliminaries}\label{np}

Most of the facts in section can be found in
\cite{Agl-You1,Agl-You3,Agl-You4,Bha,Cos2,Cos3,Cos4,Edi-Zwo,Jar-Pfl3,Pfl-Zwo}.

Let $\Bbb D\subset\C$ is the unit disc. Put
$\sigma=(\sigma_1,\dots,\sigma_n):\Bbb C^n\to\Bbb C^n,$ where
$$\sigma_k(z_1,\dots,z_n)=\sum_{1\le j_1<\dots<j_k\le
n}z_{j_1}\dots z_{j_k},\quad 1\le k\le n.$$ The open set $\Bbb
G_n=\sigma(\Bbb D^n)$ is called a symmetrized $n$-disc. Note that
$\Bbb G_n$ is a proper image of the $n$-disc $\Bbb D^n$ so it is a
pseudoconvex domain. Moreover, by \cite[Corollary 3.2]{Cos3} we
easily get that $\Bbb G_n$ is even a $c$-finite compact domain (in
particular hyperconvex), so it is a taut domain. Its Shilov
boundary is $\sigma(\Bbb T^n),$ where $\Bbb T=\partial\Bbb D$ is
the unit circle. Furthermore, the group of (holomorphic)
automorphisms of $\Bbb G_n$ admits a simple description:
$$\Aut(\Bbb G_n)=\{\sigma(h,\dots,h),\ h\in\Aut(\Bbb D)\}.$$
More general, a characterization of the proper holomorphic
mappings from $\G_n$ to itself can be found in \cite{Edi-Zwo}.

We also note that $\Bbb G_n$ is close to the balanced domains
(see Section \ref{bal}). More precisely,
$$\pi_{\lambda}(z)=(\lambda z_1,\lambda^2z_2,\dots,\lambda^n
z_n)\in\Bbb G_n,\ \lambda\in\overline{\Bbb D},z\in\Bbb G^n.$$

In fact $\Bbb G_n$ is the set of points
$(a_1,a_2,\dots,a_n)\in\Bbb C^n$ such that the zeroes of the
polynomial $f(\zeta)=\zeta^n+\sum_{j=1}^n (-1)^ja_j\zeta^{n-j},$
$a_0\neq 0,$ lie in $\Bbb D.$ Clearly $\Bbb G_1=\Bbb D.$
Furthermore, using the above description and the Cohn rule (see
Section \ref{bal}), we get that $$\Bbb
G_2=\{(s,p):|s-\overline{s}p|+|p|^2<1\}.$$

The symmetrized polydisc appears in connection with the so-called
Nevanlinna--Pick spectral problem.

Denote by $\mathcal M_n$ the set of $n\times n$ matrices of
complex coefficients. The spectral ball $\Om_n$ is defined by
$$\Om_n=\{A\in\mathcal M_n: r(A)=\max_{\lambda\in\spe(A)}|\lambda|<1\}$$
($r(A)$ and $\spe(A)$
are the spectral radius and the spectrum of $A$, respectively).

The spectral Nevanlinna--Pick problem, that we abbreviate by SNPP,
is the following:

Given $m$ different points $\lambda_1,\dots,\lambda_n\in\Bbb D$
and $m$ matrices $A_1,\dots,\\ A_m\in\Om_n$, determine whether
there exists a mapping $F\in\O(\Bbb D,\Om_n)$ that interpolates
the data, i.e. $F(\lambda_j)=A_j$ for $j=1,\dots,m.$

A nonconstructive necessary and sufficient condition for
solvability of SNPP is the solvability of the classical problem of
Nevanlinna--Pick for matrices that are similar to the given ones
(see e.g. \cite{BFT}). A more effective form of this result for
$2\times 2$ matrices can be found in \cite{Ber}. Now let us
describe an approach that reduces this problem of $n^2m$
parameters to a problem on $\Bbb G_n$ of $nm$ parameters.

We say that a matrix $A\in\mathcal M_n$ is cyclic if it has a
cyclic vector (i.e. $\CC^n=\span(v,Av,\dots,A^{n-1}v)$ for some
$v\in\CC^n$). In the appendix at the end of this section we
provide some equivalent conditions for a matrix to be cyclic. The
set of cyclic matrices in $\Omega_n$ will be denoted by $\mathcal
C_n.$

For $A\in\mathcal M_n$, put for brevity
$\sigma(A)=\sigma(\spe(A)).$ From the context it will be clear
whether we mean this $\sigma\in\O(\mathcal M_n,\C^n)$ or
$\sigma\in\O(\C^n,\C^n)$ (as defined at the beginning of this
section).

The following basic theorem for lifting a mapping from $\mathcal
O(\Bbb D,\Bbb G_n)$ to $\mathcal O(\Bbb D,\Omega_n)$ holds:

\begin{theorem}\label{np.th1} (see \cite{Agl-You1,Cos3})
Given $m$ different points $\lambda_1,\dots,\lambda_n\in\Bbb D$
and $m$ matrices $A_1,\dots,A_m\in\mathcal C_n.$ Let $f\in\mathcal
O(\Bbb D,\Bbb G_n)$ so that $f(\lambda_j)=\sigma(A_j)$ for
$j=1,\dots,m.$ Then there exists $F\in\mathcal O(\Bbb D,\Omega_n)$
so that $f=\sigma\circ F$ and $F(\lambda_j)=A_j$ for
$j=1,\dots,m.$
\end{theorem}

For arbitrary matrices $A_1,\dots,A_m\in\Om_n$ for $n\le 3$
the possibility for lifting a mapping is thoroughly discussed in \cite{NPT2}.

As $\mathcal C_n$ is a dense subset of $\Omega_n,$ this theorem
states that in the generic case SNPP is equivalent to an
interpolation problem on $\Bbb G_n$ (clearly one cannot expect a
similar result in full generality, as the spectrum does not
contain the full information for a given matrix up to a similarity
(for example for her Jordan or Frobenius form). As we noted, a
basic advantage of the second problem compared with the first one
is the smaller number of parameters. Furthermore, $\G_n$ is taut
domain, while on $\Omega_n$ one cannot apply the typical Montel
arguments, since $\Omega_n$ is not even Brody hyperbolic (it
contains complex lines). Probably the only advantage of $\Omega_n$
is that it is balanced domain; however this is compensated by the
previously mentioned fact that $\Bbb G_n$ is close to a balanced
domain.

The solution of SNPP is equivalent to finding the Lempert function
of $\Omega_n$. For cyclic matrices, Theorem \ref{np.th1} reduces
this question to finding the Lempert function of $\Bbb G_n.$
Indeed, the following simple proposition holds; we omit its proof.

\begin{proposition}\label{np.pr2}
Let $D\subset\C^n$ be a domain, $a_1,a_2\in D,$ $\lambda_1,\lambda_2\in\Bbb D.$

(i) If there exists $f\in\mathcal O(\Bbb D,D)$ so that
$f(\lambda_1)=a_1$ and $f(\lambda_2)=a_2,$ then $l_D(a_1,a_2)\le
m_{\Bbb D}(\lambda_1,\lambda_2).$

(ii) If $l_D(a_1,a_2)<m_{\Bbb D}(\lambda_1,\lambda_2),$ then there exists an $f$ as in (i).
\end{proposition}

As usual,
$m_{\Bbb D}(\lambda_1,\lambda_2)=\left|\frac{\lambda_1-\lambda_2}{1-\overline{\lambda_1}\lambda_2}\right|$
is the M\"obius distance.

Note that if for example $D$ is a taut domain, then there exist
extremal discs for $l_D$ and so the condition $l_D(a_1,a_2)\le
m_{\Bbb D}(\lambda_1,\lambda_2)$ is equivalent to the existence of
a corresponding $f.$ However, as mentioned, the spectral ball is
not such a domain. On the other hand, $\Bbb G_n$ is a taut domain.
Then Theorem \ref{np.th1} implies that
\begin{equation}\label{np2}
l_{\Omega_n}(A_1,A_2)=l_{\Bbb G_n}(\sigma(A_1),\sigma(A_2)),\quad
A_1,A_2\in\mathcal C_n
\end{equation}
(and there exists an extremal disc for $l_{\Omega_n}(A_1,A_2)$).
Note that as $\sigma\in\O(\Om_n,\G_n),$ in the general case
($A_1,A_2\in\Om_n$) we have the inequality $\ge.$

As the cyclic matrices form a dense subset of $\mathcal M_n$ and
the Kobayashi metric is continuous,
\begin{equation}\label{np5}k_{\Omega_n}(A_1,A_2)=k_{\Bbb
G_n}(\sigma(A_1),\sigma(A_2)),\quad A_1,A_2\in\Omega_n.
\end{equation}

The above considerations and the (almost) explicit finding (by the
method of geodesics) of $$l_{\Bbb G_2}(=c^\ast_{\Bbb G_2})$$ (see
the Introduction) permit a complete solution to SNPP for $n=2$
(see e.g. \cite{Cos4}). As the noncyclic $2\times 2$ matrices are
the scalar ones, for this purpose it remains just to calculate
$l_{\Omega_n}(\lambda I_n,\cdot)$ for $n=2$ ($I_n\in\M_n$ is the
unit matrix). As $\Omega_n$ is a pseudoconvex balanced domain
whose Minkowski function is the spectral radius $r,$ one has
$l_{\Omega_n}(0,\cdot)=r$ (see Proposition \ref{bala.pr1} (iv)).
Then from
\begin{equation}\label{fi}
\Phi_\lambda(A)=(A-\lambda I_n)(I_n-\overline{\lambda}A)^{-1}\in\Aut(\Omega_n),
\end{equation}
we deduce
\begin{equation}\label{np4}
l_{\Omega_n}(\lambda I_n,A)=r(\Phi_\lambda(A))=
\max_{a\in\spe(A)}m_{\Bbb D}(\lambda,a).
\end{equation}
Hence it also follows that $\Omega_2$ is an example of a
nonhyperbolic pseudoconvex balanced domain such that
$c^\ast_{\Omega_2}=k^\ast_{\Omega_2}\lneq l_{\Omega_2};$ on the
other hand, $k^\ast_{\Omega_2}(A_1,A_2)=l_{\Omega_2}(A_1,A_2)$ for
$A_1,A_2\in\mathcal C_2.$

In connection with the use of $\Phi_\lambda$ let us note that in
\cite{Ran-Whi} there was a hypothesis for the complete description
of $\Aut(\Omega_n).$ This hypothesis has been recently disproved
in \cite{Kos}. Note that for the Euclidean ball $\Bbb B_n,$ each
proper holomorphic mapping from $\Omega_n$ into itself is an
automorphism (see \cite{Zwo2}).

The approach of complex geodesics is applied in \cite{Cos3} and
for some special pairs of points from $\Bbb G_n$ for $n\ge 3.$ In
other words, the Lempert function for all these pairs of points
coincides with the Carath\'eodory function. However in Section
\ref{ine} we will see that this is not true for each pair of
points by obtaining some inequalities for the Carath\'eodory and
Kobayashi metrics on $\Bbb G_n,$ taken at the beginning. These
inequalities and the lower estimates from Section \ref{est} (for
the Carath\'eodory metric, and hence for the Kobayashi metric)
carry some information for the so-called spectral problem of
Carath\'eodory--Fej\'er. A reduction of this problem to a
corresponding problem on $\G_n$ in the spirit of Theorem
\ref{np.th1} can to be found in \cite[Theorem 2.1]{HMY}.

The easiest variant of this problem, that we abbreviate by SCFP,
is the following one:

For $A\in\Om_n$ and $B\in\M_n$, determine whether there exists a
mapping $F\in\O(\Bbb D,\Om_n)$ such that $F(0)=A$ and $F'(0)=B.$

In \cite{NPT2} SCFP is completely reduced to a problem on $\G_n$
for $n\ge 3.$

As SNPP, SCFP is also connected with finding of $\kappa_{\Om_n}.$
Similarly to (\ref{np2}), we have
\begin{equation}\label{np1}
\kappa_{\Om_n}(A;B)=\kappa_{\G_n}(\sigma(A),\sigma'_A(B)), \quad
A\in\mathcal C_n, B\in\mathcal M_n,
\end{equation}
where in this case $\sigma'_A=\sigma_{\ast,A}$ is the Fr\'echet
derivative of $\sigma$ at $A$. Furthermore,
\begin{equation}\label{np3}
\kappa_{\Omega_n}(\lambda I_n;B)=\frac{r(B)}{1-|\lambda|^2},
\end{equation}
which together with \cite[Theorem 1.1]{HMY} permits a complete
solution of SCFP for $n=2.$

We conclude the section with the fact that the Carath\'eodory
metric and Carath\'eodory distance on $\Omega_n$ can be calculated
via these on $\Bbb G_n$ (cf. (\ref{np5})).

\begin{proposition}\label{np.pr5} The following equalities hold:
$$c_{\Omega_n}(A_1,A_2)=c_{\Bbb
G_n}(\sigma(A_1),\sigma(A_2)),\ A_1,A_2\in\Omega_n,$$
$$\gamma_{\Omega_n}(A;B)=\gamma_{\Bbb G_n}(\sigma(A);\sigma'_A(B)),
\ A\in\Omega_n,B\in\mathcal
 M_n.$$
\end{proposition}

\beginproof As $\sigma\in\mathcal O(\Bbb G_n,\Bbb D),$ we have the
inequalities $\ge.$ For the reverse inequalities it suffices to
show that if $f\in\O(\Om_n,\D),$ then there exists a
$g\in\O(\G_n,\D)$ such that $f=g\circ\sigma.$ First note that if
$A,B\in\Om_n$ have identical spectra, then there exists an entire
curve $\varphi$ in $\Om_n$, passing through $A$ and $B$ (see
Proposition \ref{i.pr4} (ii)). By the Liouville Theorem (applied
for the function $f\circ\varphi$) implies $f(A)=f(B)$ and so $g$
is a well defined function. It is holomorphic, since for each
layer $\sigma^{-1}(\sigma(C))$ there is a matrix $\tilde C$ so
that $\rank\sigma'_{\tilde C}=n$ (see Proposition
\ref{cyc.equiv}).\qed

\setcounter{equation}{0}
\section{Cyclic matrices}\label{cyc}

In this brief section we provide some equivalent conditions for a
matrix to be cyclic. Part of these will be used at the end of the
chapter.

First let us recall some definitions.

Let $A\in\mathcal M_n,$ $\ad_A:X\to [A,X]$ is the adjoint mapping
of $A$ and $\mathcal C_A=\ker\ad_A$ is the centralizer of $A.$

Let $p_A(x)=x^n+a_{n-1}x^{n-1}+\dots+a_0$ is the characteristic
polynomial of $A.$ The matrix
$$
\left(
\begin{array}{ccccc}
0&0&\dots&0&-a_0\\
1&0&\dots&0&-a_1\\
0&1&\dots&0&-a_2\\
0&0&\dots&1&-a_{n-1}
\end{array}
\right).
$$
is called adjoint to $A$ (or to $p$).

\begin{proposition}\label{cyc.equiv}
For a matrix $A\in\mathcal M_n$ the following are equivalent:

\begin{enumerate}

\item $A$ has a cyclic vector.

\item $A$ is similar to its adjoint matrix (i.e. it is the Frobenius form of $A$).

\item The characteristic polynomial and the minimal polynomial of $A$ coincide.

\item Different blocks in the Jordan form of $A$ correspond to
different eigenvalues (i.e. each eigenspace is one-dimensional).

\item $\mathcal C(A)=\{M\in\mathcal M_n :M=p(A)\mbox{ for some }p\in \mathbb C[X]\}$.

\item $\dim\mathcal C_A=n$.

\item $\rank\sigma'_A=n$.

\item $\ker\sigma'_A=\Ima\ad_A.$

\end{enumerate}
\end{proposition}

A matrix with (one of) the above properties is called cyclic.
\smallskip

\beginproof  The equivalence of the properties from (1) to (6)
are well-known and can be found e.g. in \cite{Ho-Jo1,Ho-Jo2}. We
need to prove their equivalence with (7) and (8).

Observe that if $M\in\mathcal M_n^{-1}$ (i.e. $Ì$ is an invertible
matrix), then $$\sigma'_A(X)=\sigma'_{M^{-1}AM}(M^{-1}XM).$$ So to
prove that (2) implies (7), we can assume that $A$ coincides with
its adjoint matrix. Let $X=(x_{i,j})$ so that $x_{i,j}=0$ for
$1\le j\le n-1$. Then
$$\sigma'_A(X)=(-x_{n,n},x_{n-1,n},\ldots,(-1)^{n-1}x_{1,n})$$ and
consequently $\Ima\sigma'_A=\CC^n$, i.e. $\rank\sigma'_A=n$.

Let us now prove that (7) implies (4). Let $\lambda\in\Bbb C,$ and
$M_\lambda=M-\lambda I_n.$ As $p_M(x)=p_{M_\lambda}(x+\lambda)$,
there exists $\Lambda(\lambda)\in\mathcal M_n^{-1}$ so that
$\sigma(M)=\Lambda(\lambda)\sigma(M_\lambda).$ Therefore
$\rank\sigma'_A=\rank\sigma'_{A_\lambda}$.

Suppose that (7) is true and (4) is false. There exists an
eigenvalue $\lambda$ of $A$ such that $\dim\ker(A-\lambda I_n)\ge
2$. Let us complete a basis of $\ker(A-\lambda I_n)$ to a basis of
$\CC^n.$ Then the matrix $A-\lambda I_n$ is transformed to a
matrix with at least two zero columns and consequently
$\sigma_{n,n}(A-\lambda I_n+X)$ is a polynomial of degree two or
more with respect to $x_{i,j}.$ Hence
$(\sigma_{n,n})_{\ast,A_\lambda}=0$ and so
$\rank\sigma'_A=\rank\sigma'_{A_\lambda} \le n-1,$ a
contradiction.

At last let us show that $(6)+(7)\Leftrightarrow (8).$ It is
easily seen that $\Ima\ad_A\subset\ker\sigma'_A.$ Consequently
$\Ima\ad_A=\ker\sigma'_A$ if and only if these two linear spaces
have the same dimension. By the Rank theorem, this is equivalent
to $\dim\mathcal C_A=\rank\sigma'_A$. It remains to use that
$\dim\mathcal C_A\ge n\ge\rank\sigma'_A$ for each $A\in\mathcal
M_n.$\qed

\setcounter{equation}{0}
\section{$\Bbb G_n$ is not a Lu Qi-Keng domain for $n\ge 3$}\label{lqk}

In 1966 Lu Qi-Keng \cite{LQK} formulated the hypothesis that the
Bergman kernel (see Section \ref{ber} for the definition) of a
simply connected domain in $\Bbb C^n$ has no zeroes.

A domain of this property is called a Lu Qi-Keng domain. This
hypothesis is rejected in 1986 by H. P. Boas \cite{Boas1}. A
review of the role of Lu Qi-Keng domains in complex analysis,
together with various counterexamples, can be found for example in
\cite{Boas2,Jar-Pfl3}.

Using the Bell's transformation formula, (see e.g.
\cite{Jar-Pfl1}), in \cite{Edi-Zwo} the authors find an explicit
formula for the Bergman kernel of the symmetrized polydisc. This
formula implies that $\Bbb G_2$ is a Lu Qi-Keng domain.

The aim of this section is to provide an affirmative answer of
\smallskip

\noindent{\bf\cite[Problem 3.2]{Jar-Pfl3}.} {\it Does the Bergman
kernel of $\Bbb G_n$ have zeroes for $n\ge 3$?}
\smallskip

Thus $\G_n$ is the first example of a proper image of the polydisc
$\Bbb D^n,$ $n\ge 3$ that is not a Lu Qi-Keng domain, once again
showing the difference in the structure of $\Bbb G_n$ for $n=2$
and $n\ge 3.$

Now let us state the formula for Bergman kernel of $\Bbb G_n$
\cite{Edi-Zwo}:
\begin{equation}\label{l1} K_{\Bbb
G_n}(\sigma(\lambda),\sigma(\mu))=
\frac{\det[(1-\lambda_j\overline\mu_k)^{-2}]_{1\le j,k\le n}}
{\pi^n\prod_{1\le j<k\le
n}[(\lambda_j-\lambda_k)(\overline\mu_j-\overline\mu_k)]},\ \
\lambda,\mu\in\Bbb D^n.
\end{equation}
Although formally the right-hand side of (\ref{l1}) is not defined
on the whole $\Bbb G_n,$ it is continued smoothly there. From this
formula we get

\begin{proposition}\label{l.pr1} \cite[Proposition 11]{Edi-Zwo}
$\Bbb G_2$ is a Lu Qi-Keng domain.
\end{proposition}

\beginproof From (\ref{l1}) it is easily deduced that
$$K_{\Bbb G_2}(\sigma(\lambda),\sigma(\mu))=\frac
{2-(\lambda_1+\lambda_2)(\overline\mu_1+\overline\mu_2)+
2\lambda_1\lambda_2\overline\mu_1\overline\mu_2}
{\pi^2\prod_{j,k=1}^2(1-\lambda_j\overline\mu_k)^2}.$$ Consider an
automorphism of $\Bbb D$ that maps $\mu_2$ to 0. It clearly
defines an automorphism of $\Bbb G_2$ (see the beginning of
Section \ref{np}). Consequently it suffices to show that $K_{\Bbb
G_2}(\sigma(\lambda),\sigma(\mu))\neq 0$ if $\mu_2=0.$ This is
trivial since
$|(\lambda_1+\lambda_2)\overline\mu_1|<2.$\qed\smallskip

Unlike Proposition \ref{l.pr1}, we have

\begin{theorem}\label{l.th2}
$\Bbb G_n$ is not a Lu Qi-Keng domain for $n\ge 3.$
\end{theorem}

\beginproof
We prove by induction on $n\ge 3$ that:

$(\ast)$ there exist points $\lambda,\mu\in\Bbb D^n$ with pairwise
different coordinates such that
$$\Delta_n(\lambda,\mu):=\det[(1-\lambda_j\overline\mu_k)^{-2}]
_{1\le j,k\le n}=0$$ and $f_n=\Delta_n(\cdot,\lambda_2,\dots,\lambda_n,\mu_1,\dots,\mu_n)\not\equiv
0.$

{\it Base of the induction: $n=3.$} We use the following formula
(see Appendix A):
\begin{equation}\label{l2} K_{\Bbb
G_3}(\sigma(\lambda_1,\lambda_2,\lambda_3),\sigma(\mu_1,\mu_2,0))=
\frac{a(\nu)z^2-b(\nu)z+2c(\nu)}{\pi^3\prod_{1\le j\le 3,1\le k\le
2}(1-\lambda_j\overline\mu_k)^2},
\end{equation}
where $z={\overline\mu_2}/{\overline\mu_1}$ $(\mu_1\neq 0),$
$\nu_j=\lambda_j\overline\mu_1,$ $j=1,2,3,$ and
$$a(\nu)=\sigma_2(\nu)(2-\sigma_1(\nu))+\sigma_3(\nu)(2\sigma_1(\nu)-3),$$
$$b(\nu)=(\sigma_1(\nu)-2)(\sigma_2(\nu)-2\sigma_1(\nu)+3)+
3(\sigma_3(\nu)-\sigma_1(\nu)+2),$$
$$c(\nu)=\sigma_2(\nu)-2\sigma_1(\nu)+3.$$ For the fixed
point $\nu_0=(e^{i\sigma/6},e^{i\sigma/3},e^{-i\sigma/6})$ the number
$$z_0=e^{-i\sigma/4}\frac{6-3\sqrt3-\sqrt{40\sqrt3-69}}{\sqrt2(3\sqrt3-5)}$$
is a root of the equation $a(\nu_0)z_0^2-b(\nu_0)z_0+2c(\nu_0)=0$ (see
Appendix B). As $z_0\in\Bbb D,$ for $\nu\in\Bbb
D^3,$ close to $\nu_0,$ there exists $z\in\Bbb D,$ close to $z_0,$
such that $a(\nu)z^2-b(\nu)z+2c(\nu)=0.$ Now choosing $\mu_1\in\Bbb
D$ such that $|\mu_1|>|\nu_1|,|\nu_2|,|\nu_3|,$ we get points
$\lambda,\mu\in\Bbb D^3$ with pairwise different coordinates so that
$\Delta_3(\lambda,\mu)=0.$

It remains to check that $f_3\not\equiv 0.$ If this fails, then
$f_3(0)=f_3'(0)=f_3''(0)=0,$ i.e.
$$\det\left[\begin{matrix}
\overline\mu_1^j&\overline\mu_2^j&\overline\mu_3^j\\
(1-\lambda_2\overline\mu_1)^{-2}&(1-\lambda_2\overline\mu_2)^{-2}&(1-\lambda_2\overline\mu_2)^{-2}\\
(1-\lambda_3\overline\mu_1)^{-2}&(1-\lambda_3\overline\mu_2)^{-2}&(1-\lambda_3\overline\mu_3)^{-2}\\
\end{matrix}\right]=0$$
for $j=0,1,2.$ As $\mu_1,\mu_2,\mu_3$ are pairwise different, the
vectors $(1,1,1),$ $(\mu_1,\mu_2,\mu_3)$ and
$(\mu_1^2,\mu_2^2,\mu_3^2)$ are $\Bbb C$-linearly independent.
Consequently the vectors in the second and third row of the above
determinant are $\Bbb C$-linearly dependent. In particular,
$K_{\Bbb G_2}(\sigma(\lambda_2,\lambda_3),\sigma(\mu_2,\mu_3))=0,$
a contradiction.

{\it Induction Step.} Suppose that $(\ast)$ holds for some $n\ge
3.$ Choose points $\tilde\lambda_1$ and $\tilde\lambda_{n+1}$ in
$\Bbb D,$ close to $\lambda_1$ and 1, respectively (this
guarantees that the coordinates of the new points in $\Bbb
D^{n+1}$ are pairwise different), so that
$$g_{n+1}(\tilde\lambda_1,\tilde\lambda_{n+1}):=
\Delta_{n+1}(\tilde\lambda_1,\lambda_2,\dots,\lambda_n,\tilde\lambda_{n+1},
\mu_1,\dots,\mu_n,\tilde\lambda_{n+1})=0$$ and
$g_{n+1}(\cdot,\lambda_{n+1})\not\equiv 0.$ Note that
$$g_{n+1}(\tilde\lambda_1,\tilde\lambda_{n+1})=\frac{f_n(\tilde\lambda_1)}
{(1-|\tilde\lambda_{n+1}|^2)^2}+h_n(\tilde\lambda_1,\tilde\lambda_{n+1}),$$
where $h_n$ is a continuous function on $\Bbb
D\times\overline{\Bbb D}.$ As $f_n\not\equiv 0$ is a holomorphic
function, for each small $r>0$ the number $\lambda_1$ is the only
zero of $f_n$ in the closed disc $D\subset\Bbb D$ of center
$\lambda_1$ and radius $r.$ Then $m=\frac{\min_{\partial
D}|f_n|}{\max_{\partial D \times\overline{\Bbb D}}|h_n|}>0.$
Consequently
$|f_n|>(1-|\tilde\lambda_{n+1}|^2)^2|h_n(\cdot,\tilde\lambda_{n+1})|$
on $\partial D$, provided $1-|\tilde\lambda_{n+1}|^2<\sqrt m.$ Fix
one such $\tilde\lambda_{n+1}$ so that
$\tilde\lambda_{n+1}\neq\lambda_j,\mu_j,$ $1\le j\le n.$ As
$h_n(\cdot,\tilde\lambda_{n+1})$ is a holomorphic function on
$\Bbb D$, by Rouch\'e's theorem the number of zeroes of
$g_{n+1}(\cdot,\tilde\lambda_{n+1})$ in $D$ equals the
multiplicity of $\lambda_1$ as a zero of $f_n;$ in particular,
$g_{n+1}(\cdot,\tilde\lambda_{n+1})\not\equiv 0$. It remains to
choose $r$ so that $\lambda_j,\mu_j,\tilde\lambda_{n+1}\not\in D,$
$1\le j\le n,$ and a zero $\tilde\lambda_1$ of
$g_{n+1}(\cdot,\tilde\lambda_{n+1})$ in $D.$\qed
\smallskip

\noindent{\bf Remark.} The above proof shows that if $n\ge 4,$
then there exist points $(\lambda,\nu),$ close to the diagonal of
$\Bbb D^n\times\Bbb D^n$ in the sense that $\lambda_j=\mu_j>0$ for
$j=4,\dots,n$ and such that $K_{\Bbb
G_n}(\sigma(\lambda),\sigma(\mu))=0.$ On the other hand, one can
prove that $K_{\Bbb G_3}(\sigma(\lambda),\sigma(\mu))\neq 0$ if
$\lambda_3=\mu_3.$
\smallskip

\noindent{\bf Appendix A.} (\ref{l1}) implies that
\begin{equation}\label{l3}
\pi^3(\lambda_1-\lambda_2)(\lambda_1-\lambda_3)(\lambda_2-\lambda_3)\overline\mu_1
\overline\mu_2(\overline\mu_1-\overline\mu_2)K_{\Bbb
G_3}(\sigma(\lambda_1,\lambda_1,\lambda_3),\sigma(\mu_1,\mu_2,0))
\end{equation}
$$=\det\left[\begin{matrix} (1-\nu_1)^{-2}&(1-z\nu_1)^{-2}&1\\
(1-\nu_2)^{-2}&(1-z\nu_2)^{-2}&1\\
(1-\nu_3)^{-2}&(1-z\nu_3)^{-2}&1
\end{matrix}\right]$$
$$=\det\left[\begin{matrix}
(1-\nu_1)^{-2}-(1-\nu_3)^{-2}&(1-z\nu_1)^{-2}-(1-z\nu_3)^{-2}\\
(1-\nu_2)^{-2}-(1-\nu_3)^{-2}&(1-z\nu_2)^{-2}-(1-z\nu_3)^{-2}
\end{matrix}\right]$$
$$=\frac{(\nu_1-\nu_3)(\nu_2-\nu_3)z}{(1-\nu_3)^2(1-z\nu_3)^2}
\det\left[\begin{matrix}
\frac{\nu_1+\nu_3-2}{(1-\nu_1)^2}&\frac{z\nu_1+z\nu_3-2}{(1-z\nu_1)^2}\\
\frac{\nu_2+\nu_3-2}{(1-\nu_2)^2}&\frac{z\nu_2+z\nu_3-2}{(1-z\nu_2)^2}
\end{matrix}\right]$$
$$=\frac{(\nu_1-\nu_3)(\nu_2-\nu_3)z}{\prod_{1\le j\le 3,1\le k\le
2}(1-\lambda_j\overline\mu_k)^2}
\Bigl((\nu_1+\nu_3-2)(z\nu_2+z\nu_3-2)(1-z\nu_1)^2(1-\nu_2)^2$$
\begin{equation}\label{l4}
-(\nu_2+\nu_3-2)(z\nu_1+z\nu_3-2)(1-\nu_1)^2(1-z\nu_2)^2\Bigr)
\end{equation}
\begin{equation}\label{l5}
=\frac{(\nu_1-\nu_3)(\nu_2-\nu_3)z(z-1)(A(\nu)z^2-B(\nu)z+2C(\nu))}{\prod_{1\le
j\le 3,1\le k\le 2}(1-\lambda_j\overline\mu_k)^2}.
\end{equation}
To find $A(\nu)$, $B(\nu)$ and $C(\nu)$, we use that the
coefficients of $z^3,$ $z^0$ and $z$ in the large parentheses of
(\ref{l4}) are equal to
$$A(\nu)=(\nu_1+\nu_3-2)(\nu_2+\nu_3)\nu_1^2(1-\nu_2)^2-(\nu_2+\nu_3-2)(\nu_1+\nu_3)\nu_2^2(1-\nu_1)^2,$$
$$-2C(\nu)=2(\nu_2+\nu_3-2)(1-\nu_1)^2-2(\nu_1+\nu_3-2)(1-\nu_2)^2\hbox{ and}$$
$$B(\nu)+2C(\nu)=(\nu_1+\nu_3-2)(\nu_2+\nu_3+4\nu_1)(1-\nu_2)^2$$
$$-(\nu_2+\nu_3-2)(\nu_1+\nu_3+4\nu_2)(1-\nu_1)^2,$$ respectively.
Some trivial calculations show that
$$A(\nu)=(\nu_2-\nu_1)(\sigma_{3,2}(\nu)(2-\sigma_{3,1}(\nu))+\sigma_{3,3}(\nu)
(2\sigma_1(\nu)-3)),$$
$$C(\nu)=(\nu_2-\nu_1)(\sigma_2(\nu)-2\sigma_1(\nu)+3),$$
$$B(\nu)=(\nu_2-\nu_1)((\sigma_1(\nu)-2)(\sigma_2(\nu)-2\sigma_1(\nu)+3)
+3(\sigma_3(\nu)-\sigma_1(\nu)+2)).$$ To infer (\ref{l2}), it
remains to substitute these formulas in (\ref{l5}) and then to
compare (\ref{l5}) and (\ref{l3}).
\smallskip

\noindent{\bf Appendix B.} As
$$\sigma_1(\nu_0)=\frac{1+2\sqrt3+i\sqrt3}{2},\
\sigma_2(\nu_0)=\frac{2+\sqrt3+i3}{2},\
\sigma_3(\nu_0)=e^{i\sigma/3},$$ the formulas for $a(\nu),$
$b(\nu)$ and $c(\nu)$ lead to
$$a(\nu_0)=(3\sqrt3-5)e^{i\sigma/3},\ b(\nu_0)=
(6\sqrt2-3\sqrt6)e^{i\sigma/12},\
c(\nu_0)=(2\sqrt3-3)e^{-i\sigma/6}.$$ Then for $z=e^{-i\sigma/4}x$
we have $e^{i\sigma/6}(a(\nu_0)z^2-b(\nu_0)z+2c(\nu_0))$
$$=(3\sqrt3-5)x^2+(3\sqrt6-6\sqrt2)x+4\sqrt3-6=:p(x).$$ It remains to note
that the zeroes of the polynomial $p$ are equal to
$\frac{6-3\sqrt3\pm\sqrt{40\sqrt3-69}}{\sqrt2(3\sqrt3-5)}$ and the
smaller one lies in $(0,1).$

\setcounter{equation}{0}
\section{Generalized balanced domains}\label{bal}

To show that $\Bbb G_n\not\in\mathcal E$ for $n\ge 2$ (see the
Introduction), we will define the so-called generalized balanced
domains. For such domains we will find a necessary condition for
belonging to $\mathcal E$ and then we will show that $\Bbb G_n,$
$n\ge 3,$ does not satisfy this condition; for $\Bbb G_2$ the
proof is somewhat different.

Let $k_1\le k_2\le\dots\le k_n$ be natural numbers and
$$\pi_\lambda(z_1,\dots,z_n)=(\lambda^{k_1}z_1,\dots,\lambda^{k_n}z_n),\
\lambda\in\Bbb C,z\in\Bbb C^n.$$ A domain $D$ in $\Bbb C^n$ will
be called $(k_1,\dots,k_n)$-balanced or generalized balanced, if
$\pi_{\lambda}(z)\in D$ for each $\lambda\in\overline{\Bbb D},$
$z\in D.$ Put
$$h_D(z)=\inf\{t>0:\pi_{1/t}(z)\in D\},\ z\in\Bbb C^n$$
(generalized Minkowski function). The function $h_D$ is
nonnegative and upper semicontinuous,
$$h_D(\pi_\lambda(z))=|\lambda|h_D(z),\ \lambda\in\Bbb C,z\in\Bbb
C^n,$$ $$D=\{z\in\Bbb C^n:h_D(z)<1\}.$$

\noindent{\bf Example.} $ h_{\Bbb G_n}(\sigma(\xi_1,\dots,\xi_n))=
\max_{1\le j\le n}|\xi_j|.$
\smallskip

Clearly the $(1,\dots,1)$-balanced domains are exactly the usual
balanced domains. Part of their properties remain true for the
generalized balanced domains.

\begin{proposition}\label{b.pr2} Let $D$ be a generalized balanced domain.
Then $D$ is pseudoconvex exactly when $\log h_D\in\PSH(\Bbb C^n).$
Furthermore, the following are equivalent:

(i) $\log h_D\in\PSH(\Bbb C^n)\cap C(\Bbb C^n)$ and
$h_D^{-1}(0)=\{0\}$ (i.e. $D$ is a bounded domain);

(ii) $D$ is a hyperconvex domain;

(iii) $D$ is a taut domain.
\end{proposition}

\beginproof Clearly if $\log h_D\in\PSH(\Bbb C^n),$ then $D$ is pseudoconvex.

To prove the reverse implication, let $D$ be
$(k_1,\dots,k_n)$-balanced. Put $\Phi:\Bbb
C^n\owns(z_1,\dots,z_n)\mapsto (z_1^{k_1},\ldots,z_n^{k_n})\in\Bbb
C^n,$ $\tilde D=\Phi^{-1}(D)$ and $\tilde h_D=h_D\circ\Phi$. Then
$\tilde D=\{z\in\Bbb C^n:\tilde h_D(z)<1\}$ and $\tilde
h_D(\lambda z)=|\lambda|h(z)$, $\lambda\in\Bbb C$, $z\in\Bbb C^n$.
Consequently $\tilde D$ is a pseudoconvex balanced domain of
Minkowski functional $\tilde h_D$. So $\log\tilde h_D\in\PSH(\Bbb
C^n).$ On the other hand, $h_D(z)=\tilde
h_D(\root{k_1}\of{z_1},\ldots,\root{k_n}\of{z_n})$, $z\in(\Bbb
C^n)_\ast$, where the roots are arbitrarily chosen. Consequently
$\log h_D\in\PSH((\Bbb C^n)_\ast)$ By the Removable Singularities
Theorem (see e.g. \cite{Jar-Pfl1}) we conclude that $\log
h_D\in\PSH(\Bbb C^n).$

Now observe that the implication (ii)$\Rightarrow$(iii) is true
for an arbitrary domain, while (i)$\Rightarrow$(ii) is trivial,
since $\log h$ is a negative exhausting plurisubharmonic function
for $D.$ To prove (iii)$\Rightarrow$(i), we first note that $D$ is
a pseudoconvex domain and so $\log h_D\in\PSH(\Bbb C^n)$ according
to the first part of the proposition. Furthermore, if
$h_D^{-1}(0)\neq\{ 0\},$ then $h_D(z)=0$ for some $z.$ Then $D$
contains the entire curve $\Bbb
C\ni\lambda\to\pi_{\lambda}(z)\in\Bbb C^n$ and so $D$ is not even
Brody hyperbolic, a contradiction. Now suppose that $h_D$ is not
continuous. As $h_D$ is upper semicontinuous, one can find
$\varepsilon>0$ and a sequence of points $z_j$ tending to some $z$
so that $h_D(z_j)<h_D(z)-\varepsilon.$ By homogeneity of $h_D$, we
can assume that $h_D(z_j)<1<h_D(z)$ for each $j.$ Then the
holomorphic discs $\Bbb D\ni\lambda\to \pi_{\lambda}(z_j)\in D$
converge (locally uniformly) to the disc $\Bbb D\ni\lambda\to
\pi_{\lambda}(z)\in\Bbb C^n$ that is not lying completely within
$D,$ a contradiction. This proves the implication
(iii)$\Rightarrow$(i).\qed\smallskip

\noindent{\bf Remarks.} a) The above proof shows that a
generalized balanced domain is hyperbolic if exactly when it is
bounded.
\smallskip

b) In the case of a balanced domain, the implication
(iii)\\$\Rightarrow$(i) can be proven also like this. As $D$ is a
taut domain, $\gamma_D$ is a continuous function. It remains to
observe that a taut domain is hyperbolic and pseudoconvex; so
$\gamma^{-1}_D(z;\cdot)=\{0\}$ for each $z\in D$ and
$\gamma_D=h_D.$
\smallskip

The next theorem provides a necessary condition for a generalized
balanced domain in $\Bbb C^n$ to belong to the class $\mathcal E$
in terms of convexity of its intersection with a linear subspace
of $\Bbb C^n.$

\begin{theorem}\label{b.th3} Let $D\in\mathcal E$ be a
$(k_1,\dots,k_n)$-balanced domain in $\Bbb C^n.$ If $2k_{m+1}>k_n$
for some $m,$ $0\le m\le n-1,$ then the intersection
$D_m=D\cap\{z_1=\dots=z_m=0\}$ is a convex set (we put $D_m=D$ if
$m=0$).
\end{theorem}

\beginproof The proof is similar to that of \cite[Theorem 1]{Edi}.

Fix $a,b\in D_m.$ Then we can find a domain $D'\subset D$ that is
biholomorphic to a convex domain $G$, and such that $\lambda
a,\lambda b\in D'$ for $\lambda\in\overline{\Bbb D}.$ Let
$\Psi:D'\to G$ be the corresponding biholomorphic mapping. After a
linear coordinate substitution we can assume that $\Psi(0)=0$ and
$\Psi'(0)=\id.$ Put
$$g_{ab}(\lambda)=\frac{\Psi(\pi_\lambda(a))+\Psi(\pi_\lambda(b))}{2}.$$
Then $\Psi^{-1}\circ g_{ab}(\lambda)$ is a holomorphic mapping
from a neighborhood of $\overline{\Bbb D}$ in $D.$ Let
$f_{ab}(\lambda)=\pi_{1/\lambda}\circ\Psi^{-1}\circ
g_{ab}(\lambda).$ We will show that
\begin{equation}\label{b1}
\lim_{\lambda\to 0}f_{ab}(\lambda)=\frac{a+b}{2}.
\end{equation}
Then $f_{ab}(\lambda)$ is extended analytically at $0.$
Consequently $h\circ f_{ab}\in\PSH(\overline{\Bbb D})$ according
to Proposition \ref{b.pr2} and the maximum principle shows that
$$h(f_{ab}(0))\le\max_{|\lambda|=1}h(f_{ab}(\lambda))<1.$$
Hence $\frac{a+b}{2}\in D_m$ for $a,b\in D_m,$ i.e. $D_m$ is a convex set.

To prove (\ref{b1}), note that the equalities $\Psi^{-1}(0)=0$ and
$(\Psi^{-1})'(0)=\id$ imply
$$\Psi_j^{-1}\circ g_{ab}(\lambda)=g_{abj}(\lambda)+O(|g_{ab}(\lambda)|^2),\
j=1,2,\dots,n.$$ As $\Psi(0)=0,$ $\Psi'(0)=\id$ and $a,b\in D_m,$ we get
$$g_{abj}(\lambda)=\frac{a_j+b_j}{2}\lambda^{k_j}+O(|\lambda|^{2k_{m+1}}).$$
The inequalities $2k_{m+1}>k_n$ show that
$$\frac{\Psi_j^{-1}\circ g_{ab}(\lambda)}{\lambda^{k_j}}=\frac{a_j+b_j}{2}+O(|\lambda|).$$
Then (\ref{b1}) is obtained by letting $\lambda\to0.$\qed\smallskip

When $m=0$, Theorem \ref{b.th3} implies

\begin{corollary}\label{cor} A balanced domain is in the class
$\mathcal E$ exactly when it is convex.
\end{corollary}

This follows also from Corollary \ref{bala.cor0} that uses the
Lempert theorem.

The condition $2k_{m+1}>k_n$ is essential as seen from the following
\smallskip

\noindent{\bf Example.} The $(1,2)$-balanced domain
$$D=\{z\in\Bbb C^2:|z_1|^2+3|z_2+z_1^2|<1\}$$ is not convex
(for example, $(1,0),(2i,4)\in\partial D,$ while $(1/2+i,2)\not\in\overline D)$),
but it is biholomorphic to the $(1,2)$-balanced domain
$G =\{z\in\Bbb C^2:|z_1|^2+3|z_2|<4\}$ (under the mapping $(z_1,z_2)\to(z_1,z_2+z_1^2)$).
\smallskip

Clearly the symmetrized polydisc $\Bbb G_n$ is a $(1,2,\dots,n)$-balanced domain.
However Theorem \ref{b.th3} cannot be directly applied
to show that $\Bbb G_2\not\in\mathcal E$ (since $2k_1=k_2$). Anyway its proof permits us to get

\begin{proposition}\label{b.pr4} \cite[Theorem 1]{Edi}
$\Bbb G_2\not\in\mathcal E.$
\end{proposition}

\beginproof Assume the contrary. Choose an arbitrary $\varepsilon\in(0,1)$
and put $G_\varepsilon=\{z\in\Bbb C^2:h_{\Bbb
G_2}\le\varepsilon\}.$ Then we can find a domain $D_\varepsilon$
that is biholomorphic to a convex domain and so that
$G_\varepsilon\subset D_\varepsilon\subset\Bbb G_2.$ A closer
inspection of the proof of Theorem \ref{b.th3} easily shows that
there exists a constant $c_\varepsilon\in\Bbb C$ such that
$$\left(\frac{a_1+b_1}{2},\frac{a_2+b_2}{2}+c_\varepsilon(a_1-b_1)^2\right)\in
\Bbb G_2\hbox{ for each }a,b\in D_\varepsilon.$$ If
$\theta=\arg(c_\varepsilon),$ then for
$a=\sigma(\varepsilon,i\varepsilon e^{-i\theta/2}),$
$b=\sigma(\varepsilon,-i\varepsilon e^{-i\theta/2})$ we get
$c(\varepsilon)=(\varepsilon,-4|c_\varepsilon|\varepsilon^2)\in\Bbb
G_2.$ The example at the beginning of this section implies that
$$\varepsilon\frac{1+\sqrt{1+16|c_\varepsilon|}}{2}=h_{\Bbb
G_2}(c(\varepsilon))\le 1$$ and so $\lim_{\varepsilon\to
0}c_\varepsilon=0.$ Thus $\frac{a+b}{2}\in\overline{\Bbb G_2}$ for each $a,b\in\Bbb G_2,$
i.e. $\Bbb G_2$ is a convex domain. We reached a contradiction, as $(2,1),(2i,-1)\in\partial\Bbb G_2,$ while
$(1+i,0)\not\in\overline{\Bbb G_2}.$\qed\smallskip

The above proof can be easily modified to get

\begin{proposition}\label{b.pr7} If $D$ is a balanced domain, then $\Bbb G_2\times D\not\in\mathcal E.$
\end{proposition}

Recall that $\Bbb G_2\in\mathcal L.$ When we choose a
$D\in\mathcal L$ (for example convex), Proposition \ref{b.pr7}
gives the first examples of domains in $\Bbb C^n,$ $n\ge 3$ that
are in the class $\mathcal L$ but not in the class $\mathcal E.$

We will now prove that $\Bbb G_n\not\in\mathcal E$ for $n\ge 3.$
To this aim we need the following Cohn rule that permits to learn
in finitely many steps whether the zeroes of a polynomial lie in
$\Bbb D.$

\begin{proposition}\label{b.pr5} (see e.g. \cite{Rah-Sch}) The zeroes of the
polynomial $f(\zeta)=\sum_{j=0}^n a_j\zeta^{n-j}$ ($n\ge 2,$
$a_0\neq 0$) lie in $\Bbb D$ if and only if $|a_0|>|a_n|$ and the
zeroes of the polynomial
$f^\star(\zeta)=\frac{\overline{a_0}f(\zeta)-a_n\zeta^n\overline{f(1/\overline\zeta)}}{\zeta}$
lie in $\Bbb D.$
\end{proposition}

\begin{proposition}\label{b.th6} $\Bbb G_n\not\in\mathcal E$ for $n\ge 3.$
\end{proposition}

\beginproof As $\Bbb G_n$ is a $(1,2,\dots,n)$-balanced domain,
by Theorem \ref{b.th3} it suffices to prove that if
$m=\left[n/2\right]$, then the set $G_n$ of points
$(a_{m+1},\dots,a_n)\in\Bbb C^{n-m},$ for which the zeroes of the
polynomial $f_n(\zeta)=\zeta^n+\sum_{j=m+1}^n a_j\zeta^{n-j}$ lie
in $\Bbb D,$ is convex.

First we will consider the cases $n=3$ and $n=4,$ and then we will
reduce the case $n\ge 5$ to them.
\smallskip

{\it Case $n=3.$} For $f_3(\zeta)=\zeta^3+p\zeta+q$ we have
$$f_3^\star(\zeta)=\frac{f_3(\zeta)-q\zeta^3\overline f_3(1/\overline\zeta)}{\zeta}=(1-|q|^2)\zeta^2-\overline pq\zeta+p$$ and
$$f_3^{\star\star}(\zeta)=\frac{(1-|q|^2)f_3^\star(\zeta)-p\zeta^2\overline{f_3^\star}(1/\overline\zeta)}{\zeta}$$
$$=((1-|q|^2)^2-|p|^2)\zeta-\overline pq(1-|q|^2)+p^2\overline q.$$
By Proposition \ref{b.pr5} after some calculations we get
$$G_3=\{(p,q)\in\Bbb C^2:|q|<1,\ r(p,q)<0\},$$ where
$$r(p,q)=|\overline pq(1-|q|^2)-p^2\overline q|+|p|^2-(1-|q|^2)^2.$$
It is easily seen that if $q'\in(-1,1)$ and $p'=1-q'^2,$ then $
(p_1,q_1)=(p'e^{2\pi i/3},q')$ and $ (p_2,q_2)=\left(p'e^{\pi
i/3},q'e^{\pi i/2},\right)$ are boundary points for $D$ (As
$r(p',q')=0$ and $r(p,q')<0$ for $p\in(|q'|-1,p')$). Then for
$$(p_0,q_0)=\left(\frac{p_1+p_2}{2},\frac{q_1+q_2}{2}\right)=
(p'\cos\frac{\pi}{6}e^{\pi i/2},q'\cos\frac{\pi}{4}e^{\pi i/4})$$
we have
$$|\overline{p_0}q_0(1-|q_0|^2)-p_0^2\overline{q_0}|=|p_0q_0|(1-|q_0|^2+|p_0|).$$ Consequently
$$r(p_0,q_0)=(1-|q_0|^2+|p_0|)(1+|q_0|)(|p_0|+|q_0|-1).$$ So
$r(p_0,q_0)>0$ exactly when $|p_0|+|q_0|>1.$ For $q'=1/2$
we get
$$|p_0|+|q_0|=\frac{3\sqrt3+2\sqrt2}{8}>1.$$ So
$(p_0,q_0)\not\in\overline G_3,$ showing that $G_3$ is not a convex set.
\smallskip

{\it Case $n=4.$} Similarly to the previous case we get the equality
$$G_4=\{(p,q)\in\Bbb C^2:|p|+|q|^2<1,\
s(p,q)<0\},$$ where $$s(p,q)=(1-|q|^2)|\overline
pq((1-|q|^2)^2-|p|^2)-p^3\overline
q^2|+|p|^4|q|^2-((1-|q|^2)^2-|p|^2)^2.$$ It is easily seen that if
$q'\in[0,1)$ and $p'=(1-q')\sqrt{1+q'},$ then $(p_1,q_1)=(p'e^{\pi
i/2},q')\in\partial D$ and $(p_2,q_2)=(p'e^{\pi i/4},q'e^{\pi
i/3})\in\partial D$ (as $s(p',q')=0$ and $s(p',q)<0,$ if
$p\in(-p',p')$). Then for
$$(p_0,q_0)=\left(\frac{p_1+q_1}{2},\frac{p_2+q_2}{2}\right)=(p'\cos\frac{\pi}{8}e^{3\pi i/8},q'\cos\frac{\pi}{6}e^{\pi i/6})$$ we have
$$|\overline{p_0}q_0((1-|q_0|^2)^2-|p_0|^2)-p_0^3\overline {q_0}^2|=|p_0q_0|((1-|q_0|^2)^2-|p_0|^2+|p_0|^2|q_0|).$$ So
$$s(p_0,q_0)=(1-|q_0|^2)((1-|q_0|^2)(1+|q_0|)-|p_0|^2)(1+|p_0|-|q_0|^2) (|p_0|+|q_0|-1).$$ Thus $s(p_0,q_0)>0$ only when
$|p_0|+|q_0|>1.$ For $q'=2/5$ we have
$$|p_0|+|q_0|=\frac{1}{10}\left(3\sqrt{\frac{7(2+\sqrt2)}{5}}+2\sqrt3\right)>1.$$
Consequently $(p_0,q_0)\not\in\overline{G_4}$ and so $G_4$ is not a convex set.
\smallskip

{\it Case $n\ge 5.$} Let $j\in\{0,1,2\}.$ Note that the non-convex
set $G_3$ coincides with the set of points $(p,q)\in\Bbb C^2$ such
that the zeroes of the polynomial $z^jf_3(z^k),$ $k\ge1,$ lie in
$\Bbb D.$ Consequently for $n=3k+2$ and $k\ge 3$, $n=3k+1$ and
$k\ge 2,$ or $n=3k$ and $k\ge 1,$ we can view $G_3$ as an
intersection of $G_n$ with a complex plane. So $G_n$ is not a
convex set.

In the remaining cases $n=5$ and $n=8$ it suffices to observe that
the non-convex set $G_4$ coincides with the set of points
$(p,q)\in\Bbb C^2$ such that the zeroes of the polynomials
$\zeta^4 f_4(\zeta)$, respectively $f_4(\zeta^2)$, lie in $\Bbb
D,$ and then conclude the proof as above.\qed\smallskip

\setcounter{equation}{0}
\section{Notions of complex convexity}\label{con}

The main definitions and facts from this section can be found in \cite{APS,Hor2} (see also \cite{Jac2}).

Recall that a domain is called $\Bbb C$-convex, if all its
intersections with complex lines are connected and simply
connected.

We will define two other notions for complex convexity. A domain
in $\Bbb C^n$ is called linearly convex if each point in its
complement belongs to a complex hyperplane, disjoint from the
domain. If the last is true for each boundary point, then the
domain is called weakly linearly convex. The following
implications take place:
\smallskip

convexity $\Rightarrow$ $\Bbb C$-convexity $\Rightarrow$ linear convexity

$\Rightarrow$ weak linear convexity $\Rightarrow$ pseudoconvexity.
\smallskip

On the other hand, in \cite[Example 2.2.4]{APS} it is shown that a
complete Reinhardt domain, which is weakly linearly convex, is
convex. (A domain $D$ in $\Bbb C^n$ is called a complete Reinhardt
domain if for $z\in D$ the closed polydisc of center $0$ and
radius $z$ lies in $D.$) We will see that this result remains true
for balanced domains (but not for generalized balanced domains, as
shown by Theorem \ref{c.th5}).

\begin{proposition}\label{c.pr1} A weakly linearly convex balanced domain $D\subset\Bbb C^n$ is convex.
\end{proposition}

For the proof of Proposition \ref{c.pr1} we will use a
characterization of the (weakly) linearly convex domains. For a
set $D$ in $\Bbb C^n,$ containing the origin, put
$$D^\ast=\{z\in\Bbb C^n: \langle z,w\rangle\neq 1\hbox{ for each }w\in D\}$$
($\langle \cdot,\rangle$ is the Hermitian scalar product). Clearly
if $D$ is open (compact), then $D^\ast$ is compact (open).
Furthermore, $D\subset D^{\ast\ast}.$

\begin{proposition}\label{c.pr2} (see e.g. \cite{APS,Hor2})
A domain $D$ in $\Bbb C^n,$ containing the origin, is weakly
linearly convex (linearly convex) if and only if $D$ is a
component of $D^{\ast\ast}$ ($D=D^{\ast\ast}$).
\end{proposition}

\noindent{\it Proof of Proposition \ref{c.pr1}} As $D$ is a
balanced domain, then it is easily seen that $D^\ast$ is a compact
balanced set.
 Consequently $D^{\ast\ast}$ is an open balanced set and in particular a domain.

We will prove that this domain is convex. Assume the contrary.
Then there exist points $z_1,z_2\in D^{\ast\ast},$ $w\in D^\ast$
and a number $t\in(0,1)$ so that $\langle
tz_1+(1-t)z_2,w\rangle=1$. Consequently we can suppose that
$|\langle z_1,w\rangle|\ge 1.$ As $D^\ast$ is a balanced set,
$\tilde w=w/\langle w,z_1\rangle\in D^\ast$ and $\langle
z_1,\tilde w\rangle=1$ -- a contradiction.

So $D^{\ast\ast}$ is a convex domain. Since $D$ is a weakly
linearly convex domain, $D$ is a component of $D^{\ast\ast}$ and
consequently $D=D^{\ast\ast}.$\qed\smallskip

Let us note that the three notions of complex convexity are
different, but for bounded domains with $\mathcal C^1$-smooth
boundaries they coincide (in the more general case of bounded
domains this is not true). We will also mention that each $\Bbb
C$-convex domain in $\Bbb C^n$ is homeomorphic to $\Bbb C^n,$ and
each domain in $\Bbb C$ is linearly convex. Also, a Cartesian
product of (weakly) linearly convex domains is (weakly) linearly
convex. On the other hand, we have the following
\smallskip

\noindent{\bf Remark.} A Cartesian product of domains that do not
coincide with the corresponding spaces is $\Bbb C$-convex only if
both the components are convex. In particular, a Cartesian product
of $n$ simply connected non-convex domains from $\Bbb C$ is a
linearly convex domain that is biholomorphic to $\Bbb D^n,$ yet
not $\Bbb C$-convex.
\smallskip

Recall that a domain $D$ with a $\mathcal C^2$-smooth boundary is
convex (pseudoconvex), if the restriction of the Hessian (Levi
form) of its defining function on the real (complex) tangent space
at each boundary point of $D$ is a positively semidefinite
quadratic form. The following fact confirms the intermediate
character of complex convexity: a domain $D$ with a $\mathcal
C^2$-smooth boundary is $\Bbb C$-convex exactly when the
restriction of the Hessian of its defining function on the complex
tangent space at each boundary point of $D$ is a positively
semidefinite quadratic form. The last turns out to be equivalent
to the function $-2\log d_D(z)$ near $\partial D$ being of the
class of the so-called $\C$-convex functions (see e.g. \cite{APS};
the number 2 is important); here $d_D(z)=\dist(z,\partial D),$
$z\in D.$ The pseudoconvex analogue of this Proposition without a
smoothness condition is well-known. Of course we have also a
convex analogue, which is given in \cite[Proposition 7.1]{Her-McN}
for bounded domains with $C^2$-smooth boundaries. To get convinced
in this for an arbitrary domain $D,$ one can note that convexity
of $-d_D$ (or, which is the same, of $D$) trivially implies
convexity of $-\log d_D.$ The converse is also true; it suffices
to assume the contrary and then find a segment that, except for
its midpoint, lies within $D$ (see e.g. \cite[Theorem
2.1.27]{Hor2} for a more general fact).

For bounded domains with $\mathcal C^1$-smooth boundaries the
three notions of complex convexity also coincide with the
so-called weak local linear convexity (see e.g. \cite[Proposition
4.6.4]{Hor2}). A domain $D\subset\Bbb C^n$ is called a weakly
locally linearly convex if for each point $a\in\partial D$ there
exists a complex hyperplane $H_a$ through $a$ and a neighborhood
$U_a$ of $a$ so that $H_a\cap D\cap U_a=\varnothing.$ Note that
there are bounded domains that are not locally linearly convex
(see e.g. \cite{Jac2}). In \cite[p. 58]{Jac2} it is asked whether
a weakly locally linearly convex domain needs to be pseudoconvex.

The next proposition gives more than an affirmative answer to this question.

\begin{proposition}\label{c.pr0}
Let $D\subset\Bbb C^n$ be a bounded domain with the following
property: for each point $a\in\partial D$ there exists a
neighborhood $U_a$ of $a$ and a function $f_a\in\mathcal O(D\cap
U_a)$ so that $\lim_{z\to a}|f_a(z)|=\infty.$ Then $D$ is a taut
domain (in particular, pseudoconvex).
\end{proposition}

\beginproof It suffices to prove that if $\mathcal O(\Bbb D,D)\ni\psi_j\to\psi$
 and $\psi(\zeta)\in\partial D$ for some $\zeta\in\Bbb D,$ then $\psi(\Bbb D)\subset\partial D.$
Assume the contrary. Then we can easily find points
$\eta_k\to\eta\in\Bbb D$ so that $\psi(\eta_k)\in D,$ but
$a=\psi(\eta)\in\partial D.$ We can assume that $\eta=0$ and
$g_a=1/f_a$ is a bounded function on $D\cap U_a.$ Let $r\in(0,1)$
is such that $\psi(r\Bbb D)\Subset U_a.$ Then $\psi_j(r\Bbb
D)\subset U_a$ for each $j\ge j_0.$ Consequently
$|g_a\circ\psi_j|<1$ and (by passing to subsequences) we can
suppose that $g_a\circ\psi_j\to h_a\in\mathcal O(r\Bbb D,\Bbb C).$
As $h_a(\eta)=0,$ by Hurwitz's theorem $h_a=0$. This contradicts
the fact that $h_a(\eta_k)=g_a\circ\psi(\eta_k)\neq 0$ for
$|\eta_k|<r.$\qed

\begin{corollary}\label{c.pr00} A weakly locally linearly convex domain is pseudoconvex.
\end{corollary}

For the proof of this corollary it is sufficient to exhaust the
domain with bounded domains and for each boundary point to
consider the reciprocal of the defining function of the
corresponding separating hyperplane.
\smallskip

Later note that a linearly convex domain $D\subset\C^n,$
containing a complex line, is linearly equivalent to $\Bbb C\times
D',$ where $D'\subset\C^{n-1}$ \cite[Proposition 4.6.11]{Hor2}.
Indeed, we can assume that $D$ contains the $z_1$-line. As the
complement $^cD$ of $D$ is a union of complex hyperplanes not
intersecting this line, $^cD=\Bbb C\times G$ and consequently
$D=\Bbb C\times{^cG}.$

The next theorem provides some properties of the $\Bbb C$-convex
domains not containing complex lines. It generalizes a result of
T. J. Barth from \cite{Bar} for convex domains.

\begin{theorem}\label{c.pr4}
Let $D$ be a $\Bbb C$-convex domain in $\Bbb C^n$ not containing a
complex line. Then $D$ is biholomorphic to a bounded domain and is
$c$-finite compact, hence also $c$-complete ($c$ is the
Carath\'eodory distance). In particular, $D$ is hyperconvex, so it
is a taut domain.
\end{theorem}

Based on this theorem, the paper \cite{Nik-Sar} by the author and
A. Saracco includes various equivalent conditions for a
$\C$-convex domain not to contain a complex line (the convex case
is treated in \cite{Bra-Sar}).

\beginproof
For each point $z\in{^cD}$ denote by $L_z$ some complex hyperplane
through $z$ disjoint from $D.$ Let $l_z$ be a line through the
origin that is orthogonal to $L_z.$ Denote by $\pi_z$ the
orthogonal projection of $\Bbb C^n$ onto $l_z$ and put
$a_z=\pi_z(a)$ (clearly $\pi_z$ and $\pi_t$ from Section \ref{bal}
refer to different objects). The set $D_z=\pi_z(D)$ is
biholomorphic to $\Bbb D,$ since it is connected, simply connected
(see e.g. \cite[Theorem 2.3.6]{APS}) and
$\pi_z(z)\not\in\pi_z(D).$ As $D$ is a linearly convex domain not
containing complex lines, it is easily seen that there exist $n$
$\Bbb C$-independent $l_z'$ (otherwise $^cD,$ and so $D,$ would
contain a complex line that is orthogonal to each $l_z.$). We can
assume that these $l_z$ form the set $C$ of coordinate lines. Then
$D\subset G=\prod_{l_z\in C}\pi_z(D)$ and $G$ is biholomorphic to
the polydisc $\Bbb D^n$ (As the components of $G$ are simply
connected planar domains $\neq\Bbb C$, so they are biholomorphic
to $\Bbb D$ according to the Riemann theorem). Consequently $D$ is
biholomorphic to a bounded domain, so it is $c$-hyperbolic.

Later we can assume that $0\in D.$ To see that $D$
is a $c$-finite compact, it suffices to show that $\lim_{a\to z}c_D(0;a)=\infty$
 for each $z\in\partial D$ (recall that
$\infty\in\partial D,$ if $D$ is unbounded). The last statement follows from the
fact that $G\supset D$ is a $c$-finite compact domain
(being biholomorphic to $\Bbb D^n$). On the other hand, if $a\to
z\in\partial D,$ then $a_z\to\pi_z(z)\in\partial D_z$ and
consequently $c_D(0;a)\ge c_{D_z}(0;a_z)\to\infty.$\qed\smallskip

The main aim of this section is to show that $\Bbb G_2$ is a $\Bbb
C$-convex domain, which together with the fact that $\Bbb
G_2\not\in\mathcal E$ (see Proposition \ref{b.pr4}) gives a
negative answer to \cite[Problem 4]{Zna} (see the Introduction).

\begin{theorem}\label{c.th5} (i) $\Bbb G_2$ is a $\Bbb C$-convex domain.

(ii) $\Bbb G_n,$ $n\ge 3,$ is a linearly convex domain, yet not $\Bbb C$-convex.
\end{theorem}

\noindent{\bf Remarks.} a) Proposition \ref{b.pr7} implies that
$\Bbb G_2\times\Bbb C^n$ is a $\Bbb C$-convex domain that is in
the class $\mathcal L,$ but not in the class $\mathcal E$
according to Proposition \ref{b.pr7}. However we do not have a
similar example of a bounded domain in $\Bbb C^n,$ $n\ge 3.$ (The
most natural candidate is the Cartesian product, but according to
a remark above this is impossible, as the components need to be
convex).
\smallskip

b) It is easily seen that the bounded generalized balanced domains
with continuous Minkowski functional are homeomorphic to $\Bbb
C^n.$ So, Theorem \ref{c.th5} (ii) provides the first example of a
linearly convex domain, namely $\G_n$, that is homeomorphic to
$\Bbb C^n,$ $n\ge 3,$ but is not $\Bbb C$-convex, not in the class
$\mathcal L$ and not a scalar product.

c) Theorem \ref{c.th5} (ii) shows that Proposition \ref{c.pr1}
does not remain true for generalized balanced domains.\smallskip

d) Although $\Bbb G_n$ for $n\ge 3$ is not a $\Bbb C$-convex domain,
the conclusion of Theorem \ref{c.pr4} is true, i.e.
$\Bbb G_n$ is a $c$-finite compact domain. This fact follows directly
from \cite[Corollary 3.2]{Cos3} (see (\ref{i2})).
\smallskip

Let us introduce the following notation. Let $D$ be a domain in $\Bbb C^n$
containing the origin, and $0\neq a\in\partial D.$ Let
$\Gamma_D(a)$ be the set of points $z\in\C^n$ such that the hyperplane
$\{w\in\C^n:\langle z,w\rangle=1\}$ passing through $a$ does not meet $D.$

For the proof of Theorem \ref{c.th5} we will use the following
characterization of the bounded $\Bbb C$-convex domains.

\begin{proposition}\label{c.pr6} \cite[Theorem 2.5.2]{APS}
A bounded domain $D$ in $\Bbb C^n$ ($n>1$) containing the origin is
$\Bbb C$-convex if and only if for each $a\in\partial D$ the set
$\Gamma_D(a)\subset\Bbb C\Bbb P^n$ is nonempty and connected.
\end{proposition}

\noindent{\bf Remark.} Proposition \ref{c.pr6} directly implies the
mentioned fact that a $\mathcal C^1$-smooth bounded
domain $D$ in $\Bbb C^n$, $n>1,$ is $\Bbb C$-convex if and only if it is linearly convex.
\smallskip

\noindent{\it Proof of Theorem \ref{c.th5}.} (i) According to
Proposition \ref{c.pr6} we need to check that the set
$\Gamma(a)$ is nonempty and connected for each $a\in\partial D.$

First let us take a smooth boundary point $a\in\partial\Bbb G_2$;
without loss of generality it is of the form $\sigma(\mu)$, where
$|\mu_1|=1$, $|\mu_2|<1.$ Then the tangent plane to $\partial\Bbb
G_2$ at $a$ has the form $\{\sigma(\mu_1,\lambda) :\lambda\in\Bbb
C\}$ and clearly it does not meet $\Bbb G_2$. So in this case
$\Gamma(a)$ is a singleton.

Let now $a\in\partial\Bbb G_2$ be a non-smooth boundary point,
i.e. $a=\sigma(\mu)$, where $|\mu_1|=|\mu_2|=1.$ After a rotation
we can assume that $\mu_1\mu_2=1$, i.e. $\mu_2=\overline\mu_1$.
Then $\mu_1+\mu_2=2\Re\mu_1=:2x$, where $x\in[-1,1]$.

The complex lines through $a$ that meet $\Bbb G_2$ have the form
$a+\Bbb C(a-\sigma(\lambda))$, where $\lambda\in\Bbb D^2.$
Consequently the complement of $\Gamma(a)$ can be seen as the set
$$A=\{\frac{\lambda_1+\lambda_2-2x}{\lambda_1\lambda_2-1}: \lambda_1,\lambda_2\in\Bbb D\}.$$

Connectedness of $\Gamma(a)$ means simply connectedness of $A.$
Note that if $|\beta|>1$, then $\frac{z-\alpha}{z-\beta}$ maps the
unit disc $\Bbb D$ in the disc $\Bbb
D(\frac{1-\alpha\overline\beta}{1-|\beta|^2},
\frac{|\alpha-\beta|}{|\beta|^2-1})$. So the set
$\{\frac{\lambda+\lambda_1-2x}{\lambda\lambda_1-1}:\lambda\in\Bbb
D\}$ coincides with
$$A_{\lambda_1}=\Bbb D(\frac{2x-2\Re\lambda_1}{1-|\lambda_1|^2},
\frac{|2x\lambda_1-\lambda_1^2-1|}{1-|\lambda_1|^2}).$$ As
$A=\bigcup\sb{\lambda_1\in\Bbb D}A_{\lambda_1}\subset\Bbb C,$ $A$
is a simply connected set.
\smallskip

(ii) To prove the linear convexity of $\Bbb G_n,$ consider a point
$z=\sigma(\lambda)\in\Bbb C^n\setminus\Bbb G_n$. We can assume
that $|\lambda_1|\geq 1$. Then the set
$$B=\{\sigma(\lambda_1,\mu_1,\ldots,\mu_{n-1}):\mu_1,\ldots,\mu_{n-1}\in\Bbb
C\}$$ is disjoint from $\Bbb G_n$. On the other hand, it is easily
seen that $$B=\{(\lambda_1+z_1,\lambda_1z_1+z_2,
\ldots,\lambda_1z_{n-2}+z_{n-1},\lambda_1z_{n-1}):z_1,\ldots,z_{n-1}\in\Bbb
C\};$$ so $B$ is a complex hyperplane. Consequently $\Bbb G_n$ is
a linearly convex domain.

To prove that $\Bbb G_n$ is not a $\Bbb C$-convex domain for
$n\geq 3,$ we consider the points
$$a_t=\sigma(t,t,t,0,\ldots,0)=(3t,3t^2,t^3,0,\ldots,0), $$
$$b_t=\sigma(-t,-t,-t,0,\ldots,0)=(-3t,3t^2,-t^3,0,\ldots,0),\;t\in(0,1).$$
Clearly $a_t,b_t\in\Bbb G_n$. Denote by $L_t$ the complex line through $a_t$ and $b_t,$ i.e.
$$L_t=\{c_{t,\lambda}:=(3t(1-2\lambda),3t^2,t^3(1-2\lambda),0,\ldots,0): \lambda\in\Bbb C\}.$$
Suppose that $\Bbb G_n\cap L_t$ is a connected set. As $a_t=c{t,0}$ and $b_t=c{t,1},$
$c_{t,\lambda}\in\Bbb G_n$ for some $\lambda=1/2+i\tau,$
$\tau\in\Bbb R.$ Then $$c_{t,\lambda}=(-6i\tau t,3t^2,-2i\tau
t^3,0,\ldots,0).$$ Let $c_{\tau}=\sigma(\mu),$ $\mu\in\Bbb D^n.$
We can assume that $\mu_j=0$, $j=4,\ldots,n$ and $-36\tau^2t^2=(\mu_1+\mu_2+\mu_3)^2=\mu_1^2+\mu_2^2+\mu_3^2+6t^2.$
Then
$$t^2=\frac{|\mu_1^2+\mu_2^2+\mu_3^2|}{36\tau^2+6}<\frac{3}{36\tau^2+6}
\leq\frac{1}{2}.$$ Consequently $\Bbb G_n\cap L_t$ is not a connected
set for $t\in[1/\sqrt{2},1)$ and so $\Bbb G_n$ is not a $\Bbb C$-convex domain.\qed

The fact that $\Bbb G_2$ is a $\Bbb C$-convex domain is also a
consequence of a recent result by P. Pflug and W. Zwonek
\cite{Pfl-Zwo2} which also confirms the Aizenberg hypothesis. Let
first recall that a $C^2$-smooth domain $D$ in $\C^n$ is said to
be strongly linearly convex if the restriction of the Hessian of
its defining function on the complex tangent space at each
boundary point of $D$ is a positively definite quadratic form.

Considering functions of the form
$$r_{\eps}(s,p)=|s-\overline{s}p|-(1-|p|^2)^2+\eps,$$
one can show the following

\begin{theorem}\label{c.th6} (see \cite{Pfl-Zwo}) The domain
$$\Bbb G_2^\eps=\{(s,p):\sqrt{|s-\overline{s}p|+\eps}+|p|^2<1,\quad \eps\in(0,1)$$
is strongly linearly convex. Consequently, $\Bbb G_2=\Bbb G_2^0$
can be exhausted by strongly linearly convex domains and hence
$\Bbb G_2$ is a $\C$-convex domain.
\end{theorem}

Since $\Bbb G_2\not\in\mathcal E,$ we get immediately

\begin{corollary} $\Bbb G_2^\eps\not\in\mathcal E$ for $\eps>0$ small
enough.
\end{corollary}

\setcounter{equation}{0}
\section{$\Bbb G_n\not\in\mathcal L$ for $n\ge 3$}\label{ine}

As mentioned in the Introduction, $\Bbb G_2\in\mathcal L,$ i.e. the
Carath\'eodory and Lempert functions of the symmetrized polydisc coincide.
The main aim of this section is to show that this property no longer holds in
higher dimensions, i.e. $\Bbb G_n\not\in\mathcal
L$ for $n\ge 3,$ which solves \cite[Problem 1.4]{Jar-Pfl3}.

Let $n\ge 2,$ $z\in\C^n$ and $\lambda\in\overline{\Bbb D}.$
Put $$f_\lambda(z)=\frac{\sum_{j=1}^{n}jz_j
\lambda^{j-1}}{n+\sum_{j=1}^{n-1}(n-j)z_j\lambda^j}.$$ By
\cite[Theorem 3.1]{Cos3}, $z\in\Bbb G_n$ if and only if $\sup_{\lambda\in\overline{\Bbb
D}}|f_\lambda(z)|<1.$ So for the function of Carath\'eodory we have
\begin{equation}\label{i1}
c^\ast_{\Bbb G_n}(z,w)\ge m_{\Bbb G_n}(z,w):=\max_{\lambda\in\Bbb
T}|m_{\Bbb D}(f_\lambda(z),f_\lambda(w))|.
\end{equation}
Note that $m_{\Bbb G_n}$ is a distance on $\Bbb G_n$.
Furthermore, by \cite[Corollary 3.2]{Cos3} $\lim_{w\to\partial\Bbb G_n}m_{\Bbb G_n}(z,w)=1$
and consequently
\begin{equation}\label{i2}
\lim_{w\to\partial\Bbb G_n}c_{\Bbb G_n}(z,w)=\infty,
\end{equation}
i.e. $\Bbb G_n$ is a $c$-finite compact domain.

The next basic proposition is used in the proof of Propositions \ref{np.pr4} and \ref{np.pr5}.
It bears information for the zeroes of $l_{\Om_n}$ and $c^\ast_{\Om_n}.$

\begin{proposition}\label{i.pr4} Let $A,B\in\Omega_n$
and $$s(A,B)=\min_{\lambda\in \spe(A)}\max_{\mu\in \spe(B)}m_{\Bbb
D}(\lambda,\mu).$$ Then:

(i) $l_{\Bbb G_n}(\sigma(A),\sigma(B)))\le l_{\Omega_n}(A,B)\le
s(A,B);$

(ii) $l_{\Omega_n}(A,B)=0\Leftrightarrow\
c^\ast_{\Omega_n}(A,B)=0\Leftrightarrow
\spe(A)=\spe(B)\Leftrightarrow\exists\varphi\in\mathcal O(\Bbb
C,\Omega_n):\varphi(0)=A,\varphi(1)=B.$

(iii) if the eigenvalues of $A$ are equal, then the eigenvalues of $B$ are
equal $\Leftrightarrow c^\ast_{\Bbb
G_n}(\sigma(A),\sigma(B))=s(A,B)\Leftrightarrow l_{\Bbb
G_n}(\sigma(A),\sigma(B))=s(A,B).$
\end{proposition}
\vskip-2mm

\noindent{\it Proof.} (i) The inequality on the left is noted in the Introduction
(it follows from the holomorphic contractibility of the Lempert function).

To prove the inequality on the right, let $J_A$ and $J_B$ be
Jordan normal forms of $A$ and $B$, respectively. As a nonsingular
square matrix $X$ can be expressed in the form $X=e^Y,$ we get
$A=e^{Y_A}J_Ae^{-Y_A}$ and $B=e^{Y_B}J_Be^{-Y_b},$ where
$Y_A,Y_B\in\mathcal M(\Bbb C^n).$ We have
$J_a=(a_{jk})_{j,k=1}^n,$ where $\spe(A)=\{a_{11},\dots,a_{nn}\},$
$a_{j,j+1}=0,1$ and $a_{jk}=0,$ if $j<k$ or $j>k+1.$ A similar
representation is valid for $J_b=(b_{jk})_{j,k=1}^n.$ Let
$\lambda\in\spe(A)$ and $\tilde\lambda=\max_{\mu\in
\spe(B)}p_{\Bbb D}(\lambda,\mu).$ Then one can easily find
$\varphi_{jj}\in\mathcal O(\Bbb D,\Bbb D)$ so that
$\varphi_{jj}(0)=a_{jj}$ and $\varphi_{jj}(\tilde\lambda)=b_{jj}.$
Clearly we can choose $\varphi_{j,j+1}\in\mathcal O(\Bbb D,\Bbb
D)$ and $\psi\in\mathcal O(\Bbb D,\mathcal M_n)$ so that
$\varphi_{j,j+1}(0)=a_{j,j+1},$
$\varphi_{jj}(\tilde\lambda)=b_{jj}$ and $\psi(0)=Y_A,$
$\psi(\tilde\lambda)=B.$ Put $\varphi_{jk}=0,$ if $j<k$ or
$j>k+1,$ and $\varphi=(\varphi_{jk})_{j,k=1}^n.$ Then
$e^\psi\varphi e^{-\psi}\in\mathcal O(\Bbb D,\Omega_n)$ and
$\varphi(0)=A,$ $\varphi(\tilde\lambda)=B,$ which completes the
proof. Consequently $l_{\Omega_n}(A,B)\le|\tilde\lambda|$ and as
$\lambda\in\spe(A)$ was arbitrary, we get the inequality on the
right.

(ii) Clearly
$$l_{\Omega_n}(A,B)=0\Rightarrow
c^\ast_{\Omega_n}(A,B)=0=0\Rightarrow \spe(A)=\spe(B),$$ since
$c^\ast_{\Omega_n}(A,B)=c^\ast_{\Bbb G_n}(\sigma(A),\sigma(B))$
according to Proposition \ref{np.pr5}. If $\spe(A)=\spe(B),$ then
as in the proof of the inequality on the right (i), we can find
$\varphi\in\mathcal O(\Bbb C,\Omega_n)$ so that $\varphi(0)=A$ and
$\varphi(1)=B,$ leading to the implication
$\spe(A)=\spe(B)\Rightarrow l_{\Omega_n}(A,B)=0$

(iii) By $c^\ast_{\G_n}\le l_{\G_n}\le s$ (see (i) for the last one)
we get the implication $c^\ast_{\Bbb G_n}(\sigma(A),\sigma(B))=s(A,B)
\Rightarrow l_{\Bbb G_n}(\sigma(A),\sigma(B))=s(A,B).$

Later, using $\Phi_\lambda$ (see (\ref{fi})) we can assume that the
eigenvalues of $A$ are equal to $0.$

To prove that if the eigenvalues of $B$ are equal e.g. to $\mu,$
then $c^\ast_{\Bbb G_n}(0,\sigma(B))=s(0,B)(=|\mu|),$ it suffices
to construct a function $f\in\mathcal O(\Bbb G_n,\Bbb D)$ so that
$f(0)=0$ and $f(\sigma(B))=\mu.$ An example of such a function is
$f(z_1,\dots,z_n)=\frac{z_1+\dots+z_n}{n}.$

It remains to prove that if $l_{\Bbb G_n}(0,\sigma(B))=s(0,B),$ then
the eigenvalues of $B$ are equal. Let
$\spe(B)=(\nu_1,\dots,\nu_n).$ Let $\root n\of
1=\{1,\varepsilon,\dots,\varepsilon^{n-1}\}.$ For each
Blaschke product $\mathcal B$ of order $\le n$ such that
$\mathcal B(0)=0,$ consider the mapping $$\lambda\to
f_{\mathcal B}(\lambda)=\sigma(\mathcal B(\root
n\of{\lambda}),\mathcal B(\varepsilon\root n\of
{\lambda}),\dots,\mathcal B(\varepsilon^{n-1}\root
n\of{\lambda}))$$ (where $\root n\of{\lambda}$ is arbitrarily chosen).
It is easily seen that $f_{\mathcal B}\in\mathcal
O(\Bbb D,\Bbb G_n).$ We need the following

\begin{lemma}\label{i.l5} Let $\delta_1,\dots,\delta_n\in\Bbb T$
be pairwise different points. Then for aritrary points
$\nu_1,\dots,\nu_n\in\Bbb D$ there exists $\beta\in\Bbb D$ and a Blaschke
product $\mathcal B$ of order $\le n$ such that
$$\mathcal B(0)=0,\mathcal B(\delta_1\beta)=\nu_1,\dots,\mathcal
B(\delta_n\beta)=\nu_n.$$
\end{lemma}

\beginproof Let $S$ be the set of all $\beta\in\Bbb D$ such that the
classical Nevanlinna--Pick problem
with data
$(0,0),(\delta_1\beta,\nu_1),\dots,(\delta_n\beta,\nu_n)$ has a
solution, i.e. there exists $f\in\mathcal O(\Bbb D,\Bbb D)$ so
that $f(\delta_j\beta)=\nu_j,$ $1\le j\le n.$ Recall that this
condition is equivalent to the positive semidefiniteness of
$A(\beta)=[a_{j,k}(\beta)]_{j,k=1}^n,$ where $
a_{j,k}(\beta)=\frac{1-\nu_j\overline\nu_k}{1-\delta_j\overline\delta
_k|\beta|^2}$ (see e.g. \cite{Gar}). Note that the function
$a_{j,k},$ $j\neq k,$ is bounded on $\Bbb D.$ On the other hand,
$\lim_{\beta\to\Bbb T}a_{j,j}(\beta)=+\infty.$ Consequently the
matrix $A(\beta)$ is positively semidefinite when $\beta$ is close
to $\Bbb T.$ We can assume that $0\not\in S$ (otherwise put
$\mathcal B=\id$) and then $S$ is a proper nonempty closed subset
of $\Bbb D$ (consisting of circles). So it has a boundary point
$\beta_0\in\Bbb D.$ Consequently the number $m=\hbox{rank
}A(\beta_0)$ is not maximal, i.e. $m<n+1.$ Thus the corresponding
problem of Nevanlinna--Pick has a unique solution, which is a
Blaschke product of order $m$ (see e.g. \cite{Gar}).\qed\smallskip

The lemma for $\delta_j=\varepsilon^j,$ $1\le j\le n,$ implies that
$$l_{\Bbb G_n}(0,\sigma(B))\le|\beta|^n.$$ It remains to prove that
if $|\beta|^n\ge|\nu_j|$ for each $1\le j\le n,$ then
$\nu_1=\dots=\nu_n.$ After a rotation we can assume that
$$\mathcal B(z)=z\frac{a_0z^k+a_1z^{k-1}+\dots+a_k}{\overline a_k
z^k+\overline a_{k-1}z^{k-1}+\dots+\overline a_0},$$ where
$a_0=1$ and $k\le n-1.$ As $|\nu_j|\ge|\mathcal
B(\varepsilon^j\beta)|,$ we get $|\beta|^n\ge |\mathcal
B(\varepsilon^j\beta)|,$ i.e.
$$|\beta|^{n-1}|{\overline
a_k(\varepsilon^j\beta)^k+\overline
a_{k-1}(\varepsilon^j\beta)^{k-1}+\dots+\overline a_0}|$$
$$\ge|a_0(\varepsilon^j\beta)^k+a_1(\varepsilon^j\beta)^{k-1}+\dots+a_k|.$$
By squaring we get
$$|\beta|^{2n-2}(\sum_{s=0}^k|a_s|^2|\beta|^{2s}+2\Re\sum_{0\le
p<s\le k}a_p\overline
a_s\beta^s\overline\beta^p\varepsilon^{j(s-p)})\ge$$
$$\sum_{s=0}^k|a_s|^2|\beta|^{2(k-s)}+2\Re\sum_{0\le p<s\le
k}a_p\overline
a_s\beta^{k-p}\overline\beta^{k-s}\varepsilon^{j(s-p)}.$$
Adding these inequalities for $j=1,\dots,n$ yields
$$|\beta|^{2n-2}\sum_{s=0}^k|a_s|^2|\beta|^{2s}\ge
\sum_{s=0}^k|a_s|^2|\beta|^{2(k-s)},$$ i.e.,
$$\sum_{s=0}^k|a_s|^2(|\beta|^{2(n+s-1)}-|\beta|^{2(k-s)}|)\ge
0.$$ As $k\le n-1,$ then $k-s<n+s-1$ if $s>0$ and so
$a_s=0.$ On the other hand, $a_0=1\neq 0$ and consequently $k=n-1.$
We thus got $\mathcal B(z)=z^n$ and hence $\nu_1=\dots=\nu_n.$ \qed
\smallskip

Let $e_1,\dots,e_n$ be the standard basisi in $\Bbb C^n$ and $X=\sum_{j=1}^n X_je_j.$ Put
$$\tilde f_{\lambda}(X)=\frac{\sum_{j=1}^{n}jX_j \lambda^{j-1}}{n}\hbox{\ \
and\ \ }\rho_n(X)=\max_{\lambda\in\Bbb T}|\tilde f_\lambda(X)|.$$ By
(\ref{i1}) we get the following estimate for the Carath\'eodory metric of $\G_n:$
$$\gamma_{\Bbb G_n}(0;X)\ge\lim_{\Bbb C_\ast\ni t\to 0}\frac{p_{\Bbb G_n}(0,tX)}{|t|}=\rho_n(X).$$

Let $L_{k,l}=\span(e_k,e_l).$ Clearly if $X \in L_{k,l}$,
$k\neq l$, then $$\rho_n(X)=\frac{k|X_k|+l|X_l|}{n}.$$

As noted, one of the basic results that motivate the discussion of the symmetrized
bidisc is the fact that
$\Bbb G_2\in\mathcal L.$ More precisely (see \cite{Agl-You3}), $$l_{\Bbb
G_2}=k^\ast_{\Bbb G_2}=c^\ast_{\Bbb G_2}=m_{\Bbb G_2};$$ in particular $\kappa_{\G_2}=\gamma_{\G_2}.$
\smallskip

The next proposition shows that $\Bbb G_n$ does not have similar properties for $n\ge 3.$

\begin{theorem}\label{i.th1} (i) If $k$ divides $n,$ then $\kappa_{\Bbb
G_n}(0;e_k)=\rho_n(e_k).$ Consequently, if also $l$ divides $n,$ then
$\hat\kappa_{\Bbb G_n}(0;X)=\gamma_{\Bbb G_n}(0;X)=\rho_n(X)$ for
$X\in L_{k,l}.$

(ii) If $n\ge 3$ and $X\in L_{1,n}\setminus(L_{1,1}\cup L_{n,n}),$
then $\kappa_{\Bbb G_n}(0;X)>\rho_n(X).$

(iii) If $k$ does not divide $n,$ then $\gamma_{\Bbb
G_n}(0;e_k)>\rho_n(e_k).$
\end{theorem}

As $\G_n$ is a taut domain, Theorem \ref{der.th1} and (\ref{car})
imply

\begin{corollary}\label{i.cor2} If $n\ge 3,$ then
$$l_{\Bbb G_n}(0,\cdot)\gneq k^\ast_{\Bbb G_n}(0,\cdot)\ge
c^\ast_{\Bbb G_n}(0,\cdot)\gneq m_{\Bbb G_n}(0,\cdot).$$
\end{corollary}

\noindent{\bf Remark.} Clearly $\Bbb G_2$ is a domain that is biholomorphic to
${\Bbb G_{2n}}\cap L_{n,2n}.$ Then, unlike
Theorem \ref{i.th1}, for $z,w\in L_{n,2n}$ we have $$m_{\Bbb
G_{2n}}(z,w)\le m_{\Bbb G_{2n}}(z,w)\le l_{\Bbb G_2}(z,w)=m_{\Bbb
G_2}(z,w)\le m_{\Bbb G_{2n}}(z,w)$$ and so $l_{\Bbb
G_{2n}}(z,w)=m_{\Bbb G_{2n}}(z,w).$
\smallskip

\begin{corollary}\label{i.cor3} The convex hull of the spectral unit ball
$\Omega_n$ is $$\hat\Omega_n=\{A\in\mathcal M_n(\Bbb
C):h_{\hat\Omega_n}(A)=|\tr A|/n<1\}.$$
\end{corollary}

\beginproof By (\ref{np5}) we get $$h_{\hat\Omega_n}(A)=\mathcal
Dk_{\Omega_n}(0;A)=\lim_{t\to 0}\frac{k_{\Bbb
\Omega_n}(0;tA)}{|t|}=\lim_{t\to 0}\frac{k_{\Bbb
G_n}(0,\sigma(tA))}{|t|}.$$ As
$$\sigma(tA)=(t\tr A+o(t),o(t),\dots,o(t))$$ and $\Bbb G_n$ is a
taut domain, Theorem \ref{i.th1} (i) implies that the last limit
equals $$\hat\kappa_{\Bbb G_n}(0;(\tr A)e_1)=|\tr
A|/n.\qed$$\smallskip

\noindent{\bf Remark.} Corollary \ref{i.cor3} can be proven algebraically, too.
\smallskip

\noindent{\it Proof of Theorem \ref{i.th1}.} (i) For $1\le j\le n$
and $\zeta\in\Bbb D$ put $\varphi_j(\zeta)=0,$ if $k$ does not
divide $j,$ and $\varphi_j(\zeta)=\binom{n/k}{j/k}\zeta^{j/k},$ if
$k$ divides $j.$ As the zeroes of the polynomial
$(1+\zeta^k)^{n/k}$ lie in $\Bbb D,$ we get
$\varphi=(\varphi_1,\dots,\varphi_n)\in\mathcal O(\Bbb D,\Bbb
G_n).$ Furthermore, $\varphi'(0)=ne_k/k$ and so $$\kappa_{\Bbb
G_n}(0;e_k)\le n/k=\rho_n(e_k).$$ The opposite inequality is
straightforward.
\smallskip

(ii) First note that if $\lambda\in\Bbb T,$ then
$\pi_{\lambda}\in \Aut(\Bbb G_n).$ Furthermore, $\kappa_{\Bbb
G_n}(0;\lambda X)=\kappa_{\Bbb G_n}(0;X).$ These two facts imply that we
can assume that $X_1,X_n>0.$

As
\begin{multline*}\kappa_{\Bbb G_n}(0;X)\ge\kappa_{\Bbb
G_{n-1}}(p_{n,1}(0);p_{n,1}'(0)(X))\\ =\kappa_{\Bbb
G_{n-1}}\left(0; \frac{n-1}{n}X_1e_1+X_ne_{n-1}\right),
\end{multline*}
by induction on $n$ we get $\kappa_{\Bbb
G_n}(0;X)\ge\kappa_{\Bbb G_3}(0;Y),$ where
$$Y=3X_1e_1/n+X_ne_3=Y_1e_1+Y_3e_3.$$ Suppose that
$\kappa_{\Bbb G_n}(0;X)=\rho_n(X).$ Then
$$\rho_n(X)\ge\kappa_{\Bbb G_3}(0;Y)\ge\rho_3(Y)=\rho_n(X)$$ and
consequently $\kappa_{\Bbb G_3}(0;Y)=\rho_3(Y).$ As $\Bbb
G_3$ is a taut domain, there exists an extremal disc for
$\kappa_{\Bbb G_3}(0;Y)$ of the form
$$\varphi(\zeta)=(\zeta\varphi_1(\zeta),\zeta\varphi_2(\zeta),
\zeta\varphi_3(\zeta)),$$ where $\varphi'(0)=Y/\rho_3(Y),$
\begin{equation}\label{i3}
\varphi_1(0)=\frac{Y_1}{3(Y_1+3Y_3)},\ \varphi_2(0)=0,\
\varphi_3(0)=\frac{Y_1}{3(Y_1+3Y_3)}.
\end{equation}
Note that $f_\lambda\circ\varphi\in\mathcal O(\Bbb D,\Bbb
D)$ and $f_\lambda\circ\varphi(0)=0$ for each
$\lambda\in\overline{\Bbb D}.$ For $\lambda\in\overline{\Bbb D}$ and
$\zeta\in\Bbb D$ put
$$g_{\lambda}(\zeta)=\frac{f_\lambda(\varphi(\zeta))}{\zeta}=
\frac{\sum_{j=1}^{3}j\varphi_j(\zeta)
\lambda^{j-1}}{3+2\zeta\varphi_1(\zeta)\lambda+\zeta\varphi_2(\zeta)\lambda^2}.$$
We have $g_{\lambda}\in\mathcal O(\Bbb D,\overline{\Bbb D})$
according to the Schwarz--Pick lemma. By (\ref{i3}) we get
$g_{\pm1}(0)=1$ and so $g_{\pm1}\equiv 1$ by the maximum principle, i.e.
$$\varphi_1(\zeta)\pm2\varphi_2(\zeta)+3\varphi_3(\zeta)=
3\pm2\zeta\varphi_1(\zeta)+\zeta\varphi_2(\zeta).$$ These two equalities imply that
$$\varphi_2(\zeta)\equiv\zeta\varphi_1(\zeta)\hbox{ and }
\varphi_3(\zeta)\equiv 1+\frac{\zeta^2-1}{3}\varphi_1(\zeta).$$

Let $\psi(\zeta)=\varphi_1(\zeta)/3.$ Now by
$g_\lambda\in\mathcal O(\Bbb D,\overline{\Bbb D})$ for
$\lambda\in\Bbb T$ we get $$\left|\frac{\psi(\zeta)+2\lambda
\zeta\psi(\zeta)+\lambda^2(1+(\zeta^2-1)\psi(\zeta))} {1+2\lambda
\zeta\psi(\zeta)+\lambda^2\zeta^2\psi(\zeta)}\right|\le 1$$
$$\Leftrightarrow\left|\frac{\psi(\zeta)(1+\lambda
\zeta)^2+\lambda^2(1-\psi(\zeta))}{\psi(\zeta)(1+\lambda\zeta)^2
+1-\psi(\zeta)}\right|\le 1$$
$$\Leftrightarrow\hbox{Re}(\psi(\zeta)(1-\overline{\psi(\zeta)})((\overline{\lambda}+\zeta)^2
-(1+\lambda \zeta)^2))\le 0.$$ If $\lambda=x+iy,\zeta=ir,r\in\Bbb
R,a=\hbox{Re}(\psi(\zeta))-|\psi(\zeta)|^2,b=\hbox{Im}(\psi(\zeta)),$
then $$y(a(2r-y(r^2+1))+bx(1-r^2))\le 0,\ \forall\ x^2+y^2=1.$$
Then for $x=0$ we get $a\ge 0.$ On the other hand,
leaving $y\to 0^+$ yields $-2ar\ge (1-r^2)|b|.$
Consequently $a=b=0$ if $r>0.$ Then by the uniqueness principle
we get $\psi\equiv0$ or $\psi\equiv1.$ So $X_1=0$ or
$X_n=0$ -- a contradiction.\qed
\smallskip

(iii) Let $\root k\of 1=\{\xi_1,\dots,\xi_k\}.$ For
$z\in\overline{\Bbb G_n}$ and $\lambda\in\overline{\Bbb D}$ such
that the denominator of the first formula below is nonzero, put
$$g_z(\lambda)=\lambda f_\lambda(z)=\frac{\sum_{j=1}^{n}jz_j
\lambda^j}{n+\sum_{j=1}^{n-1}(n-j)z_j\lambda^j},$$
$$g_{z,k}(\lambda)=\frac{\sum_{j=1}^kg_z(\xi_j\lambda)}{k\lambda^k}.$$
The equalities $\sum_{j=1}^k\xi_j^m=0,$ $m=1,\dots,k-1,$ and the
Taylor formula show that $g_{z,k}$ is extended analytically at
$0.$ More precisely, $g_{z,k}(0)=P_k(z),$ where $P_k$ is
polynomial such that $\frac{\partial P_k}{\partial z_k}(0)=k/n$
and $$P_k(tw_1,t^2w_2,\dots,t^nw_n)=t^kP(w),\
t,w_1,w_2,\dots,w_n\in\Bbb C.$$ The maximum principle implies that
$g_{z,k}\in\mathcal O(\Bbb D,\overline{\Bbb D})$. In particular,
$|P_{k}(z)|\le 1.$ To prove the desired inequality $\gamma_{\Bbb
G_n}(0;e_k)>\rho_n(e_k),$ it suffices to show that $|P_k(z)|<1$
for $z\in\overline{\Bbb G_n}.$ Assume the contrary. Then
$P_k(z)=e^{i\theta}$ for some $\theta\in\Bbb R$ and
$z\in\overline{\Bbb G_n}.$ The maximum principle and the
inequality of imply that $g_z(\xi_j\lambda)=
e^{i\theta}\lambda^k,$ $\lambda\in\Bbb T,$ $1\le j\le k.$ For
$\xi_j=1$ we get
$$\sum_{j=1}^{n}jz_j\lambda^j=
e^{i\theta}(n\lambda^{k}+\sum_{j=1}^{n-1}(n-j)z_j\lambda^{k+j}).$$
Comparing the coefficients of the corresponding powers of $\lambda,$
we get $z_k=e^{i\theta}n/k,$
$z_{n+1-k}=\dots=z_{n-1}=0$ and $$(k+j)z_{k+j}=e^{i\theta}(n-j)z_j,\
1\le j\le n-k.$$ These equalities imply that
$z_{kl}=e^{i\theta}\binom{n/k}{l},$ $1\le l\le[n/k].$ On the other hand,
as $k$ does not divide $n,$ $n-k<k[n/k]<n$ and so
$z_{k[n/k]}=0,$ a contradiction.

\setcounter{equation}{0}
\section{Estimates for $\gamma_{\Bbb G_{2n+1}}(0;e_2)$}\label{est}

One of the aims of this section is to measure the quantities in
the inequality of Theorem \ref{i.th1} (iii) in the simplest case.
More precisely, we will find $\gamma_{\Bbb G_{2n+1}}(0;e_2)$ with
an error of $o(n^{-3}).$ To this aim we use that $\gamma_{\Bbb
G_n}(0;e_j)$ solves an extremal problem for a class of
polynomials. This observation, combined with computer checks,
allows us to show that the Carath\'eodory and Kobayashi metrics do
not coincide on $\Bbb G_3,$ thereby sharpening Theorem
\ref{i.th1}. Probably our approach can be applied for obtaining
the same result for $\Bbb G_n,$ $n\ge 4.$

Let $n,k\in\Bbb N,$ $k\le n.$ Note that $$\kappa_{\Bbb
G_n}(0;e_k)\le\kappa_{\Bbb G_{k[n/k]}}(0;e_k)=1/[n/k].$$
Consequently $$k/n\le\gamma_{\Bbb
G_n}(0;e_k)\le\kappa_{\Bbb G_n}(0;e_k)\le 1/[n/k].$$ in particular,
$$\lim_{n\to\infty}n\gamma_{\Bbb
G_n}(0;e_k)=\lim_{n\to\infty}n\kappa_{\Bbb G_n}(0;e_k)=k.$$

Let now $n\ge 3$ be an odd number. Then $2/n<\gamma_{\Bbb
G_n}(0;e_2)$ by Theorem \ref{i.th1} (iii). On the other hand,
$$\gamma_{\Bbb G_n}(0;e_2)\le\kappa_{\Bbb
G_n}(0;e_2)\le\frac{2}{n-1}.$$

We will later improve both estimates. For the upper estimate we
need the following. Let $D\subset\Bbb C^n$ be a
$(k_1,\dots,k_n)$-balanced domain. Denote by $\mathcal P_j$ the
set of polynomials $P$ such that $\sup_D|P|\le 1$ and
$P\circ\pi_\lambda=\lambda^{k_j}P,$ $\lambda\in\Bbb C.$ Put
$\mathcal L_j=\span (e_j,\dots,e_l),$ where $l\ge j$ is the
greatest index such that $k_l=k_j.$ The proof of Theorem
\ref{i.th1} (iii) easily implies that

\begin{proposition}\label{e.pr1} If $D\subset\Bbb C^n$ is
$(k_1,\dots,k_n)$-balanced domain and $X\in\mathcal L_j,$ $1\le
j\le n$, then $\gamma_D(0;X)=\sup\{|P'(0)X|:P\in\mathcal P_j\}.$
\end{proposition}

\noindent{\bf Remarks.} a) This proposition directly implies
that if $D$ is balanced domain, then
$$\gamma_D(0;X)=\sup\{|L'(0)X|:\sup_D|L|\le 1,\ L\hbox{ -- linear
function}\}$$ and so $\gamma_D(0;\cdot)=\gamma_{\hat D}(0;\cdot)$
(for the last one see also Proposition \ref{bala.pr1} (i)).

b) Another corollary is the formula
\begin{equation}\label{e1}
\gamma^{-1}_{\Bbb G_n}(0;e_2)=\inf_{c\in\Bbb C}\max_{z\in\partial
G_n}|z_2+cz_1^2|.
\end{equation}
In spite of this formula, $\gamma_{\Bbb G_{2n+1}}(0;e_2)$
is hard to calculate (see Lemma \ref{e.l4} for $n=1$).
\smallskip

c) If $n$ is an even number, the extremal polynomials for
$\gamma_{\Bbb G_n}(0;e_2)=2/n$ can differ not only by a constant
of an absolute value 1. For example, after some easy calculations,
by the proof of Theorem \ref {i.th1} we get the
polynomial $2z_2/n-(n-1)z_1^2/n^2,$ but the polynomial $(2z_2-z_1^2)/n$
is also extremal.
\smallskip

\begin{proposition}\label{e.pr2} If $n\ge 3$ is an odd number, then
$$\frac{2}{n}\left(1+\frac{2}{(n-1)(n+2)}\right)<\gamma_{\Bbb
G_n}(0;e_2)<\frac{2}{n}\left(1+\frac{2}{(n-1)(n+1)}\right).$$
\end{proposition}

\beginproof {\it Lower estimate.} Let us first see that the polynomial
$$P_n(z)=\frac{n-1}{2(n+1)}z_1^2-z_2$$ satisfies the equality
$$\max_{\partial \Bbb G_n}|P_n|=M_n:=\frac{(n-1)(n+2)}{2(n+1)}.$$
So if $$ g_n(t)=\frac{1}{2}\sum_{j=1}^n
t_j^2-\frac{1}{n+1}(\sum_{j=1}^n t_j)^2,\ t\in\Bbb C^n,$$ then
$\max_{\Bbb T^n}|g_n|=M_n.$ To prove the last one, let
$M_n^\ast=\max_{\Bbb T^n}|g_n|.$ As
$g_n(e^{i\theta}t)=e^{2i\theta}g_n(t)$ for each $\theta\in\Bbb R,$
$t\in\Bbb C^n,$ there exists point $u\in\Bbb T^n$ such that
$g_n(u)=M_n^\ast.$ Putting $u_j=x_j+iy_j,$ $x_j,y_j\in\Bbb R,$
$1\le j\le n,$ we get
$$M_n^\ast=\mbox{Re}(g_n(u))=\frac{1}{2}\sum_{j=1}^n
(x_j^2-y_j^2)+\frac{1}{n+1}((\sum_{j=1}^ny_j)^2-(\sum_{j=1}^nx_j)^2)$$
$$\le\frac{1}{2}\sum_{j=1}^n(x_j^2-y_j^2)+\frac{1}{n+1}(n\sum_{j=1}^ny_j^2-
(\sum_{j=1}^nx_j)^2)$$
$$=\frac{(n-1)n}{2(n+1)}+\frac{1}{n+1}(\sum_{j=1}^n
x_j^2-(\sum_{j=1}^nx_j)^2)$$ by the Cauchy--Schwarz inequality and
the equalities $y_1^2=1-x_1^2,\dots,y_n^2=1-x_n^2.$ The last
expression is a linear function for each $x_j.$ Consequently it is
maximal for $1$ and/or $-1.$ As $n$ is odd,
$$M_n^\ast=\frac{(n-1)n}{2(n+1)}+\frac{n-1}{n+1}=M_n,$$ with maximum
attained only if $[n/2]$ or $[n/2]+1$ among the numbers $t_j$ are
equal to some $t_0\in\Bbb T,$ and the rest to $-t_0.$

Using this fact one can easily prove that if $\varepsilon>0$
is sufficiently small and
$$g_{n,\varepsilon}(t)=g_n(t)+\varepsilon\sum_{j=1}^n t_j^2-\varepsilon(n+1)(\sum_{j=1}^n t_j)^2,\ t\in\Bbb C^n,$$ then
$\max_{\Bbb T^n}|g_{n,\varepsilon}|<M_n.$ Consequently
for $$P_{n,\varepsilon}=\frac{n-1-2n(n+1)\varepsilon}{2(n+1)}z_1-(1+2\varepsilon)
z_2$$ one has the inequality $\max_{\partial \Bbb
G_n}|P_{n,\varepsilon}|<M_n,$ showing that $$\gamma_{\Bbb
G_n}(0;e_2)>\frac{1}{M_n}=\frac{2}{n}\left(1+\frac{2}{(n-1)(n+2)}\right).$$

{\it Upper estimate.} By (\ref{e1}) we need to prove that
if $c\in\Bbb C,$ then $$m_{n,c}:=\max_{z\in\partial
G_n}|z_2+cz_1^2|>\frac{n(n^2-1)}{2(n^2+1)}.$$ The coefficients
of the polynomials $(t-1)^n$ and $(t-1)(t^2-1)^{\frac{n-1}{2}}$
give two points $z\in\partial\Bbb G_n$ such that $z_1=n$,
$z_2=\frac{n(n-1)}{2}$ and $z_1=1,z_2=\frac{1-n}{2}$, respectively. Then
$$2m_{n,c}\ge\max\{|n-1-2c|,|n(n-1)+2cn^2|\}$$ and consequently
$$2(n^2+1)m_{n,c}\ge|n^2(n-1)-2cn^2|+|n(n-1)+2cn^2|$$ $$\ge
n^2(n-1)+n(n-1)=n(n^2-1).$$ This means that
$m_{n,c}\ge\frac{n(n^2-1)}{2(n^2+1)}.$ Suppose that we have an
equality. Then $ c=-\frac{(n-1)^2}{2(n^2+1)}.$ On the other hand,
the coefficients of the polynomial $(t-i)(t-1)^{n-1}$ give a point
$z\in\partial\Bbb G_n$ such that
$z_1=n-1+i,z_2=\frac{(n-1)(n-2)}{2}+(n-1)i$ and
$\left|z_2-\frac{(n-1)^2}{2(n^2+1)}z_1^2\right|>\frac{n(n^2-1)}{2(n^2+1)}$,
a contradiction.\qed\smallskip

Let now $D$ be a domain in $\Bbb C^n,$ $z\in D$ and $k\in\Bbb N.$
Denote by $\hat\gamma_D^{(k)}(z;X)$ the largest pseudonorm
not exceeding the $k$-th Carath\'eodory pseudometric
$$\gamma_D^{(k)}(z;X)=
\sup\{|f^{(k)}_z(X)|:f\in\O(D,\D),\ord_z f\ge k\},$$ where
$$f^{(k)}_z(X)=\sum_{|\alpha|=k}\frac{D^\alpha
f(z)X^\alpha}{\alpha!}.$$ (One can define similarly the $k$-th Kobayashi pseudometric,
which is essentially different from
the Kobayashi pseudometric $\kappa_D^{(k)}$ of order $k$, defined in Section \ref{lem-der}).

As $\gamma_D(z;\cdot)$ is a pseudonorm,
$$\gamma_D\le\hat\gamma_D^{(k)}\le\hat\kappa_D.$$ Also note that,
since the family $O(\Bbb G_3,\Bbb D)$ is normal,
the argument in the proof of Theorem \ref{kb.th} shows the existence
of $m\le 2n-1$ $\Bbb R$-linearly independent vectors
$X_1,\dots,X_m\in\Bbb C^n$ of sum $X$, so that
$$\hat\gamma^{(k)}_D(z;X)=\sum_{j=1}^m \gamma^{(k)}_D(z;X).$$

\begin{theorem}\label{e.th3}
$\hat\gamma^{(2)}_{\Bbb G_3}(0;e_2)>\gamma_{\Bbb G_3}(0;e_2).$
In particular, $\hat\kappa_{\Bbb G_3}(0;e_2)>\gamma_{\Bbb G_3}(0;e_2)$
and consequently $k_{\Bbb G_3}(0,\cdot)\neq c_{\Bbb G_3}(0,\cdot).$
\end{theorem}

This theorem follows from the two lemmas below.

\begin{lemma}\label{e.l4}$\gamma_{\Bbb G_3}(0;e_2)\le C_0:=\sqrt\frac{8}{13\sqrt
{13}-35}=0,8208\dots.$
\end{lemma}

\begin{lemma}\label{e.l5} $\hat\gamma_{\Bbb G_3}^{(2)}(0;e_2)\ge
C_1=\sqrt{0,675}=0,8215\dots.$
\end{lemma}

\noindent{\it Proof of Lemma \ref{e.l4}.} By
(\ref{e1}) we need to show that if $c\in\Bbb C,$ then
$$\max_{z\in\partial\Bbb G_3}|z_2-cz_1^2|^2\ge C_0^{-2}.$$
It suffices to show this inequality for $c\in\Bbb R.$
Indeed, for each $z\in\partial\Bbb G_3$ we have $\overline
z\in\partial\Bbb G_3$ and so $$2\max_{z\in\partial\Bbb
G_3}|z_2-cz_1^2|\ge\max_{z\in\partial\Bbb G_3}(|z_2-cz_1^2|+
|\overline z_2-c\overline z_1^2|)$$ $$\ge\max_{z\in\partial\Bbb
G_3}|2z_2-(c+\overline c)z_1^2|=2\max_{z\in\partial\Bbb
G_3}|z_2-\mbox{Re}(c)z_1^2|.$$

Let now $c\in\Bbb R.$ Then $$\max_{z\in\partial\Bbb
G_3}|z_2-cz_1^2|^2\ge\max_{\varphi\in[0,2\pi)}|1+2e^{i\varphi}-c(2+e^{i\varphi})^2|^2$$
$$=\max_{\varphi\in[0,2\pi)}(4c(4c-1)\cos^2\varphi+4(10c^2-7c+1)\cos\varphi+25c^2-22c+5).$$
Put $$f_c(x)=4c(4c-1)x^2+4(2c-1)(5c-1)x+25c^2-22c+5,\
x\in[-1,1].$$ If $c\not\in\Delta=\left(\frac{1}{6},\frac{5-\sqrt{17}}{4}\right),$ then
$$\max_{x\in[-1,1]}f_c(x)=\max\{f_c(-1),f_c(1)\}\ge
\left(\frac{9-\sqrt{17}}{4}\right)^2>\frac{1}{C_0^2}.$$ Otherwise
$$\max_{x\in[-1,1]}f_c(x)=f_c\left(\frac{10c^2-7c+1}{2c(1-4c)}\right)
=\frac{(3c-1)^3}{c(4c-1)}=:g(c)$$ and it remains to see that
$\min_{c\in\Delta}g(c)=g\left(\frac{\sqrt{13}-1}{12}\right)=\frac{1}{C_0^2}.$\qed\smallskip

\noindent{\bf Remark.} Let $ c_0=\frac{\sqrt{13}-1}{12}$ and
$M=\max_{z\in\partial\Bbb G_3}|z_2-c_0z_1^2|.$ As in the
proof of Proposition \ref{e.pr2} we have $$M=\max_{z\in\partial\Bbb
G_3}\hbox{Re}(z_2-c_0z_1^2)=\max_{\alpha,\beta,\gamma\in\Bbb R}
h(\alpha,\beta,\gamma),$$ where
$$h(\alpha,\beta,\gamma)=(1-2c_0)(\cos(\alpha+\beta)+\cos(\beta+\gamma)
+\cos(\gamma+\alpha))$$
$$-c_0(\cos2\alpha+\cos2\beta+cos2\gamma).$$ The computer calculations
show that the critical points of $h$ (up to a permutation of the
variables) are of the form ($k\pi,l\pi,m\pi)$ or
$(\pm\alpha_0+j\pi/2+2k\pi,\pm\alpha_0+j\pi/2+2l\pi,\pm\gamma_0+j\pi/2+2m\pi),$
$k,l,m\in\Bbb Z,$ $j=0,1,2,3.$ Then the proof of Lemma \ref{e.l4}
implies that $ M=C_0^{-1},$ i.e. $\gamma_{\Bbb G_3}(0;e_2)=C_0.$
\smallskip

\noindent{\it Proof of Lemma \ref{e.l5}.} Let
$$f(z)=0,675z_2^2-0,291z_2z_1^2+0,033z_1^4.$$ We will check first that
$\max_{z\in\partial\Bbb G_3}|f(z)|<1$ by reducing the check to
finitely many points, and then using a computer program for them.
Put $\theta=(\theta_1,\theta_2),\ \theta_1,\theta_2\in[0,2\pi),$
$$g_1(\theta)=1+e^{i\theta_1}+e^{i\theta_2},\
g_2(\theta)=e^{i(\theta_1+\theta_2)}+e^{i\theta_1}+e^{i\theta_2},$$
$$g(\theta)=0,675g_2^2(\theta)-0,291
g_2(\theta)g_1^2(\theta)+0,033g_1^4(\theta).$$ We have to prove that
$\max|g(\theta)|<1.$ Let
$$d(\theta,\tilde\theta)=\max\{|\theta_1-\tilde\theta_1|,
|\theta_2-\tilde\theta_2|\}.$$ As
$|e^{i\theta_j}-e^{i\tilde\theta_j}|\le|\theta_j-\tilde\theta_j|,$
$j=1,2,$ we get
$$|g_1(\theta)-g_1(\tilde\theta)|\le2d(\theta,\tilde\theta),\
|g_2(\theta)-g_2(\tilde\theta)|\le4d(\theta,\tilde\theta).$$
Then the inequalities $|g_1|\le 3,|g_2|\le 3$ imply
$$|g(\theta)-g(\tilde\theta)|\le(0,675\cdot 24+0,291\cdot
72+0,033\cdot 216)
d(\theta,\tilde\theta)=44,28d(\theta,\tilde\theta).$$

Let now $\theta_1,\theta_2$ vary in the interval
$[0;6,2832]\supset[0,2\pi]$ with a step of $4\cdot 10^{-5}.$
The results of the corresponding computer program (see Appendix C)
show that $|g(\theta)|\le 0,999$ for the variable $\theta=(\theta_1,\theta_2).$
(In fact these results lead to the
hypothesis that $\max |g(\theta)|=0,999,$ with a maximum attained at
the points $(0,\pi),$ $(\pi,0)$ and $(\pi,\pi).$) Then by
the inequalities $|g(\theta)-g(\tilde\theta)|\le
44,28d(\theta,\tilde\theta)$ and $\frac{2}{44,28}\cdot
10^{-3}>4\cdot 10^{-5}$ we easily get $\max
|g(\theta)|<1.$

Of the above it follows that if $X\in\span(e_1,e_3),$ then
$$\gamma^{(2)}_{\Bbb G_3}(0;e_2+X)\ge|f^{(2)}_0(e_2+X)|/2= |f^{(2)}_0(e_2)|/2=C_1.$$

On the other hand, recall that there exist five vectors
$X_1,\dots,X_5\in\Bbb C^3$ (some of them can be zero vectors)
of sum $e_2$ and such that $\hat\gamma^{(2)}_{\Bbb
G_3}(0;e_2)=\sum_{j=1}^5\gamma^{(2)}_{\Bbb G_3}(0;X_j).$ As
$\gamma^{(2)}_{\Bbb G_3}(0;X_j)\ge C_1|\langle X_j,e_2\rangle|,$
we get $\hat\gamma^{(2)}_{\Bbb G_3}(0;e_2)\ge C_1.$\qed\smallskip

\noindent{\bf Remark.} An important moment in the above proof is
finding a polynomial of the form $f(z)=az_2^2+bz_2z_1^2+cz_1^4$
such that $\max_{\partial\Bbb G_3}|f|\le 1$ and $\sqrt a>C_0.$
Computer experiments show that the maximal value of $a$ is
$0,676\dots,$ i.e. very close to the number $0,675/0,999.$
\smallskip

Finally let us note that $\gamma^{(2)}_{\Bbb G_n}(0;\cdot)$ is not a norm.

\begin{proposition}\label{e.pr6} If $X_1,X_n\in\Bbb C,$ then
$$\gamma^{(2)}_{\Bbb
G_n}(0;X_1e_1+X_ne_n)\ge\sqrt{\frac{n+1}{2}\gamma_{\Bbb
G_n}(0;e_2)|X_1X_n|}.$$ In particular, as $\gamma_{\Bbb
G_3}(0;e_2)>2/3$ and $\gamma_{\Bbb G_n}(0;e_n)\ge2/n,$ it follows that
$$\gamma^{(2)}_{\Bbb G_n}(0;ne_1+e_n)>2=\hat\kappa_{\Bbb
G_n}(0;ne_1+e_n)=\gamma^{(2)}_{\Bbb
G_n}(0;ne_1)+\gamma^{(2)}_{\Bbb G_n}(0;e_n),\ n\ge 3.$$
\end{proposition}

\beginproof
Let $t_1,\dots,t_n\in\Bbb D.$ Consider
$\sum_{k=1}^n\frac{t_k^{n+1}}{n}$ as a function $f$ of
$\sigma(t).$ Then $f\in\mathcal O(\Bbb G_n,\Bbb D),$
$\mbox{ord}_0f=2$, and the Waring formula (see e.g. \cite{Waer}))
implies that the coefficient of $z_1z_n$ equals
$(-1)^{n-1}\frac{n+1}{n}.$ So $$\gamma^{(2)}_{\Bbb
G_n}(0;X_1e_1+X_ne_n)\ge\left|\frac{f_0^{(2)}(X)}{2}\right|
=\sqrt{\frac{n+1}{n}\gamma_{\Bbb G_n}(0;e_2)|X_1X_n|},$$ where
$X=X_1e_1+X_ne_n.$ As $\gamma_{\Bbb G_n}(0;e_2)=2/n$ for an even
number $n,$ we get the proposition for such an $n.$

On the other hand, Proposition \ref{e.pr1} implies that if
$2C_n:=\gamma_{\Bbb G_n}(0;e_2),$ then there exists a $c_n$ such that
$P(z)=2C_nz_2-c_nz_1^2$ is an extremal function for $\gamma_{\Bbb
G_n}(0;e_2).$ For an odd $n=2k-1$ we replace $t_1,\dots,t_n$
by $t_1^k,\dots,t_n^k.$ Thus we get the function
$$(C_n-c_n)\left(\sum_{j=1}^nt_j^k\right)^2
-C_n\sum_{j=1}^nt_j^{2k}.$$ Consider this function as a function
$g$ of $\sigma(t).$ Then $g\in\mathcal O(\Bbb G_n,\Bbb D),$
$\mbox{ord}_0 g=2,$ and the coefficient of $z_1z_n$ equals $-(n+1)C_n.$
Consequently
$$\gamma^{(2)}_{\Bbb G_{n}}(0;X)\ge\left|\frac{g_0^{(2)}(X)}{2}\right|=\sqrt{\frac{n+1}{2}
\gamma_{\Bbb G_n}(0;e_2)|X_1X_n|}.\qed$$

\noindent{\bf Appendix C.}
\smallskip
\begin{small}
\begin{verbatim}
      Program niki
      implicit real*8 (a-h,o-z)
      implicit integer*4 (i-n)
      complex*16 g0,g1,g2

      data c1,c2,c3 /0.675D0, -0.291D0, 0.033D0/
      data e,o,t1d,t2d,t1u,t2u/1.0D0,3*0.0D0,2*6.2832D0/
      write(*,102)
  200  continue
      read(*,*,ERR=201,END=201) s
      N1=(t1u-t1d)/s
      N2=(t2u-t2d)/s
      gu=-1D30
      do i1=0,N1
         t1=t1d+FLOAT(i1)*s
         do i2=0,N2
            t2=t2d+FLOAT(i2)*s
            g0=DCMPLX(DCOS(t1)+DCOS(t2),DSIN(t1)+DSIN(t2))
            g1=g0+DCMPLX(e,o)
            g2=g0+DCMPLX(DCOS(t1+t2),DSIN(t1+t2))
            g = CDABS(c1*g2**2+c2*g2*g1**2+c3*g1**4)
            if (g.GT.gu) then
               gu=g
               t1g=t1
               t2g=t2
            endif
         enddo
      enddo
      write(*,100) s,gu,t1g,t2g
      goto 200
  201  continue
      write(*,101)
      stop
  100 format(1x,2f20.15,2f15.10)
  101 format(8x,' step ',15x,' g-max',9x,'tita-1',9x,'tita-2')
  102 format(8x,' step ' )
      end
\end{verbatim}
\end{small}\begin{footnotesize}\begin{verbatim}
----------------- Results ------------------------------------------
 step                g-max               tita-1         tita-2

 0.001000000000000   0.998999998608272   3.1420000000   3.1420000000
 0.000400000000000   0.998999999688699   3.1414000000   3.1414000000
 0.000100000000000   0.998999999999547   3.1416000000   3.1416000000
 0.000040000000000   0.998999999999547   3.1416000000   3.1416000000
 0.000010000000000   0.998999999999941   3.1415900000   3.1415900000
 0.000004000000000   0.998999999999985   3.1415940000   3.1415940000
 0.000001000000000   0.998999999999999   3.1415930000   3.1415930000
 0.000000400000000   0.999000000000000   3.1415928000   3.1415928000
 0.000000100000000   0.999000000000000   3.1415927000   3.1415927000
 0.000000040000000   0.999000000000000   3.1415925600   3.1415925600
\end{verbatim}\end{footnotesize}

\setcounter{equation}{0}
\section{Continuity of $l_{\Omega_n}(A,\cdot)$}\label{cont}

As mentioned in the Introduction, the continuous dependance of
SNPP on the data (a necessary condition for reduction to an
analogous problem on $\G_n$) is linked with the continuity of the
function $l_{\Om_n}.$ The aim of this section is to describe all
matrices $A\in\Om_n$ such that $l_{\Om_n}(A,\cdot)$ is a
continuous function.

First recall that
$$l_{\Omega_n}(A,B)\ge l_{\Bbb
G_n}(\sigma(A),\sigma(B)).$$ Furthermore, if $A,B\in\CY_n$ (i.e. they are
cyclic matrices in $\Om_n$), then we have equality (see \ref{np2})) and
so $l_{\Omega_n}$ is a continuous function on the open
set $\CY_n\times\CY_n.$ In general,
we have equality if and only if $l_{\Omega_n}$ is a continuous function
in $(A,B).$ To this aim it suffices to use that
$l_{\Bbb G_n}$ is a continuous function and the set
$\CY_n\times\CY_n$ (where we have equality) is dense in
$\Om_n\times\Om_n.$

In \cite{Tho-Tra2} the authors consider matrices $B\in\Omega_n$ such that
$l_{\Omega_n}(A,.)$ is a continuous function at $B$ for each
$A\in\Om_n$. They hypothesize that this is true for each
$B\in \mathcal C_n$ and confirm this for $n\le 3$
\cite[Proposition 1.4]{Tho-Tra2}. Using the results from
Section \ref{zeko} for the continuity of $\kappa_{\Omega_n}(A,.)$,
the converse proposition is proven for each
$n$ (see \cite[Theorem 1.3]{Tho-Tra2}).

We first prove the following

\begin{proposition}\label{np.pr4} If $\lambda\in\Bbb D$ and $A\in\mathcal C_n,$
then the following are equivalent:

(i) the eigenvalues of $A$ are equal;

(ii) $l_{\Omega_n}$ is a continuous at the point $(A,\lambda I_n);$

(iii) $l_{\Omega_n}(A,\cdot)$ is continuous at the point $\lambda I_n.$
\end{proposition}

\beginproof The implication (ii)$\Rightarrow$(iii) is trivial.
For the rest of the proof we may assume that $\lambda=0,$
applying $\Phi_\lambda$ (see (\ref{fi})).

We will now show that (i)$\Rightarrow$(ii). Let the eigenvalues of $A$
be equal to $a.$ If $A_j\to A$ and $B_j\to 0,$ then
$$l_{\Omega_n}(A_j,B_j)\ge c^\ast_{\Bbb G_n}(\sigma(A_j),\sigma(B_j))\to c^\ast_{\Bbb
G_n}(\sigma(A),0)=|a|=l_{\Omega_n}(A,0).$$ (the last two
equalities follow from Proposition \ref{i.pr4} (iii) and
\ref{np4}), respectively. So the function $l_{\Omega_n}$ is lower
semicontinuous at the point $(0,B).$ As it is (always) upper
semicontinuous, it is continuous at this point.

It remains to prove that (iii)$\Rightarrow$(i). As $\mathcal
C_n$ is a dense subset of $\Omega_n,$ we can find a sequence
$\mathcal C_n\supset(B_j)\to 0.$ By (\ref{np4}) and (\ref{np2})
we get $$r(A)=l_{\Om_n}(A,0)\leftarrow
l_{\Omega_n}(A,B_j)=l_{\G_n}(\spe(A),\spe(B_j))\to
l_{\G_n}(\spe(A),0).$$ Proposition \ref{i.pr4} (iii)
implies that the eigenvalues of $A$ are equal.\qed
\smallskip

Unlike the above proposition, (\ref{np4}) implies that $l_{\Om_n}(A,\cdot)$
is a continuous function for each scalar
matrix $A\in\Om_n$

As we noted, if $A\in\Om_n$ ($n\ge 2$), then the following are equivalent:

(i) the function $l_{\Om_n}$ is continuous at $(A,B)$ for each $B\in\Om_n;$

(ii) $l_{\Omega_n}(A,\cdot)=l_{\mathbb{G}_n}(\sigma(A),\sigma(\cdot)).$

\noindent Also consider the condition

(iii) $A\in\mathcal C_2$ has (two) equal eigenvalues.

\noindent By \cite[Theorem 8]{Cos4}, (iii) implies (ii).
Theorem \ref{septh} says that the scalar matrices and those satisfying (iii)
are the only ones, for which $l_{\Om_n}(A,\cdot)$
is a continuous function. Then Proposition \ref{np.pr4} implies that (iii)
follows from (i). So the assertions (i), (ii) and (iii)
are equivalent.

\begin{theorem}\label{septh} If $A\in\Om_n,$ then $l_{\Om_n}(A,\cdot)$
is a continuous function if and only if $A$ is a scalar
matrix or $A\in\mathcal C_2$ has two equal eigenvalues.
\end{theorem}

\beginproof Applying $\Phi_\lambda$ and \begin{equation}\label{gi}
\Psi_P(X)=P^{-1}XP, \quad P\in\mathcal M_n^{-1}, X\in\mathcal M_n,
\end{equation}
we can assume that $0$ is an eigenvalue of $A$ with a maximal number of
Jordan blocks and the matrix is in Jordan form.

It suffices to prove that $l_{\Om_n}(A,\cdot)$ is not a continuous
function, if $A$ has a nonzero eigenvalue or $A\in\Om_n$ is a nonzero
nilpotent matrix and $n\ge 3$.

In the first case let $d_1\ge\dots\ge d_k$ be the number of Jordan
blocks that correspond to the different eigenvalues
$\lambda_1=0,\lambda_2,\dots,\lambda_k.$ We will prove that the
function $l_{\Om_n}(A,\cdot)$ is not continuous at $0.$ It is
easily seen that $A$ can be expressed as blocks $A_1,\dots A_l$
(of dimensions $n_1,\dots,n_l$) so that the eigenvalues of $A_1$
are equal to $0$ and the remaining blocks are cyclic with at least
two different eigenvalues ($A_1$ is missing if $d_1=d_2$). By
Proposition \ref{i.pr4} (iii), there exists a sequence of matrices
$(A_{i,j})_j\to 0,$ $1\le i\le l$ such that
$\sup_{i,j}l_{\Om_{n_i}}(A_i,A_{i,j}):=m<r(A).$ Forming $A_j$ from
the blocks $A_{1,j},\dots,A_{l,j},$ it follows that
$l_{\Om_n}(A,A_j)\le\max_i l_{\Om_{n_i}}l(A_i,A_{i,j})\le
m<l_{\Om_n}(A,0),$ meaning that $l_{\Om_n}(A,\cdot)$ is not a
continuous function at $0.$

Let now $A\neq 0$ be a nilpotent matrix. Then $A= (a_{ij})_{1\le
i,j\le n},$ where $a_{ij}=0$ for $j\neq i+1$. Let
$r=\mbox{rank}(A)\ge 1$. Following the proof of \cite[Proposition
4.1]{Tho-Tra2}, let
$$F_0=\{1\} \cup \left\{ j \in \{2,\dots,n\} : a_{j-1,j}=0 \right\}
:=\{ 1=b_1 < b_2 < \dots < b_{n-r} \},$$
and $b_{n-r+1}=n+1$. Put $d_i = 1+ \# \left( F_0 \cap
\{(n-i+2),\dots,n\}\right) $. As $A\neq 0$ is a nilpotent matrix,
it has a Jordan form such that
$a_{n-1,n}=1$ and $1=d_1=d_2\le d_3\le\dots\le d_n= \#F_0 = n-r$,
$d_{j+1}\le d_j+1$.

In \cite[Proposition 4.1, Corollary 4.3]{Tho-Tra2} there is a
necessary and sufficient condition for lifting of discs from $\O(\Bbb D,\G_n)$
to ones from $\O(\Bbb D,\Om_n),$ passing through a
cyclic and nilpotent matrix. They easily imply that for each $C\in\mathcal C_n$,
$$l_{\Om_n}(A,C)=h_{\G_n}(0,\sigma(C)):=\inf\{|\alpha|:\exists
\psi\in\HH(\D,\G_n):\psi(\alpha)=\sigma(C)\},$$ where
$$\HH(\D,\G_n)=\{\psi\in\O(\D,\G_n):\mbox{ord}_0\psi_j\ge d_j,
\ 1\le j\le n\}.$$

Note that $d_j\le j-1$ for $j\ge 2.$ Let $m=\min_{j\ge2}\frac{d_j}{j-1}$
and choose $k$ such that $\frac{d_k}{k-1}=m$.
If $m=1$, then ${d_j}= {j-1}$ for each $j\ge 2$ and so in this case
for $n\ge 3$ one can take $k=3$.

Let $\lambda$ be a sufficiently small positive number,
$b=k\lambda^{k-1}$ and $c=(k-1)\lambda^k.$ Then $\lambda$ is a
double zero of the polynomial $\Lambda(z) = z^{n-k}(z^k-bz+c)$ of
zeroes in $\mathbb D.$ Let $B$ be a diagonal matrix with a
characteristic polynomial $P_B(z)=\Lambda(z)$.

Suppose that the function $l_{\Om_n}(A,\cdot)$ is continuous at $B.$ Then
$$l_{\Om_n}(A,B)=h_{\G_n}(0,\sigma(B))=:\alpha.$$

\begin{lemma}
\label{sumder} If $l_{\Om_n}(A,B)=\alpha$, then there exists
$\psi\in\HH(\D,\G_n)$ so that $\psi(\alpha)=\sigma(B)$ and
$$\sum_{j=1}^n\psi'_j(\alpha)(-\lambda)^{n-j}=0.$$
\end{lemma}

\beginproof Similarly to the proof of \cite[Proposition 4.1]{Tho-Tra2},
let $\varphi \in \O (\D, \Om_n)$ and $\tilde\alpha\in\D$ so that
$\varphi(0)=A$ and $\varphi(\tilde\alpha)=B$. By \cite[Corollary
4.3]{Tho-Tra2} we have $\tilde\psi=\sigma \circ \varphi \in
\HH(\D,\G_n)$.

Now let us examine
$\sigma_n(\varphi(\zeta))-\sigma_n(B)=\sigma_n(\varphi(\zeta))$
near $\zeta =\tilde\alpha$. We can assume that the first two
diagonal elements of $B$ are equal to $\lambda$. If
$\varphi_\lambda (\zeta)=\varphi(\zeta)-\lambda I_n$, then the
first two columns of $\varphi_\lambda (\alpha)$ are zero.
Consequently $\sigma_n \circ \varphi_\lambda =
\det\varphi_\lambda$ has a zero of order at least $2$ at $\alpha$.
On the other hand, $\G_n$ is a taut domain, which easily provides
the required $\psi.$\qed

\begin{lemma}
\label{alphsm} We have $\alpha^m\lesssim\lambda$. Furthermore,
if $m=1$ and $n\ge 3$, then $\alpha^{2/3}\lesssim\lambda$. In particular,
always $\alpha\ll\lambda$.
\end{lemma}

\beginproof Note that there exists $\eps>0$ so that for $\lambda<\eps$ the mapping
$\zeta\to(0,\dots,0,k(\eps\zeta)^{d_k},(k-1)\lambda(\eps\zeta)^{d_k},0,\dots,0)$
is a competitor for $h_{\Om_n}(A,B).$ Consequently
$(\eps\alpha)^{d_k}\le\lambda^{k-1}$, i.e.
$\alpha^m\lesssim\lambda$.

If $m=1$ and $n\ge k=3,$ by considering the mapping
$\zeta\to(0,3\lambda^{1/2}\eps\zeta,2(\eps\zeta)^2,0\dots,0)$
we get $(\eps\alpha)^2\le\lambda^3$.\qed
\smallskip

To finish the proof of the theorem, put
$\psi_j(\zeta)=\zeta^{d_j}\theta(\zeta)$; the condition in Lemma
\ref{sumder} becomes
\begin{equation}
\label{condS} a\frac{(-\lambda)^n}{\alpha}+S=0,
\end{equation}
where $a=(k-1)d_k-kd_{k-1}$ and
$S=\sum_{j=1}^n\alpha^{d_j}\theta'_j(\alpha)(-\lambda)^{n-j}.$
Note that $a\neq 0.$ Indeed, if $m<1,$ then $d_k=d_{k-1}$ and
consequently $a=-d_k,$ while if $m=1$, then $a=(k-1)(k-1) -k (k-2)
= 1.$ As $\mathbb G_n$ is a bounded domain, the Cauchy
inequalities imply $|\theta'_j(\alpha)|\lesssim 1.$

By \ref{alphsm} and by the choice of $k$ it follows that for each $j$,
$$
\alpha^{d_j}\lesssim \lambda^{(k-1)d_j/d_k}  \le \lambda^{j-1}
\le\lambda^{n-1}.$$ So $S\lesssim\lambda^{n-1}.$ Once again by
Lemma \ref{alphsm}, $\alpha\ll\lambda,$ contradicting (\ref{condS}).\qed

\setcounter{equation}{0}
\section{Zeroes of $\kappa_{\Omega_n}$}\label{zeko}

Recall that the spectral Carath\'eodory--Fej\'er problem of order
1 (SCFP) is reduced to the calculation of the Kobayashi metric
$\kappa_{\Omega_n}$ of $\Om_n.$ Furthermore, if $A\in\mathcal C_n$
(i.e. $A$ is a cyclic matrix), then (see (\ref{np5}))
$$\kappa_{\Omega_n}(A;B)=\kappa_{\Bbb G_n}(A;\sigma'_A(B)).$$
In particular,
$\kappa_{\Omega_n}(A;B)=0\Leftrightarrow\sigma'_A(B)=0.$

On the other hand, by Proposition \ref{cyc.equiv} $\sigma'_A(B)=0$
exactly when there exists $Y\in\mathcal M_n$ so that $B=[Y,A].$
Consequently, if $A\in\mathcal C_n$ and $\sigma'_A(B)=0,$
considering $\zeta\to e^{\zeta Y}Ae^{-\zeta Y},$ we find even an
entire curve $\varphi:\Bbb C\to \Omega_n$ so that  $\varphi(0)=A$
and $\varphi'(0)=B.$ In general, if $\kappa_{\Omega_n}(A;B)=0$ (we
do not assume $A\in\mathcal C_n$), then SCFP for an arbitrary disc
instead of the unit one has a solution. Therefore it is important
to know the zeroes of $\kappa_{\Omega_n}.$ This also bears
information for the discontinuity of this function (hence also of
SCFP).

Recall that for the Carath\'eodory metric of $\Omega_n$,
things are much simpler (see Proposition \ref{np.pr5}):
$$\gamma_{\Omega_n}(A;B)=\gamma_{\Bbb G_n}(\sigma(A);\sigma'_A(B))$$
and so $\gamma_{\Omega_n}(A;B)=0\Leftrightarrow\sigma'_A(B)=0.$

To formulate the results in this section we need to introduce some notions.

For $A\in\Omega_n$ denote by $C_A$ the tangent cone (see
\cite[p. 79]{Chi} for this notion) to the isospectral (analytic) set
$$L_A=\{C\in\Omega_n:\spe(C)=\spe(A)\},$$
i.e.
$$C_A=\{B\in\mathcal M_n:\exists 0<c_j\to 0, C_j\in L_A\mbox{ so that }
c_j(C_j-A)\to B\}.$$

Note that $L_A$ is smooth at $D,$ if $D\in\mathcal C_n.$
Then $C_A=\ker\sigma'_D$ and as $\dim \ker\sigma'_D=n^2-n$
(see Proposition \ref{cyc.equiv}), $C_A$ is an analytic set and
$\dim C_A=\dim L_A =n^2-n$ by \cite[Corollary, p. 83]{Chi}.
If $A\not\in\mathcal C_n$, then $\dim \ker\sigma'_A>n^2-n$ (see
Proposition \ref{cyc.equiv}) and so $C_A\subsetneq\ker\sigma'_A.$ Thus
$$C_A=\ker\sigma'_A\Leftrightarrow A\in\mathcal C_n.$$

The next theorem characterizes $C_A$ as the set of the "generalized"\
tangent vectors at $A$ to an entire curve in $\Omega_n$
passing through $A$ (in particular, this curve is contained in $L_A$).

\begin{theorem}\label{zeko.th1} Let $A\in\Omega_n$ and $B\in\mathcal M_n.$
Then there exists $m\in\Bbb N$ ($m\le n!$) and $\varphi\in\mathcal
O(\Bbb C,\Omega_n)$ so that $\varphi(0)=A,\
\varphi'(0)=\dots=\varphi^{(m-1)}(0)=0,\ \varphi^{(m)}(0)=B$ only
if $B \in C_A$.
\end{theorem}

We are not including the proof of this theorem, due to its length
and the use of results about analytic sets that are beyond the
scope of the dissertation. In can be found in the paper
\cite{Nik-Tho1} of the author and P.~J.~Thomas.

Theorem \ref{zeko.th1} shows that $C_A$ is contained in the set
of zeroes of the singular Kobayashi metric
$\kappa^s_{\Omega_n}(A;\cdot).$ Recall that (see \cite{Yu3})
$$\begin{aligned}\kappa^s_{\Omega_n}(A;B)=\inf\{|\alpha|:\exists
m\in\Bbb N,\varphi\in&\mathcal O(\Bbb
D,\Omega_n):\\&\hbox{ord}_0(\varphi-z)\ge m,\
\alpha\varphi^{(m)}(0)=m!X\}.
\end{aligned}$$

Now we define another cone $C'_A \subset \mathcal M_n,$
$A\in\Omega_n.$

For a function $g,$ holomorphic near $A$, and for $X$ in a neighborhood of $A$,
put $g(X)-g(A)=g_A^\ast(X-A)+\cdots$, where $g_A^\ast$ is the homogeneous polynomial
of lowest nonzero degree in the expansion of $g$ near $A.$ Put
$$C'_A=\{B\in\mathcal M_n:f_A^\ast(B)=0\mbox{ for each }f\in\O(\Om_n,\Bbb D)\}.$$

Note that
$$C_A\subset C'_A\subset\ker\sigma'_A;$$ the first inclusion
is proven e.g. in \cite[p. 86]{Chi}), and the second one
follows from the facts that each $f\in\O(\Om_n,\Bbb D)$ is constant on $L_A$
(by the Liouville theorem) and that
$$\ker\sigma'_A=\{(\sigma_j)_A^\ast=0 \mbox{ for all }j
\mbox{ such that } \deg (\sigma_j)_A^\ast=1 \}.$$

Also, each of these three sets $S$ is invariant under automorphisms of $\Omega_n.$

The cone $C'_A$ coincides with $C_A$ for $n=2$ and $n=3$ (for the
last fact see Proposition \ref{zeko.pr3} below and the remarks
preceding it). We do not know whether this is true for each $n.$

In the most trivial case of a non-cyclic matrix, namely a scalar one,
$C'_A=C_A$ is the set of zero-spectrum matrices, while $\ker\sigma'_A$
is the set of zero-trace matrices.

Note that $\kappa^s_{\Omega_n}\ge\gamma^s_{\Omega_n},$
where $\gamma^s_{\Omega_n}=\sup_{m\in\NN}\gamma^{(m)}_{\Omega_n}$
is the singular Carath\'eodory metric of $\Omega_n$ (see Section
\ref{est} for the definition of $\gamma^{(m)}$).

Theorem \ref{zeko.th1} implies that
$$B\in C_A\Rightarrow\kappa^s_{\Omega_n}(A;B)=0\Rightarrow
\gamma^s_{\Omega_n}(A;B)=0\Leftrightarrow B\in C'_A$$
(the last equivalence is trivial). In particular,
$$\kappa_{\Omega_n}(A;B)=0\Rightarrow B\in C'_A.$$

\begin{proposition}\label{zeko.nec}
If $A\in\Omega_n \setminus \mathcal C_n,$ then $C'_A
\neq\ker\sigma'_A$.
\end{proposition}

\beginproof As $A\in\Omega_n\setminus\mathcal C_n,$ at least two of
the eigenvalues of $A$ are equal, for example to $\lambda.$
Applying $\Phi_\lambda$ (see (\ref{fi})) and $\Psi_P$ (see
(\ref{gi})) we can assume that $\lambda=0$ and that $A$ is in a Jordan form.
In particular,
$$
A = \left( \begin{array}{cc} A_0 & 0 \\ 0 & A_1 \end{array} \right),
$$
where $A_0 \in \mathcal M_m$, $2\le m \le n$, $\spe(A_0)=\{0\}$,
$A_1 \in \mathcal M_{n-m}$, $0 \notin \spe(A_1)$.

Later, there exists a set $J \subsetneq \{2, \dots, m\}$,
possibly empty, such that $a_{j-1,j}=1$ for $j \in J$, and all other
elements $a_{ij}$ are equal to $0$ for $1\le i, j \le m$.
Put $0\le r = \# J = \mbox{rank} A_0 \le m-2.$ Let
$$B = \left(
\begin{array}{cc} B_0 & 0 \\ 0 & 0 \end{array} \right) \in
\mathcal M_n,$$ where $B_0 = (b_{ij})_{1\le i, j \le m}$ so that
$b_{j-1,j}=-1$ for $j \in \{2, \dots, m\} \setminus J$,
$b_{m1}=1$, and $b_{ij}=0$ otherwise.

As $\sigma_m/\binom{n}{m}\in\O(\Om_n,\Bbb D),$ it suffices to prove the following

\begin{lemma}
\label{B} $(\sigma_m)_A^\ast(B)=1,$ but $\sigma'_A(B)=0$.
\end{lemma}

\beginproof
First let us calculate $\sigma_j(A_0+hB_0),$ $1\le j \le m,$
$h\in\Bbb C.$ By developing along the first column, we get
$$
\det(tI-(A_0+h B_0))=t^m+(-1)^{m-1}h^{m-r}.
$$
Comparing the coefficients on both sides leads to
\begin{equation}\label{h}
\sigma_j(A_0+hB_0)=\left\{\begin{array}{ll}
0,&1\le j\le m-1\\
h^{m-r},&j=m\end{array}.\right.
\end{equation}

Now we need a general formula for $\sigma_j$. For a given matrix
$M=(m_{ij})_{1\le i, j \le n}$ and a set $E\subset \{1, \dots,
n\}$ denote by $\delta_E(M)$ the determinant of the matrix
$(m_{ij})_{i,j\in E} \in \mathcal M_{\#E}$. For convenience put
$\delta_\varnothing(M)=\sigma_0(M):=1$. Then
\begin{equation}
\label{sigmexp} \sigma_j (M) = \sum_{E\subset \{1, \dots, n\},
\#E=j} \delta_E(M).
\end{equation}
The block structure of our matrices implies that
$$
\delta_E(A+hB)= \delta_{E\cap \{1, \dots, m\}} (A_0 + h B_0)
\delta_{E\cap \{m+1, \dots, n\}} (A_1).
$$
So
\begin{multline*}
\sigma_j (A+hB) =\sum_{\max(0,j-n+m)\le k \le \min(m,j)} \left(
\sum_{E'\subset \{1, \dots, m\}, \#E'=k} \delta_{E'} (A_0 + h B_0)
\right) \times\\
\times \left( \sum_{E''\subset \{m+1, \dots, n\}, \#E''=j-k} \delta_{E''} (A_1)
\right)=
\\\sum_{\max(0,j-n+m)\le k \le \min(m,j)}\sigma_k(A_0+hB_0)\sigma_{j-k}(A_1).
\end{multline*}
By $(\ref{h})$ we get $\sigma_j(A+hB)= S_1 + S_2$, where
$$S_1=\left\{\begin{array}{ll}
\sigma_j(A_1),&j\le n-m\\
0,&\hbox{otherwise}\end{array},\right.\quad
S_2=\left\{\begin{array}{ll}
h^{m-r} \sigma_{j-m}(A_1),&j\ge m\\
0,&\hbox{otherwise}\end{array}.\right.
$$
In particular,
$$\sigma_j(A)=\left\{\begin{array}{ll}
\sigma_j(A_1),&j\le n-m\\
0,&\hbox{otherwise}\end{array}.\right.$$ Then
$$
\sigma_j (A+hB) - \sigma_j (A) = \left\{\begin{array}{ll}
h^{m-r} \sigma_{j-m}(A_1),&j\ge m\\
0,&\hbox{otherwise}\end{array}\right.
$$
As $m-r\ge 2,$ we get $\sigma'_A(B)=0$, but
$(\sigma_m)_A^\ast(B)=1.$\qed
\smallskip

The main corollary from Proposition \ref{zeko.nec} and the
implication preceding it form the fact that that SCFP does not
depend continuously on the data (so it cannot be reduced to a
similar problem on the symmetrized polydisc).

\begin{corollary}\label{zeko.disc}
If $A\in\Omega_n\setminus \mathcal C_n$ and $B\in\ker\sigma'_A\setminus{C'_A},$
then
$$\kappa_{\Omega_n}(A;B)>0=\lim_{\mathcal C_n\ni A'\to
A}\kappa_{\Omega_n}(A';B).$$
\end{corollary}

When $A$ is scalar matrix, we know more (cf. Proposition \ref{np.pr4}):

\begin{proposition}\label{e} For $B\in\mathcal M_n$ and $t\in\Bbb D$
the following are equivalent:

(i) the eigenvalues of $B$ ar equal;

(ii) the function $\kappa_{\Omega_n}$ is continuous at the point
$(tI_n;B);$

(iii) the function $\kappa_{\Omega_n}(\cdot;B)$ is continuous at
the point $tI_n.$
\end{proposition}

\beginproof The implication (ii)$\Rightarrow$(iii) is trivial.
For the rest of the proof we can assume that $t=0$
(applying $\Phi_t$).

We will now prove that (i)$\Rightarrow$(ii). Let the eigenvalues of
$B$ be equal to $0.$ If $A_j\to 0$ and $B_j\to B,$ then
$$\kappa_{\Omega_n}(A_j;B_j)\ge\kappa_{\Bbb
G_n}(\sigma(A_j);\sigma'_{A_j}(B_j))\to\kappa_{\Bbb
G_n}(0;\sigma'_0(B))=$$
$$\kappa_{\Bbb G_n}(0;(\tr B)e_1)=|b|=\kappa_{\Omega_n}(0;B)$$
(the last two equalities follow from Theorem \ref{i.th1}
(i) and (\ref{np3}), respectively).

Thus the function $\kappa_{\Omega_n}$ is lower semicontinuous
at the point $(0;B).$ As it is (always) upper semicontinuous,
it is continuous at this point.

It remains to prove that (iii)$\Rightarrow$(i). As $\mathcal
C_n$ is a dense subset of $\Omega_n,$ we can find a sequence
$\mathcal C_n\supset(A_j)\to 0.$ Then (\ref{np3}),
(\ref{np1}) and Theorem \ref{i.th1} (i) imply that
$$r(B)=\kappa_{\Omega_n}(0;B)\leftarrow
\kappa_{\Omega_n}(A_j;B)=$$ $$\kappa_{\Bbb
G_n}(\sigma_(A_j);\sigma'_{A_j}(B_j))\to\kappa_{\Bbb
G_n}(0;\sigma'_0(B))=|\tr B|/n.$$ So $r(B)=|\tr B|/n,$
i.e. the eigenvalues of $B$ are equal.\qed
\smallskip

Now let us formulate the following hypothesis for the zeroes of $\kappa_{\Omega_n}.$
\smallskip

\noindent{\bf Hypothesis.} $\kappa_{\Omega_n}(A;B)=0$ if and only
if there exists a $\varphi\in \mathcal O(\mathbb C,\Omega_n)$ so
that $\varphi(0)=A$ and $\varphi'(0)=B$. In particular, if
$\kappa_{\Omega_n}(A;B)=0,$ then $B\in C_A.$
\smallskip

Note that there are matrices $B\in C_A$ such that $\kappa_{\Omega_n}(A;B)\neq 0$
(see Proposition \ref{zeko.pr3} (ii)
and Corollary \ref{zeko.cor4}).

In some cases the above hypothesis can be checked.

The remarks at the beginning of this section imply that this hypothesis is
true for cyclic matrices.

Also, as the zeroes of $\kappa_{\Omega_n}(0;\cdot)$ are exactly the zero-spectrum
matrices and this set of matrices is a union of complex lines through the origin,
the hypothesis is true for scalar matrices.

As the non-cyclic matrices $A$ in $\Omega_2$ are only the scalar
ones, we can choose $m=1$ in Theorem \ref{zeko.th1} for $n=2$;
then $C_A$ coincides with the zeroes of
$\kappa_{\Omega_2}(A;\cdot),$ as well as with the set of matrices
$B=\varphi'(0)$ for some entire curve $\varphi$ in $\Omega_2$. (On
the other hand, $\ker\sigma'_A=\{B\in\mathcal M_2:\tr B=0\}$.) So
we have a complete description of the set of zeroes of
$\kappa_{\Omega_2}$ and the above hypothesis is true for $n=2.$

Now let us consider the set of zeroes of
$\kappa_{\Omega_3}(A;\cdot),$ when $A$ is a non-cyclic and
non-scalar matrix (we will confirm the hypothesis for $n=3$, too).
The use of the automorphisms $\Phi_\lambda$ and $\Psi_P$ of
$\Om_3$ reduces the problem to the following two cases:
$$A=A_t:=\left(\begin{array}{ccc}
0&0&0\\0&0&0\\0&0&t\\
\end{array}\right),\ t\in\Bbb D_\ast,\quad
A=\tilde A:=\left(\begin{array}{ccc}
0&0&0\\0&0&1\\0&0&0\\
\end{array}\right).$$
It is easily seen that
$$C'_{A_t}\subset C''_{A_t}:=\{B\in\mathcal M_3:
\sigma^\ast_{A_t}(B)\}=$$
$$\{B\in\mathcal M_3:b_{33}=b_{11}+b_{22}=b_{11}^2+b_{12}b_{21}=0\},$$
$$C'_{\tilde A}\subset C''_{\tilde A}:=\{B\in\mathcal M_3:
\sigma^\ast_{\tilde A}(B)\}=$$
$$\{B\in\mathcal M_3:b_{11}+b_{22}+b_{33}=b_{32}=b_{12}b_{31}=0\}$$
(for example, to check the second equality, we observe that if
$B_\eps=\tilde A+\eps B+o(\eps),$ then $\tr B_\eps=\eps \tr
B+o(\eps),$ $\sigma_2(B_\eps)=-\eps b_{32}+o(\eps)$ and $\det
B_\eps=\eps^2(b_{12}b_{31}-b_{11}b_{32})+o(\eps^2)$). As the
tangent cones are closed, the next proposition shows in particular
that $C_{A_\lambda}=C'_{A_\lambda}=C''_{A_\lambda}$ and $C_{\tilde
A}=C'_{\tilde A}=C''_{\tilde A}.$

\begin{proposition}\label{zeko.pr3} (i) If $B\in C''_{A_t}$ ($t\neq 0$), then
there exists a $\varphi\in\mathcal O(\Bbb C,\Omega_3)$ such that
$\varphi(0)=A_t$ and $\varphi'(0)=B.$

(ii) Let $B\in C''_{\tilde A}.$ Then there exists a
$\varphi\in\mathcal O(\Bbb C,\Omega_n)$ so that
$\varphi(0)=\tilde A$ and $\varphi'(0)=B$ only if $b_{11}=0$ and
$b_{12}\neq b_{31}.$ Otherwise $\kappa_{\Omega_3}(\tilde
A;B)=1.$
\end{proposition}

As $\kappa_{\Omega_3}(A;B)>0$ for $B\not\in C'_A,$ this
proposition and the remarks preceding it give a complete
description of the set of zeroes of $\kappa_{\Omega_3},$ thereby
confirming the hypothesis for $n=3.$
\smallskip

\beginproof (i) Let first $B\in
C'_{A_t}.$ We express $B$ in the form $B=X+[Y,A_t],$ where $X$ is
such that $\psi(\zeta)=A_t+\zeta X\in L_{A_t}$ for each
$\zeta\in\Bbb C.$ Then $\varphi(\zeta)=e^{\zeta
Y}\psi(\zeta)e^{-\zeta Y}$ has the required properties.

It is easily calculated that $\psi(\CC)\subset L_{A_t}$ exactly when
$\spe(X)=0$ and $x_{11}+x_{22}=x_{11}^2+x_{12}x_{21}=0.$ On the other hand,
$$[Y,A_t]=t\left(\begin{array}{ccc}
0&0&y_{13}\\0&0&y_{23}\\-y_{31}&-y_{32}&0\\
\end{array}\right).$$

So we can choose
$$X=\left(\begin{array}{ccc}
b_{11}&b_{12}&0\\b_{21}&b_{22}&0\\0&0&0\end{array}\right),\quad
Y=t^{-1}\left(\begin{array}{ccc}
0&0&b_{13}\\0&0&b_{23}\\-b_{31}&-b_{32}&0\end{array}\right).$$

(ii) Let first $B\in\CC'_{\tilde A}.$ If $b_{11}=0$ or
$b_{12}\neq b_{31},$ it suffices to find (as above) $X$ and $Y$ so that
$B=X+[Y,\tilde A]$ and $\tilde A+\zeta X\in L_{\tilde
A}$ for each $\zeta\in\Bbb C.$ The last
condition means that the eigenvalues of $X$ are zeroes and $x_{32}=x_{12}x_{31}=0.$
On the other hand,
$$[Y,\tilde A]=\left(\begin{array}{ccc}
0&0&y_{12}\\-y_{31}&-y_{33}&y_{22}-y_{33}\\0&0&y_{32}\\
\end{array}\right).$$

Let us suppose that $b_{31}=0$ (when $b_{12}=0$ the calculations are analogical).
We have to choose $X$ of the form
$$X=\left(\begin{array}{ccc} b_{11}&b_{12}&b_{13}-y_{12}
\\b_{21}+y_{31}&b_{22}+y_{32}&b_{23}-y_{22}+y_{33}\\0&0&-b_{11}-b_{22}-y_{32}\\
\end{array}\right)$$
so that $\det X=0$ and $\sigma_2(X)=0,$ i.e. $DT=0$ and $D=T^2,$
where
$$
D= \left| \begin{array}{cc} b_{11}&b_{12}
\\b_{21}+y_{31}&b_{22}+y_{32} \end{array} \right|, \quad
T= b_{11}+b_{22}+y_{32}.
$$
These two conditions are true only if $$y_{32}=-b_{11}-b_{22},\quad y_{31}=\left\{\begin{array}{ll}
-b_{21},&b_{11}=0\\
-b_{21}-\frac{b_{11}^2}{b_{12}},&b_{12}\neq 0
\end{array}.\right.$$

It remains to show that if $b_{11}\neq 0$ and $b_{12}=b_{31}=0,$
then $\kappa_{\Omega_3}(\tilde A;B)=1.$ We can assume that
$b_{11}=1.$ Put $\tilde B=\diag(1,e^{2\pi i/3},e^{4\pi i/3}.$ As
above, we can choose $\tilde B$ and $Y$ so that $B=\tilde
B+[Y,A_t].$ Let $\alpha>0$ and $\varphi\in\mathcal O(\alpha\mathbb
D,\Omega_3)$ be such that $\varphi(0)=A_t$ and $\varphi'(0)=B.$
Putting $\tilde\varphi(\zeta)=e^{-\zeta Y}\varphi(\zeta)e^{\zeta
Y},$ we have $\tilde\varphi\in\mathcal O(\alpha\D,\Omega_3),$
$\tilde\varphi(0)=\tilde A$ and $\tilde \varphi'(0)=\tilde B.$ So
$\kappa_{\Omega_3}(\tilde A;B)\ge\kappa_{\Omega_3}(\tilde A;\tilde
B).$ The converse inequality follows similarly. It remains to
apply Proposition \ref{exam} from the next section.

\begin{corollary}\label{zeko.cor4} For each $n\ge 3 $ there exist $A\in\Omega_n$
and $B\in C_A$ so that $\kappa_{\Omega_n}(A;B)>0.$
\end{corollary}

\beginproof Put
$$\tilde A=\left(\begin{array}{ccc}
0&0&0\\0&0&1\\0&0&0\\
\end{array}\right),\
\tilde B_\eps=\left(\begin{array}{ccc}
1&\eps&0\\0&-1&0\\0&0&0\\
\end{array}\right),$$
$$A=\left(\begin{array}{cc}\tilde A&O\\O&O
\end{array}\right),\
B_\eps=\left(\begin{array}{cc}\tilde B_\eps&O\\O&O
\end{array}\right).$$
As in the proof of Proposition \ref{zeko.pr3} (ii) it follows that

$\bullet$ $\kappa_{\Omega_n}(A;B_0)>0;$

$\bullet$ for $\eps\neq 0$ there exists $\varphi_\eps\in\mathcal
O(\Bbb C,\Omega_n)$ such that $\varphi_\eps(0)=A$ and $\varphi_\eps'(0)=B_\eps.$

Then $B_\eps\in C_A,$ $\eps\neq 0,$ so $B_0\in C_A.$\qed

\setcounter{equation}{0}
\section{The Kobayashi metric vs. the Lempert function}\label{k-l}

As an application of part of the above considerations, in this
section we will provide an example showing that, in general, the
Kobayashi pseudometric of a pseudoconvex domain is not equal to
the weak "derivative"\ of the Lempert function. The pseudoconvex
domain will be the spectral ball $\Om_3\subset\C^9$ (that is also
a balanced non-taut unbounded domain).

Recall that the Kobayashi metric of a taut domain $D\subset\C^n$
coincides with the "derivative"\ of the Lempert function (see
Section \ref{lem-der}):
$$\kappa_D(z;X)=\lim_{t\rightarrow 0,z'\to z,X'\to X}\frac{l_D(z',z'+tX')}{|t|}.$$
On the other hand, Proposition \ref{der.pr2} states that
\begin{equation}\label{k-l1}
\kappa_D(z;X)\ge\CD l_D(z;X):=\limsup_{t\rightarrow 0,z'\to
z,X'\to X}\frac{l_D(z',z'+tX')}{|t|}
\end{equation}
for an arbitrary domain $D\subset\C^n.$

The aim of this section is to show that the inequality
\begin{equation}\label{k-l2}
\kappa_D(z;X)\ge\wdtl\CD l_D(z;X):=\limsup_{t\rightarrow
0}\frac{l_D(z,z+tX)}{|t|}.
\end{equation}
is strict in the general case (of a pseudoconvex domain).
\smallskip

Put
$$A=\left(\begin{array}{ccc}0&0&0\\0&0&1\\0&0&0\end{array}\right)\text{\ \ and\ \ }
B_t=\left(\begin{array}{ccc}1&0&0\\0&\omega&0\\0&3t&\omega^2\end{array}\right).$$
where $\omega=e^{2\pi i/3}$. Let $B=B_0.$

\begin{proposition}\label{Prop-Ex} The following inequality holds:
$$\kappa_{\Om_3}(A;B)>0=\wdtl\CD l_{\Om_3}(A;B).$$

Moreover, if $t_j\to 0$ and $C_j\to B$ ($C_j=(c_{k,l}^j)$) so that
$\liminf_{j\to\infty}|c_{3,2}^j/t_j-3|>0,$ then
$$\lim_{j\to\infty}\frac{l_{\Om_3}(A,A+t_jC_j)}{|t_j|}=0.$$
\end{proposition}

\noindent{\bf Remark.} As $\kappa_D$ and $l_D$ have the product
property (see Section \ref{desc}), in general the inequality
(\ref{k-l2}) is strict for pseudoconvex domains in $\C^n$ for
$n\ge 9$ (for example for $\Om_3\times\D^k$). In fact the proof
below shows that $\wdtl\CD l_{\wdtl\Om_3}(A;B)=0,$ where
$\wdtl\Om_3$ is the set of zero-trace matrices in $\Om_3$.
Consequently the inequality in (\ref{k-l2}) is strict for the
pseudoconvex domain $\wdtl\Om_3\subset\C^8.$
\smallskip

\noindent{\bf Question.} It would be interesting to find an
example of a lower dimension, as well as to see whether in general
the inequality (\ref{k-l1}) is strict (the last question was posed
at the end of Section \ref{lem-der}).
\smallskip

Recall that there exist matrices $\wdtl B\to B$ so that
$\kappa_{\Om_3}(A;\wdtl B)=0$ (see Proposition \ref{zeko.pr3} (ii));
in particular, the function $\kappa_{\Om_3}(A;\cdot)$ is not continuous
at $B.$

Also note that the condition
$\liminf_{j\to\infty}|c_{3,2}^j/t_j-3|>0$ in Proposition
\ref{Prop-Ex} is essential, as seen by the following

\begin{proposition}\label{exam} $\kappa_{\wdtl\Om_3}(A;B)=
\lim_{t\rightarrow 0}\frac{l_{\Om_3}(A,A+tB_t)}{|t|}=1.$

In particular,
$$1=\kappa_{\wdtl\Om_3}(A;B)=\kappa_{\Om_3}(A;B)=\CD
l_{\wdtl\Om_3}(A;B)=\CD l_{\Om_3}(A;B).$$
\end{proposition}

For the proof of Proposition \ref{Prop-Ex} we will use the following
special case of \cite[Proposition 4.1, Corollary
4.3]{Tho-Tra2} (see also \cite{NPT2} for more general facts).

\begin{lemma}\label{lift-1} Let $M\in\Om_3$ is a cyclic matrix and
$\phi\in\O(\D,\G_3)$ is such a mapping that $\phi(0)=0$ and
$\phi(\alpha)=\sigma(M)$ ($\alpha\in\D$). Then there exists e
$\psi\in\O(\D,\Om_3)$ such that $\psi(0)=A,$ $\psi(\alpha)=M$ and
$\phi=\sigma\circ\psi$ exaclty when $\phi_3'(0)=0$.

In particular,
$$
l_{\Om_3}(A,M)=
\inf\{|\alpha|:\exists\phi\in\O(\D,\G_3):\phi(0)=0,\phi(\alpha)=\sigma(M),\phi'_3(0)=0\}$$
and (as $\G_3$ is a taut domain) there exists an extremal disc for
$l_{\Om_3}(A,M)$.
\end{lemma}

\beginproof If such a $\psi$ exists, then one directly calculates
$\phi_3'(0)=(\sigma_3\circ\psi)'(0)=0.$

Conversely, let $\phi_3'(0)=0$. Put
$$\tilde\psi(\zeta):=\left(\begin{array}{ccc}0&\zeta&0\\0&0&1\\
\phi_3(\zeta)/\zeta&-\phi_2(\zeta)&\phi_1(\zeta)\end{array}\right),\quad\zeta\in\D.$$
Then $\tilde\psi(0)=A$ and $\phi=\sigma\circ\tilde\psi.$
Furthermore, $e_3=(0,0,1)$ is a cyclic vector for
$\tilde\psi(\zeta)$ when $\zeta\neq 0.$ So $\tilde\psi(\alpha)$ is
a cyclic matrix with the same spectrum as the cyclic matrix $M$
and consequently they are similar (to their adjoint matrix) by
Proposition \ref{cyc.equiv}. Then we can express $M$ in the form
$M=e^S\tilde\psi(\alpha)e^{-S}$ for some $S\in\M_3.$ It remains to
put $\psi(\zeta)=e^{\zeta S/\alpha}\tilde\psi(\zeta)e^{-\zeta
S/\alpha}.$\qed
\smallskip

\noindent{\it Proof of Proposition \ref{Prop-Ex}.} By
Proposition \ref{exam} we only need to check that
$$
\lim_{j\to\infty}\frac{l_{\Om_3}(A,A+t_jC_j)}{|t_j|}=0
$$
in the conditions for $c_{3,2}^j$.

Suppose the contrary. Then we may assume that
$$
\frac{l_{\Om_3}(A,A+t_jC_j)}{|t_j|}\to a>0.
$$

{\it Step 1.} Suppose that there exists a subsequence
(for brevity we still use the index $j$) such that all matrices
$A+t_jC_j$ are cyclic and belong to $\Om_3$. By some calculations we get
$$
\sigma(A+t_jC_j)=(t_jf_1(C_j),t_jf_2(C_j),t_j^2f_3(C_j))=:(a_j,b_j,c_j),
$$
where $f_1(C_j)\to 0$, $f_2(C_j)\to 0$ and $f_3(C_j)\to 0$.

Put $$\phi_j(\zeta)=(\zeta a_j/r_j, \zeta b_j/r_j,\zeta^2
c_j/r_j^2), \quad\zeta\in\D,
$$
where $r_j=\max\{3|a_j|,3|b_j|,\sqrt{3|c_j|}\}$. Then
$\phi_j\in\O(\D,\G_3),$ with $\phi_j(0)=0$, $\phi_{j,3}'(0)=0$ and
$\phi_j(r_j)=\sigma(A+t_jC_j)$. Lemma \ref{lift-1} implies that
$$
l_{\Om_3}(A,A+t_jC_j)/|t_j|\leq r_j/|t_j|\to 0,
$$
a contradiction.

{\it Step 2.} Suppose that all matrices $A+t_jC_j$ are non-cyclic.
Then their minimal polynomials have degrees less than
3 (see Proposition \ref{cyc.equiv}). Consequently these degrees are
equal to 2 for all sufficiently large $j.$ Hence
$$(A+t_jC_j)^2+x_j(A+t_jC_j)+y_jI_3=0,$$
where $x_j, y_j\in\C.$ We get $9$ equations (for the components);
denote them by $E_{k,\ell}^j$, where $k$ and $\ell$
are the indices of the row and the column, respectively.
The equation $E_{2,3}^j$
gives $x_j/t_j\to 1$. Using this in $E_{1,1}^j,$ we get $y_j/t_j^2\to -2$.
Finally, $E_{2,2}^j$ implies
$c_{3,2}^j/t_j\to 2-\omega-\omega^2=3$ -- a contradiction.

This completes the proof.\qed
\smallskip

\noindent{\it Proof of Proposition \ref{exam}.} As
$A+\zeta B\in\wdtl\Om_3$ for each $\zeta\in\D,$ we get
$\kappa_{\wdtl\Om_3}(A;B)\le 1.$

It remains to show that
$$\liminf_{t\rightarrow 0}\frac{l_{\Om_3}(A,A+tB_t)}{|t|}\ge 1.$$

Note that $A+tB_t$ is similar to the matrix
$D_t=\mbox{diag}(t,t-2t)$ and consequently $l_{\Om_3}(A,A+tB_t)
=l_{\Om_3}(A,D_t)$ (we already applied this argument several times).

Suppose that $t_j\to 0$ so that $l_{\Om_3}(A,D_{t_j})/|t_j|\to
c<1$.

Let $\psi_j\in\O(\D,\Om_3)$ be a disc such that $\psi_j(0)=A$,
$\psi(\alpha_j)=D_{t_j}$ and $|\alpha_j|/|t_j|\to c$. Put
$\phi_j=\sigma\circ\psi_j.$ Direct calculations lead to $\phi'_{j,3}(0)=0$ and $$
\phi_{j,3}'(\alpha_j)-t_j\phi_{j,2}'(\alpha_j)+t_j^2\phi_{j,1}'(\alpha_j)=0.
$$
Expressing $\phi_j$ in the form
$$
\phi_j(\zeta)=(\zeta\theta_{j,1}(\zeta),\zeta\theta_{j,2}(\zeta),
\zeta^2\theta_{j,3}(\zeta)),
$$
the last equality becomes
\begin{equation}\label{prel}
t_j^3=\alpha_j^2(\alpha_j\theta_{j,3}'(\alpha_j)-t_j\theta_{j,2}'
(\alpha_j)+t_j^2\theta_{j,1}'(\alpha_j))
\end{equation}
(we use that $\theta_{j,1}(\alpha_j)=0,$
$\theta_{j,2}(\alpha_j)=-3t_j^2/\alpha_j$ and
$\theta_{j,3}(\alpha_j)=-2t_j^3/\alpha_j^2$). As $\G_3$ is a taut
domain, by passing to subsequences we can assume that
$\phi_j\to\phi=(\zeta\rho_1,\zeta^2\rho_2,\zeta^3\rho_3)\in\O(\D,\G_3)$
and $\rho_1(0)=0.$ Then the equation (\ref{prel}) shows that if
$k=1/c,$ then
$$\rho_3(0)=k^3+k\rho_2(0).$$

Proposition \ref{b.pr2} (see also \cite[Proposition 16]{Edi-Zwo})
implies that
$$h_{\G_3}(z)=\max\{|\lambda|:\lambda^3-z_1\lambda^2+z_2\lambda-z_3=0\}$$
is a (logarithmically) plurisubharmonic function and
$\G_3=\{z\in\C^3:h_{\G_3}(z)<1\}$ ($h_{G_3}$ is the function of
Minkowski of the $(1,2,3)$-balanced domain $\G_3$). As
$$|\zeta| h_{\G_3}(\rho_1(\zeta),\rho_2(\zeta),\rho_3(\zeta))=
h_{\G_3}(\phi(\zeta))<1,\quad \zeta\in\D,$$ the maximum principle
for plurisubharmonic functions implies
$h_{\G_3}(\rho_1,\rho_2,\rho_3)\leq 1$ on $\D$. In particular,
$h_{\G_3}(\rho_1(0),\rho_2(0),\rho_3(0))\leq 1$. Consequently all the
three zeroes of the polynomial
$P(\lambda)=\lambda^3-\rho_1(0)\lambda^2+\rho_2(0)\lambda-\rho_3(0)$,
lie in $\overline\D$. On the other hand,
$P(\lambda)=(\lambda-k)(\lambda^2+k\lambda+k^2+\rho_2(0)),$ with
$k>1$ -- a contradiction.\qed

\chapter{Estimates and boundary behavior of invariant metrics on
$\C$-convex domains} \label{chap.bound}

\section{Estimates for the Carath\'eodory and Kobayashi metrics}\label{car-kob}

The aim of this section is obtaining estimates for the Kobayashi
and Carath\'eodory metrics on $\C$-convex domains in terms of the
distance to the boundary of the corresponding direction. These
results generalize similar statements for bounded smooth
$\C$-convex domains of finite type, whose original proofs are
quite hard (see \cite{Lie}).

For a point $z$ from a domain $D\subset\C^n$ and a vector
$X\in(\C^n)_\ast$, we denote by $d_D(z;X)$ the distance from $z$
to $\partial D$ in the direction of $X$, i.e.
$$d_D(z;X)=\sup\{r>0:\Delta_X(z,r)\subset D\},$$
where $$\Delta_X(z,r)=\{z+\lambda X:|\lambda|<r\}.$$ Clearly
$$\dist(z,\partial D)=:d_D(z)=\inf_{||X||=1}d_D(z;X).$$
If $d_D(z;X)=\infty,$ i.e. $D$ contains the line through $z$ in the direction of
direction $X,$ then
$$\gamma_D(z;X)=\kappa_D(z;X)=0.$$

First recall the following result for convex domains.

\begin{proposition}\label{bound0} \cite{BP} Let $D\subset\C^n$ be a
convex domain. If $d_D(z;X)<\infty,$ then
$$1/2\le\gamma_D(z;X)d_D(z;X)=\kappa_D(z;X)d_D(z;X)\le
1.$$
\end{proposition}

\beginproof The upper estimate holds for each domain $D,$ as
$\D_X(z,d_D(z,X))\subset D.$ For the lower estimate consider an
(open) supporting half-space $\Pi$ of $D$ for a boundary point of
the type $z+\lambda X.$ Then
$$\gamma_D(z;X)\ge
\gamma_{\Pi}(z;X)=\frac{||X||}{d_{\Pi}(z;X)}=\frac{||X||}{d_D(z;X)}.$$
It remains to note that the equality in the proposition follows
from the Lempert theorem (see e.g. \cite{Lem1,Lem2}).\qed
\smallskip

The constants $1/2$ and $1$ cannot be improved, as seen from the
examples of a half-space and a ball.

Now we will establish a similar result for $\C$-convex domains.

\begin{proposition}\label{bound1} Let $D\subset\C^n$ be a $\C$-convex domain.
If $d_D(z;X)<\infty,$ then
$$1/4\le\gamma_D(z;X)d_D(z;X)\le\kappa_D(z;X)d_D(z;X)\le
1.$$
\end{proposition}

The constant $1/4$ is the best possible in the plane, as seen in
the example with the image $D=\CC\setminus[1/4,\infty)$ of $\D$
for the K\"obe transformation $z\to z/(1+z)^2.$

\begin{corollary}\label{bound2} For each $\CC$-convex domain $D\subset\CC^n$,
we have $\kappa_D\le 4\gamma_D.$
\end{corollary}

This is another argument supporting the hypothesis that
$\kappa_D=\gamma_D$ for each $\CC$-convex domain $D\subset\CC^n$
(a weaker variant of \cite[\bf Problem 4']{Zna}; see the
Introduction).
\smallskip

\noindent{\it Proof of Proposition \ref{bound1}}. We can assume
that $||X||=1.$ Let $l$ be the complex line through $z$ with
direction $X$ and $a\in l\cap\partial D$ so that
$||z-a||=d_D(z;X).$ Consider a complex hyperplane $H$ through $a$
not intersecting $D$ and denote by $G$ the projection of $D$ onto
$l$ in the direction of $H.$ Note that $G$ is a simply connected
domain (see e.g. \cite[Theorem 2.3.6]{APS} or \cite[Theorem
2.3.6]{Hor2}), $a\in\partial G$ and
$d_D(z;X)=||z-a||=d_G(z):=\dist(z,\partial G).$ It remains to
apply the K\"obe $1/4$ theorem to get
$$\gamma_D(z;X)\ge\gamma_G(z;1)\ge\frac{1}{4 d_G(z)}.$$ Indeed, if
$f:\D\to G$ is a conformal mapping such that $f(0)=z,$ by the
K\"obe theorem $G$ contains the disc of center $z$ and radius
$|f'(0)|/4.$ So $|f'(0)|\le 4d_G(z)$ and then
$$1=\gamma_{\D}(0;1)=\gamma_G(f(0);f'(0))=|f'(0)|\gamma_G(z;1)\le
4d_G(z)\gamma_G(z;1)$$ and the result follows. \qed
\smallskip

Recall that if a $\Bbb C$-convex domain in $\C^n$ contains a
complex line, then it is linearly equivalent to the Cartesian
product of $\Bbb C$ and a $\Bbb C$-convex domain in $\Bbb C^{n-1}$
(see Section \ref{con}).

In view of this it is natural to ask about the boundary behavior
of the metrics in the directions, where there are (linear) discs
in the boundary in these directions.

More precisely, for a boundary point $a$ of a domain
$D\subset\CC^n$ we denote by $L_a$ the set of all vectors
$X\in\CC^n$ such that there exists a $\varepsilon>0$ so that
$\partial D\supset\Delta_X(a,\varepsilon).$ The following result
is an application of Proposition \ref{bound1}.

\begin{proposition}\label{bound3} Let $a$ be a boundary point of a
$\CC$-convex domain $D\subset\CC^n.$

(i) Then
$$\lim_{z\to a}\gamma_D(z;X)=\infty$$
locally uniformly on $X\not\in L_a.$

(ii) If $\partial D$ is $\mathcal C^1$-smooth at $a,$ then $L_a$
is a linear space. In addition, for each non-tangent cone
$\Lambda$ with vertex
$a$*\footnote{*$\Lambda=\{z\in\C^n:c||z-a||<|\mbox{pr}_a(z)|\},$
where $c\in(0,1)$ and $\mbox{pr}_a$ is the projection onto the
internal normal to $\partial D$ at $a.$} we have
$$\limsup_{\Lambda\ni z\to a}\kappa_D(z;X)<\infty$$ locally uniformly on $X\in L_a.$
\end{proposition}

The proof of this proposition, as well as of a part of the next
ones, will be based on the following geometrical property of the
weakly linearly convex domains (see also \cite{Zna-Zna}).

\begin{lemma}\label{bound11} Suppose that the weakly linearly convex
domain $G\subset\Bbb C^n$ contains the $n$ unit discs lying in the coordinate lines.
Then $G$ contains the convex hull of these discs $E=\{z\in\Bbb C^n:\sum_{j=1}^n|z_j|<1\}.$
\end{lemma}

\beginproof For each $\eps\in (0,1)$ there exists a $\delta>0$ so that
$$
X_\eps=\bigcup_{j=1}^n\Big(\delta\Bbb D\times\dots\times\delta\Bbb
D\times\underbrace{\eps\Bbb D}_{j-\mbox{\small th place}}\times
\delta\Bbb D\times\dots\times\delta\Bbb D\Big)\subset G.$$ Note that
$\widehat X_\eps\subset G,$ where $\widehat X_\eps$
is the least linearly convex set that contains $X_\eps$.
In addition,
$$
\widehat X_\eps=\{z\in\CC^n|\forall b\in\CC^n:<z,b>=1\ \exists
a\in X_\eps:<a,b>=1\}.
$$
(see e.g. \cite[p. 7]{APS} or \cite[Proposition 4.6.2]{Hor2}.
Then $\widehat X_\eps$ is a balanced domain and as it is linearly
convex, it is convex (see Proposition \ref{c.pr1}).
Consequently
$$E_\eps=\{z\in\Bbb C^n:\sum_{j=1}^n|z_j|<\eps\}\subset\widehat X_\eps\subset
G$$ and for $\eps\to 1$ we get the desired proposition.\qed
\smallskip

\noindent{\bf Remark.} The same arguments show that $G$
contains the convex hull of each of its balanced sub-domains. In particular,
the maximal balanced sub-domain of $G$ is convex
(see also \cite{Zna-Zna}).
\smallskip

\noindent{\it Proof of Proposition \ref{bound3}.} (i) Assuming the
contrary, we can find $r>0$ and sequences $D\supset(z_j)\to a,\
\CC^n\supset(X_j)\to X\not\in L_a$ such that
$4r\gamma_D(z_j;X_j)\le1.$ By Proposition \ref{bound1},
$d_D(z_j;X_j)\ge r.$ Then $\Delta_{X_j}(z_j,r)\subset
D_r=D\cap\Bbb B_n(a,2r)$ for each sufficiently large $j.$ Note
that $D_r$ is a (weakly) linearly convex open set. By Proposition
\ref{c.pr0}, it is taut. Therefore $\Delta_X(a,r)\subset\partial
D,$ a contradiction.

(ii) As $\partial D$ is $\mathcal C^1$-smooth, for each two linearly
independent vectors $X,Y\in L_a$ one can find a neighborhood
$U$ of $a$ and a number $\varepsilon>0$ so that
$\Delta_X(z,\varepsilon)\subset D$ and $\Delta_Y(z,\varepsilon)\subset D$
for $z\in D\cap U\cap\Lambda.$
By Lemma \ref{bound11}, $\Delta_{X+Y}(z,\varepsilon')\subset D$ for some $\varepsilon'>0.$
As in (i) we get $\Delta_{X+Y}(a,\varepsilon')\subset\partial D.$ Consequently
$L_a$ is linear space. Then, choosing a basis for $L_a$ and applying
Lemma \ref{bound11}, we find a neighborhood $U$ of $a$ and a number
$c>0$ so that $\Delta_X(z,c)\subset D$ for each $z\in
D\cap U\cap\Lambda$ and each unit vector $X\in L_a.$ Now
the required estimate follows from Proposition \ref{bound1}.
\smallskip

\noindent{\bf Remark.} The condition for smoothness is redundant, if $D$ is a
convex domain. Indeed, in this case it is clear that $L_a$ is a linear space. Also, if
$\Delta\subset\partial D,$ then for each point $b\in D$ and each
number $t\in(0,1]$ we have $tb+(1-t)\Delta\subset D.$ So we can replace
$\Lambda$ by an arbitrary cone with vertex $a$ having as a base an arbitrary compact set of $D.$

\setcounter{equation}{0}
\section{Types of boundary points}\label{type}

The aim of this section is finding estimates for the behavior of invariant
metrics of $\C$-convex domains near a boundary point
depending on the multitype of this point.

Let $a$ be a ($\mathcal C^\infty$-)smooth boundary point of a
domain $D\subset\CC^n.$ Denote by $m_a$ the (D'Angelo) type of
$a,$ i.e. the maximal order of tangency of $\partial D$ at $a$
with (nontrivial) analytic discs through $a$ (see e.g. the
Ph.~D.~thesis \cite{Nik0} of the author, we will refer to it
several times in this chapter):
$$m_a=\sup_{\gamma}\frac{\ord_a(r\circ\gamma)}{\ord_a\gamma},$$ where
$\gamma$ varies over all analytic discs through $a,$ while $r$ is a smooth
defining function of $D$ near $a$ (this definition
depends on $r$). By requiring $\gamma^{\ord_a\gamma}(a)=X,$ for a given
vector $X\in(\C^n)_\ast,$ we define the number $m_{a,X}.$

The point $a$ is said to be of finite type, if $m_a<\infty.$ A bounded domain
$D$ is said to be of finite type if all its boundary
points are of finite type.

Replacing the analytic discs by complex lines, we define the linear type
$l_a$ of $a.$ We can also define the number $l_{a,X}$ as the
order of tangency of $\partial D$ at $a$ to the line through $a$ in the direction of $X.$

Then $l_{a,X}\le m_{a,X}$ and $l_a\le m_a.$ Note that if
$l_{a,X}<\infty,$ then $X\not\in L_a.$

\begin{proposition}\label{type4} Let $a$ be a smooth boundary point of a
$\Bbb C$-convex domain $D\subset\CC^n$ and let
$X\in(\CC^n)_\ast$ so that $l_{a,X}<\infty.$ Denote by $n_a$
the internal normal to $\partial D$ at $a.$ Then there exists a
neighborhood $U_X$ of $a$ and a constant $c_X\ge 1$ so that
$$c_X^{-1}d_D(z)\le d_D(z;X)^{l_{a,X}}\le c_X d_D(z),\quad z\in D\cap U_X\cap
n_a.$$
\end{proposition}

\beginproof We can assume that $\Re z_1<0$ is the internal normal to
$\partial D$ at $a=0.$ Let $r(z)=\Re z_1+o(|z_1|)+\rho('z)$ be a smooth
defining function of $D$ near $0.$

For each sufficiently small $\delta>0$ we have
$\delta=d_D(\delta_n),$ where $\delta_n=(-\delta,'0).$ Put
$L_\delta(\zeta)=-\delta_n+\zeta X,$ $\zeta\in\CC^n.$

We consider two cases.

1. $l_{a,X}=1.$ This means that $X_1\neq 0.$ Then
$r(L_\delta(\zeta))=-\delta+\Re(\zeta X_1)+o(|\zeta|).$
Consequently $L_\delta(\zeta)\in D,$ if
$|\zeta|<\frac{\delta}{2|X_1|}$ and $\delta$ is sufficiently small.
This proves the left-hand side inequality.

The right hand-side inequality follows from the inequality
$r(L_\delta(2\delta/X_1))>0,$ which holds for each small
$\delta>0.$

2. $l_{a,X}\ge 2.$ This means that $X_1=0.$ Then
$r(L_\delta(\zeta))=-\delta+\rho(\zeta 'X ).$ As $\rho(\zeta
'X)\le c|\zeta|^l$ for some $c_1>0,$ we get
$L_\delta(\zeta)\in D,$ if $c_1|\zeta|^l<\delta.$
This proves the left-hand side inequality.

To prove the right hand-side inequality, we need to find a $c_2>0$
so that for each small $\delta>0$ there exists $\zeta$
such that $|\zeta|^l=c_2^{-1}\delta$ and $\rho(\zeta'X)\ge\delta.$ As
$D$ is a (weakly) linearly convex domain, it follows that
$\rho(\zeta'X)=h(\zeta)+o(|\zeta|^l)\ge 0,$ where
$$h(\zeta)=\sum_{j+k=l}a_{jk}\zeta^j\overline{\zeta}^k\not\equiv 0.$$
Homogeneity of $h$ implies $h\ge 0.$ In addition, as
$h\not\equiv 0,$ we can find $\zeta$ so that $|\zeta|=1$ and $h(\zeta)>c_2>0.$
Now the constant $c_2$
has the required properties for all small $\delta>0.$\qed
\smallskip

Combining Proposition \ref{bound1} and \ref{type4}, we directly get
the following generalization
(in an easy way) of the main result in \cite{Lee} that deals with convex domains.

\begin{corollary}\label{type5} For the notations from Proposition
\ref{type4} we have that if $z\in D\cap U_X\cap n_a,$ then
$$(4c_X)^{-1}(d_D(z))^{-1/l_{a,X}}\le\gamma_D(z;X)\le\kappa_D(z;X)
\le c_X(d_D(z))^{-1/l_{a,X}}.$$
\end{corollary}

The main result in \cite{McN1} (see also \cite{Boas-Str}) states that
$m_a=l_a$ for each convex domain. As an application of Corollary \ref{type5}
we will easily show
something more even for an arbitrary $\CC$-convex domain.

\begin{proposition}\label{type6} If $a$ is a smooth boundary point of a
$\CC$-convex domain $D\subset\CC^n,$
then $m_{a,X}=l_{a,X}$ for each vector $X\neq 0.$

In particular, $m_a=l_a.$
\end{proposition}

\beginproof It suffices to prove that $m_{a,X}\le l_{a,X},$ if
$l_{a,X}<\infty.$ By Corollary \ref{type5} we have
$$\limsup_{D\cap n_a\ni z\to
a}\kappa^s_D(z;X)d^{1/l_{a,X}}\ge\liminf_{D\cap n_a\ni z\to
a}\gamma_D(z;X)d^{1/l_{a,X}}>0$$ (see Section \ref{zeko} for the definition of
$\kappa^s_D$) and the desired inequality follows from
\cite[Corollary]{Yu3}.\qed
\smallskip

The following result is important if the boundary is not real-analytic near a
boundary point of infinite type.

\begin{proposition}\label{type7} If $a$ is a $\mathcal C^1$-smooth
boundary point of a $\CC$-convex domain $D\subset\CC^n,$ then $\partial D$ does not
contain analytic discs through $a$ exactly when $L_a=\{0\}$ (i.e.
$\partial D$ does not contain linear discs through $a$).
\end{proposition}

\beginproof We use the notations from the proof
of Proposition \ref{type4}. It suffices to show that if
$\varphi:\Bbb D\to\partial D$ is an analytic disc, for which
$\varphi(0)=0,$ then $L_a\neq\{0\}.$ As $\partial D$ is $\mathcal
C^1$-smooth near $a,$ there exists $c>0$ so that
$\varphi_\delta(\zeta)=-\delta_n+\varphi(\zeta)\in D$ for
$\delta<c$ and $|\zeta|<c.$ Put $m=\ord_0\varphi$ and
$X=\frac{\varphi^{(m)}(0)}{m!}.$ Then
$\gamma_D(\delta_n;X)\le\kappa_D^s(\delta_n;X)\le1/c$ and as in
the proof of Proposition \ref{bound3} it follows that
$\Delta_X(a,c/4)\subset\partial D.$\qed
\smallskip

\noindent{\bf Remark.} In the case of a convex domain the smoothness condition
is redundant, as seen in the argument of the last remark in the previous section.
\smallskip

Now we will discuss the so-called multitypes of a smooth boundary point
$a$ of a domain $D\subset\C^n.$ For each $k=1,\dots,n$ put
$$m_a^k=\inf_L \sup_{\gamma}\frac{\ord_a(r\circ\gamma)}{\ord_a\gamma},$$ where
$S$ varies over all hyperplanes through $a$ with dimension $k,$ while
$\gamma$ varies over all analytic discs in $S$ that pass through
$a$ (see e.g. \cite{Nik0}). By replacing the analytic discs by complex lines we
define $l_a^k.$ For $k=n$ these numbers
coincide with $m_a$ and $l_a$, respectively. Clearly $l_a^1=m_a^1=1$ and
$l_a^k\le m_a^k.$
The D'Angelo multitype of $a\in\partial D$ is defined as the nondecreasing
$n$-tuple of numbers $M_a=(m_a^1,\dots,m_a^n).$
The D'Angelo linear type $L_a$ is defined in a similar way.
We can also define the Catlin multitype $\tilde M_a=(\tilde m_a^1,\dots,\tilde m_a^n)$
and the Catlin linear multitype
$\tilde L_a=(\tilde l_a^1,\dots,\tilde l_a^n)$ (see e.g. \cite{Yu1}). Note that
$$l_a^k\le\tilde l_a^n\le\tilde m_a^k\le m_a^k.$$

The main result in the voluminous work \cite{Yu1} states that
$\tilde L_a=M_a$ (and so $=\tilde M_a$) for each convex domain.
Using \cite{Yu1} and other nontrivial facts, in \cite{Con} this
equality is proven for $\C$-convex domains.

As a corollary from Proposition \ref{type6} (that we proved easily),
we can get the above results and even to strengthen them a bit.

\begin{proposition}\label{type9} If $a$ a is smooth boundary point of the $\CC$-convex
domain $D\subset\CC^n,$ then $L_a=M_a.$
\end{proposition}

\beginproof We can assume that $a=0.$ We have to show that $l^k_0\ge
m_0^k,$ if $l^k_0<\infty$ and $k>1.$ Let $l^k_0$ be attained for
some $S$ and a line $s\in S.$ If $S$ is orthogonal to the complex
normal $N_0$ to $\partial D$ at $0,$ we consider the subspace $S'$
generated by $N_0$ and a subspace of $S$ of codimension 1,
containing $s.$ Then $D_k=D\cap S'\subset\C^k$ is a $\C$-convex
domain, which is smooth near $0.$ Let $m_{0,k}$ and $l_{0,k}$ be
the type and the linear type of the point $0\in\partial D_k,$
respectively. Then $l_{0,k}=l^k_0,$ as if a line $s'\subset S'$ is
not orthogonal to $N_0,$ then $\ord_0(r\circ s')=1\le l^k_0.$ It
remains to use that $m^k_0\le m_{0,k}$ and $m_{0,k}=l_{0,k}$ by
Proposition \ref{type6}.\qed
\smallskip

Let us mention that a pseudoconvex point $a$ of finite type, for which
$\tilde M_a=M_a$, is called semiregular (see \cite{Die-Her}).
Thus each smooth point of finite type of a $\C$-convex domain is semiregular.

\setcounter{equation}{0}
\section{Estimates for the Bergman kernel and the Bergman metric}\label{ber}

In this section we will prove some estimates for the Bergman
kernel and the Bergman metric of a $\Bbb C$-convex domain
$D\subset\CC^n$ not containing complex lines. The constants in
these estimates depend only on $n.$ The estimate for the Bergman
metric is in the spirit of these for the Carath\'eodory and
Kobayashi metrics from Section \ref{car-kob}. As a corollary we
get that these three metrics are comparable with constants
depending only on $n.$

First recall the definitions for Bergman kernel and Bergman metric for a domain $D\subset\C^n.$
For them and other basic facts see e.g. \cite{Jar-Pfl1}.

Denote by $L_h^2(D)$ The Hilbert space of the square-integrable
holomorphic functions $f$ in $D$. This space has a (unique)
reproducing kernel $\wdtl K_D(z,w)$ -- the Bergman kernel. For
brevity, its restriction $K_D(z)=\wdtl K_D(z,z)$ to the diagonal
is also called Bergman kernel; further we will mainly work with
$K_D.$ It is well-known that $K_D$ is a solution to the following
extremal problem: $$K_D(z)=\sup\{|f(z)|^2:f\in
L_h^2(D),\;\|f\|_D\le1\},$$ where $\|\cdot\|$ is the $L^2$-norm.
If $K_D(z)>0$ for some $z\in D$, then the quadratic form
$$\sum_{j,k}\sum_{j,k=1^n}\frac{\partial^2}{\partial z_j\partial\overline{z_k}}\log K_D(z)X_j\overline{X_k},\quad
X\in\C^n,$$ is positively semidefinite and its square root $B_D(z;X)$ is called Bergman metric.
It also solves an extremal problem:
$$B_D(z;X)=\frac{M_D(z;X)}{\sqrt{K_D(z)}},
$$
where $M_D(z;X)=\sup\{|f'_z(X)|:f\in
L_h^2(D),\,\|f\|_D=1,\;f(z)=0\}$.

Recall that the Carath\'eodory metric does not exceed the Bergman
metric (if the latter one is defined):
$$\gamma_D\le B_D.$$

There are the following transformation rules for the Bergman kernel
and the Bergman metric: if $f:G\to D$ is a biholomorphism
between domains in $\C^n,$ then
$$K_D(f(z),f(w))\Jac f(z)\overline{\Jac f(w)}=K_G(z,w),$$
$$B_D(f(z);f'_z(X))=B_G(z;X).$$
Note that unlike the Carath\'eodory and Kobayashi metrics,
the Bergman metric is not monotone under domain inclusions.
Anyway is is the quotient of two monotone invariants, $M_D$ and $K_D.$

This will help to attain the main goal of this section, namely to
show the converse inequality
to $\gamma_D\le B_D$ up to a constant depending only on $n.$

\begin{theorem}\label{ber1} There exists a constant
$c_n>0$ depending only on $n,$ so that for each $\Bbb C$-convex
domain $D\subset\Bbb C^n$, not containing a complex
line,*\footnote{*Under this assumption $D$ is biholomorphic to a
bounded domain (see Proposition \ref{c.pr4}), so $B_D$ is
defined.} we have the inequality
$$1/4\le B_D(z;X)d_D(z;X)\le c_n.$$
\end{theorem}

By Propositions \ref{bound0} and \ref{bound2}, and by the
inequality $\gamma_D\le B_D,$ we get

\begin{corollary}\label{ber2} There exists a constant
$c_n\ge 1$ depending only on $n,$ so that for each $\Bbb C$-
convex domain $D\subset\Bbb C^n$, not containing a complex line, we have
$$\kappa_D/4\le B_D\le c_n\gamma_D.$$

If $D$ is a convex domain, then the constant 4 can be replaced by 1.
\end{corollary}

The first results, similar to Theorem \ref{ber1} and to Theorem
\ref{ber3} below (for the Bergman kernel $K_D$), refer to bounded
smooth convex domains of finite type \cite{Chen,McN2,McN3}.
Unfortunately the geometric construction there (see also
\cite{McN1}) has a flaw, as we will observe in the next section.
These results are later proven for bounded smooth $\C$-convex
domains of finite type \cite{Blu} using a correct geometric
construction from \cite{Hef1,Hef2,Con} and the paper
\cite{Nik-Pfl3} of the author and P. Pflug (see also
\cite{Die-For3}). Note that the constants in the corresponding
estimates depend on the domains.

Now let us show the most general form of this construction.

Let $D\subset\C^n$ be a domain not containing a complex line. For
a point $z$ we choose $a^1\in\partial D$ so that
$d_1:=||a^1-z||=d_D(z).$ Put $H_1=z+\span(a^1-z)^\bot$ and
$D_1=D\cap H_1.$ Let $a^2\in\partial D_1$ so that
$d_2:=||a^2-z||=d_{D_1}(z).$ Put $H_2=z+\span(a^1-z,a^2-z)^\bot,$
$D_2=D\cap H_2$ and so on. Thus we get an orthonormal basis of the
vectors $e_j=\frac{a^j-z}{||a^j-z||},$ $1\le j\le n,$ which will
be called minimal (for $D$ at $z$) and positive numbers $d_1\le
d_2\le\dots\le d_n$ (the basis and the numbers are not uniquely
determined).

Put  $$p_D(z)=d_1\dots d_n.$$

The lower estimate for the Bergman kernel $K_D$ via $p_D$ in the
next theorem is a main point in the proof of Theorem \ref{ber1},
but is also of independent interest.

\begin{theorem}\label{ber3} Let $D\subset\Bbb C^n$ be a
$\Bbb C$-convex domain not containing a complex line. Then
$$\frac{1}{(16\pi)^n}\le K_D(z)p_D^{2}(z)\le\frac{(2n)!}{(2\pi)^n}.$$
In addition, the lower estimate is precise for $n=1,$ while the
upper estimate is exact for each $n$ (even for convex domains);
the inequality is strict for $n\ge 2.$

In addition, if $D$ is a convex domain not containing complex
lines, the lower estimate can be improved by replacing the number
16 by 4. In this case the estimate is precise for each $n.$
\end{theorem}

\beginproof{\it The upper estimate.} We can assume that $z=0.$ By Lemma
\ref{bound11},
$$D\supset G=\{z\in\CC^n:\sum_{j=1}^n\frac{|z_j|}{d_j}<1\}.$$
Consequently $G\cup\Bbb B_n(0,d_1)\subset D$ and so
$$K_D(0)\le K_{G\cup
\Bbb B_n(0,d_1)}(0)\le K_G(0)=K_E(0)/p_D^2(0),$$ where
$$E=\{z\in\CC^n:\sum_{j=1}^n|z_j|<1\}$$ (here we applied the transformation
rule for the Bergman kernel to the dilatation of the coordinates
$(z_1,\dots,z_n)\to (z_1/d_1,\dots,z_n/d_n)$). As $E$ is a
complete Reinhardt domain, $K_E(0)=\vol(E)^{-1}.$ It is easily
calculated that this volume equals
$\frac{(2\pi)^n}{(2n)!}{(2\pi)^n},$ thereby proving the upper
estimate.

It is precise for $n=1,$ as seen in the example of the unit disc
and its center. If $n\ge 2$, then $G$ does not contain $\Bbb
B_n(0,d_1)$ so the second inequality above is strict (since the
volume of $G$ is less than that of $G\cup \Bbb B_n(0,d_1)$).

To finish the discussion about the upper estimate, it remains to
show that it is precise for $n\ge 2.$ For $m\in\Bbb N$ put
$b_j=j^m$ for $1\le j\le n$. Let $B_j=\Bbb B_n(0,b_j)\cap
H'_{j-1},$ where $H'_{j-1}=\{0\}\times \CC^{n-j+1}.$ Denote by $T$
the convex hull of the union $\cup_{j=1}^n B_j$ and the domain
$\{z\in\CC^n:\sum_{j=1}^n|z_j|/b_j<1\}$. It is not hard to see
that $b_j=\dist(0,\partial(T\cap H'_{j-1}))$.

Later if $\Psi(z)=(z_1/b_1,\ldots,z_n/b_n)$, then $\Psi(T)$ is the
convex hull of the union $S=\cup_{k=1}^n\Psi(B_k)$ and $E$. For
each $k>j$ we have $b_k/b_j\to\infty$ if $m\to\infty$.
Consequently for an arbitrary $\lambda>1$ one can find an $m$ such
that $S\subset\lambda E$. As $\lambda E$ is a convex domain, it
contains $\Psi(T)$. So
$$K_T(0)(b_1\ldots b_n)^2=K_{\Psi(T)}(0)\ge
K_{\lambda E}(0)= \frac{K_E(0)}{\lambda^{2n}}$$ and as
$\lambda>1,$ the upper estimate is precise.
\smallskip

{\it The lower estimate.}*\footnote{*The geometric proof is close
to that of Proposition \ref{c.pr4}.} After a translation and a
rotation, we can assume that $z=0$, $H_j=\{0\}\times \CC^{n-j}$
($j=1,\dots,n-1$) and
$a^j=(0,a^j_j,0)\in\CC^{j-1}\times\CC\times\CC^{n-j}$
($j=1,\dots,n$) so that $d_j=|a_j^j|$.

As $D$ is a $\CC$-convex domain, there exists a hyperplane
$a^j+W_{j-1}$ through $a^j$ that is disjoint from $D$. By our
construction the ball in $H_1$ of center $0$ and radius $a^2_2$
lies in $D$ and so $W_1\cap H_1$ is orthogonal to $a^2$, i.e.
$W_1\cap H_1\subset\{0\}\times\CC^{n-2}$. Consequently $W_1$ is
defined by the equation $\alpha_{1,1}z_1+z_2=0$. The same argument
shows that the equation of $W_j$ for $j=0,\dots,n-1$ is the
following one:
$$
\alpha_{j,1}z_1+\cdots+\alpha_{j,j}z_j+z_{j+1}=0.
$$

Let $F:\CC^n\to\CC^n$ be the linear mapping of matrix $A$ having
as vector-rows $(\alpha_{j,1},\dots,\alpha_{j,j},1,0,\dots,0),$
$j=0,\dots,n-1.$ Then $G=F(D)$ is also a $\CC$-convex domain ($G$
was another domain in the proof of the upper estimate). Note that
$K_D(0)=K_G(0),$ as $\det A=1.$ Put $G_j=\pi_j(G)$, where $\pi_j$
is the projection onto the $j$-th coordinate plane. Then $G_j$ is
a simply connected domain (see e.g. \cite{APS}) and $G\subset
G_1\times\cdots\times G_n$. Consequently
\begin{equation}\label{br1}
K_D(0)\geq K_{G_1\times\cdots\times G_n}(0)=K_{G_1}(0)\cdots
K_{G_n}(0).
\end{equation}
As $G_j\neq\C$ is a simply connected domain, it is biholomorphic
to $\D$ and $\sqrt{\pi K_\D(0)}=1=\gamma_\D(0;1)$ implies that
\begin{equation}\label{br2}
\sqrt{\pi K_{G_j}(0)}=\gamma_{G_j}(0;1).
\end{equation}
On the other hand,
\begin{equation}\label{br3}
\gamma_{G_j}(0;1)\geq \frac{1}{4d_{G_j}(0)}
\end{equation}
by the K\"obe $1/4$ theorem (this argument was already used in the
proof of Proposition \ref{bound1}).

Later, $F(a^j)\in\partial G$, and the $j$-th coordinate of this
point is $a_j^j.$ In addition, the hyperplane $\{z\in\Bbb
C^n:z_j=a_j^j\}$ does not intersect $G$. Consequently
$a_j^j\in\partial G_j$; in particular
$$d_j=|a_j^j|\geq d_{G_j}(0).$$
This together with (\ref{br1}), (\ref{br2}) and (\ref{br3}) proves
the lower estimate.

Note that the constant 16 is the best possible for $n=1,$ as shown
by the example of the image $D=\CC\setminus[1/4,\infty)$ of $\D$
under the K\"obe transformation  $z\to z/(1+z)^2$ (already used in
Section \ref{car-kob}). This example is not applicable for $n\ge
2$ in a trivial way, since a $\C$-convex Cartesian product of
domains that are different from $\C$ is necessarily convex (see
the first remark in Section \ref{con}).

For the lower estimate in the case of a convex domain is
sufficient to note that $G_j$ are convex domains, so the number 4
can be replaced by 2 due to Proposition \ref{bound0}.

Finally note that in this case the constant 4 is the best
possible, as seen in the example of a Cartesian product of
half-planes.\qed
\smallskip

Using the lower estimate in Theorem \ref{ber3}, we can now prove Theorem \ref{ber2}.
\smallskip

\noindent{\it Proof of Theorem \ref{ber2}}. We will use the geometric
configuration from the proof of Theorem
\ref{ber3}.

Let $X\in(\CC^n)_\ast.$ First we will find an upper estimate for
$M_D(0;X).$ Fix a $k\in J=\{j:X_j\neq 0\}$. Then
$$\Psi_k(z)=(z_1-\frac{X_1}{X_k}z_k,\ldots,z_{k-1}-\frac{X_{k-1}}{X_k}z_{k},
z_k,z_{k+1}-\frac{X_{k+1}}{X_k}z_k,\ldots,z_n-\frac{X_n}{X_k}z_k)$$
is a linear mapping of Jacobian 1 and
$$Y^k:=\Psi_k(X)=(0,\ldots,0,X_k,0,\ldots,0).$$ Let $\Delta_j$ is
the disc in the $j$-th coordinate plane of center $0$ and radius
$d_j$ for $j\neq k$, and radius $d'_k=|X_k|d_D(0,X)$ for $j=k$.
Then $\Delta_j\subset D_k=\Psi_k(D)$ and by Lemma \ref{bound11},
$$D_k\supset E_k=\{z\in\CC^n:\frac{|z_k|}{d'_k}+\sum_{j=1,j\neq
k}^n\frac{|z_j|}{d_j}<1\}.$$ Consequently
$$M_D(0;X)=M_{D_k}(0;Y^k)\le
M_{E_k}(0;Y^k)=\frac{C_nd_{k,D}(0)}{|X_k|p_D(0)d_D^2(0,X)},$$
where $C_n:=M_E(0;e_1)=\sqrt{\frac{(2(n+1))!}{6(2\pi)^n}}$
(the last is calculated directly, as $E$ is a complete Reinhardt domain),
and $e_1$ is the first basis vector.

From this estimate and the lower estimate in Theorem \ref{ber3} it follows that
\begin{equation}\label{br5}B_D(0;X)=\frac{M_D(0;X)}{\sqrt{K_D(0)}}\le
\frac{c'_nd_{k,D}(0)}{|X_k|d_D^2(0,X)}, \quad 1\le k\le n,
\end{equation}
where
$c'_n=(4\sqrt{\pi})^nC_n=2^n\sqrt{\frac{2^{n-1}(2(n+1))!}{3}}.$

It remains to note that Lemma \ref{bound11} implies the inequality
\begin{equation}\label{br6}
\frac{1}{d_D(0,X)}\le\sum_j^n\frac{|X_j|}{d_j}
\end{equation}
and then put $c_n=nc'_n.$\qed
\smallskip

The above results and their proofs allow us to understand the
boundary behavior of any of the considered metrics $F_D$ -- of
Carath\'eodory, Kobayashi or Bergman -- of a $\C$-convex domain
that does not contain a complex line, in terms of minimal bases.
This strengthens some results from \cite{Chen,McN2,McN3,Blu,Lie},
dealing with bounded smooth domains of finite type; the constants
there depend on the domain (the first three works even refer to
convex domains).

\begin{proposition}\label{ber4} There exists a constant
$c_n\ge 1,$ depending only on $n,$ so that for each $\Bbb
C$-convex domain $D\subset\Bbb C^n$, not containing a complex line, we have
$$c_n^{-1}\le F_D(z;X)\left(\sum_j^n\frac{|\langle X,e_j(z)\rangle|}
{d_j(z)}\right)^{-1}\le c_n.$$ (Here $e_j(z)$ are the basis
vectors of a minimal basis of $D$ at $z,$ and $d_j(z)$ are the
corresponding numbers.)
\end{proposition}

\beginproof By (\ref{br5}) and the inequality
$$B_D(z;X)\ge\frac{1}{4d_D(z;X)}$$
we get
$$\frac{|X_j(z)|}{d_j(z)}\le\frac{4c_n'}{d_D(z)}.$$
So
$$\frac{1}{d_D(z;X)}\le\sum_j^n\frac{|X_j(z)|}{d_j(z)}\le\frac{4c_n}{d_D(z;X)},$$
where $c_n=nc_n'.$ Then (\ref{br5}) and (\ref{br6}) show that
$$(16c_n)^{-1}\le F_D(z;X)\left(\sum_j^n\frac{|X_j(z)|}{d_j(z)}\right)^{-1}\le
c_n.\qed$$
\smallskip

The following result is in the spirit of Proposition \ref{ber4}, but it deals with a fixed basis.
As each boundary point of a bounded $\C$-convex domain is semiregular (see the end of Section
\ref{type}), it directly follows from \cite{Yu4,BSY} and \cite[Theorem 3.3.1]{Nik0}.

\begin{proposition}\label{ber5} Let $a$ be a boundary point of finite
type of a bounded smooth $\C$-convex domain $D\subset\C^n.$ Denote
by $\tilde M_a=(m_1,\dots,m_n)$ the Catlin multitype of $a$. Then
there exists a linear basis change with the following property:
for each non-tangent cone $\Lambda$ with vertex $a$ there exists a
constant $c>0$ so that for an arbitrary vector $X=(X_1,\dots,X_n)$
in the new basis we have
$$c^{-1}\le \liminf_{\Lambda\ni z\to
a}F_D(z;X)\left(\sum_j^n\frac{|X_j|}{(d_D(z))^{1/m_j}}\right)^{-1}$$
$$\le\limsup_{\Lambda\ni z\to
a}F_D(z;X)\left(\sum_j^n\frac{|X_j|}{(d_D(z))^{1/m_j}}\right)^{-1}\le
c.$$

In addition, for the Bergman kernel we have
$$c^{-1}\le\liminf_{\Lambda\ni z\to a}K_D(z)
(d_D(z))^{-2q}\le\limsup_{\Lambda\ni z\to
a}K_D(z)(d_D(z))^{-2q}\le c,$$ where $q=1/m_+\dots+1/m_n.$
\end{proposition}

The basis change can be chosen to be linear because $\tilde L_a=\tilde M_a.$

Note that Proposition \ref{ber5} implies a more precise variant of
Proposition \ref{type4} in the case of finite type, namely that
for each vector $X\in(\C^n)_\ast$ there exists $j=1,\dots,n$ so
that $l_{a,X}=m_j.$

Finally let us mention that by a result from \cite{Con} the
quotient $(d_{j,D})^m/d_D$ is bounded near $a$ (recall that
$m=m_n$ is the type of $a$) and then Proposition \ref{ber4}
provides a neighborhood $U$ of $a$ and a constant $c>0$ so that
$$\kappa_D(z;X)\ge\frac{c||X||}{(d_D(z))^{1/m}},\quad z\in D\cap U.$$

\setcounter{equation}{0}
\section{Maximal basis. A counterexample}\label{count}

To get estimates for the Bergman kernel and the Bergman metric for
$\C$-convex domains, in the previous section we introduced a basis
(called minimal) with an origin at a given point from a given
domain. As mentioned, in the special case of smooth convex domain
of finite type, in \cite{Chen,McN2,McN3} it is introduced a
similar basis (that we call maximal). The minimal and maximal
bases can be considered in the context of the so-called extremal
bases (see \cite{Char-Dup}). Many other important results, like
those connected with the linear type, with the
$\overline{\partial}$-problem or with domains with noncompact
groups of automorphisms (see e.g.
\cite{McN1,McN-Ste1,McN-Ste2,Gau1}), use in an essential way the
properties of the maximal basis, and most of all one extremal
property, satisfied also by the minimal basis. In general, this
property means that the vectors from the basis are orthogonal to
the corresponding hyperplanes (see the Introduction; the details
are given below). Unfortunately exactly this property of the
maximal basis turns out to be wrong (the hints for corresponding
proofs are based on an incorrect application of the method of
Lagrange multipliers).

The main aim of this section is providing a counterexample
for the extremal property of the maximal basis.

Now we define the notion "maximal basis". Let $D\subset\C^n$ be a
domain, not containing a complex line. For point $q\in D$ we
choose a unit vector $a_1\in\C^n$ so that
$$s_1:=d_D(q;a_1)=d_D(q).$$ Then we choose a unit vector
$a_2\in\span(a_1)^\perp$ so that
$$s_2:=d_D(q;a_2)=\sup d_D(q;a),$$
where the supremum is taken over all the unit vectors $a\in
\span(a_1)^\bot.$ On the next step we choose a unit vector
$a_3\in\span(a_1,a_2)^\perp$ so that
$$s_3:=d_D(q;a_3)=\sup d_D(q;a),$$
where the supremum is taken over the unit vectors $a\in
\span(a_1,a_2)^\bot.$ Continuing the procedure, we get a
orthonormal basis $a_1,\dots,a_n$ that will be called maximal (for
$D$ at $q$) and a sequence of positive numbers $s_2\ge\dots\ge
s_n\ge s_1\ge 0$ (but they are not uniquely determined). Note
that, unlike the minimal basis, after the first step the
corresponding distances are chosen maximal (rather than minimal).

Assume now that $D$ is a convex domain that is smooth near a
boundary point $p_1$ (of finite type). Let $r$ be a locally
defining function. Now we will describe the extremal property
mentioned in the Introduction. For $q\in D$ on the internal normal
to $\partial D$ at $p_1$, sufficiently close to $p_1$, we consider
a coordinate system defined by the maximal basis at $q$, i.e. we
put $q=0$ and express each $z\in\CC^n$ in the form $z=\sum_{j=1}^n
w_ja_j$. We choose $p_k\in\partial D,$ $k=2,\dots,n$ so that
$p_k=\lambda_ka_k,$ where $|\lambda_k|=s_k.$ Many of the works
mentioned in the Introduction (see e.g. \cite[Proposition 2.2
(ii)]{Chen}, \cite[Proposition 3.1 (i)]{McN1}, \cite[Proposition
2.1 (iii)]{McN2}) note that
$$
(\ast)\quad\frac{\partial r(p_k)}{\partial w_j}=0,\quad j=k+1,\dots,n.
$$
This means that
$$(\ast\ast)\quad T_{p_k}^\CC(\partial D)\cap\span(a_1,\dots,a_k)^\bot=\span(a_{k+1},\dots,a_n).$$
Note that the minimal basis has a very essential property; we started by it when obtaining
the lower estimate for the Bergman kernel in Theorem \ref{ber3}.

However now we will demonstrate a counterexample in $\CC^3$ for
the property $(\ast)$ of the maximal bases at the points from an
internal normal to the boundary of a domain in $\C^3$ (in $\C^2$
this property clearly holds).

Let $0<\beta_2<\beta_1<1$. Put
$$
D=\{z\in\CC^2\times\CC:\rho(z)+|z_3|^2<1\},
$$
where $\rho(z)=x_1^2+\beta_1y_1^2+x_2^2+\beta_2y_2^2.$ Note that
$D$ is a strictly (pseudo)convex domain with a real-analytic
boundary. Let $q=(0,0,\delta),$ where $0<\delta<1.$ The
construction of a maximal basis of $D$ at $q$ leads to
$s_1=1-\delta$ and $a_1=(0,0,1)$. From the next step we get the
domain
$$
D_\delta=\{z\in\CC^2:\rho(z)<1-\delta^2\}.
$$
Note that after a homothety $D_\delta$ goes to $D_0$ and so we can
examine only $D_0$. For the second vector $a_2$ from the maximal
basis we have $a_2=(b,0),$ where $b \in\CC^2$. Then
$\span(a_1,a_2)^\bot$ is generated by $(-\overline b_2, \overline
b_1, 0)$. Put $$\mathcal T=\{b\in\CC^2:\frac{\partial
\rho(b)}{\partial z_1}(-\overline b_2)+ \frac{\partial
\rho(b)}{\partial z_2}(\overline b_1)=0\}.$$

\begin{lemma} \label{lemma}\label{count1}
$\mathcal T=\{b\in\CC^2: b_1=0 \text{ or } b_2=0 \text{ or } \Im
b_1=\Im b_2=0\}$.
\end{lemma}

\beginproof
Some elementary calculations show that $b\in\mathcal T$ exactly when
$$\begin{aligned}(\beta_1-\beta_2)\Im b_1\Im b_2&=0\\
(1-\beta_1)\Im b_1\Re b_2&=(1-\beta_2)\Im b_2\Re b_1,
\end{aligned}$$
and the result follows.\qed
\smallskip

Let $p_2\in\partial D_0$ so that
$$\frac{d_{D_0}(0;p_2)}{\|p_2\|}=s_2=\sup_{\|a\|=1}d_{D_0}(0;a).$$
The following result shows that the property $(\ast\ast)$,
equivalent of $(\ast)$, is not true at the points on the internal
normal to $D$ at $(1,0,0)$ (formally we must also consider the
case $\delta<0,$ but then the closest point is $(-1,0,0)$ and the
situation is similar.)

\begin{proposition}$p_2\not\in\mathcal T.$\label{count2}
\end{proposition}

\beginproof
Let $b\in\mathcal T$ be a unit vector. Note that $\rho(re^{i\alpha}b)<1$
for each $\alpha\in\RR$ if and only if
$r^2R(b)<1,$ where $R(b)=\max_{\alpha\in\RR}\rho(e^{i\alpha}b).$ Consequently
$d_{D_0}(0;b)=1/\sqrt{R(b)}$. Let $b=(e^{i\phi_1}\cos\Theta,e^{i\phi_2}\sin\Theta)$, where
$0\leq\Theta<2\pi$ and $0\le\phi_1,\phi_2\le\pi/2.$ By Lemma
\ref{count1} there are three possibilities for $b$:

$\bullet$ $\Theta=0$ or $\Theta=\pi$: $\rho(e^{i\alpha}b)=\cos^2(\alpha+\phi_1)+\beta_1\sin^2(\alpha+\phi_1).$

$\bullet$ $\Theta=\pi/2$ or $\Theta=3\pi/2$: $\rho(e^{i\alpha}b)=\cos^2(\alpha+\phi_2)+\beta_2\sin^2(\alpha+\phi_2).$

$\bullet$ $\phi_1=\phi_2=0$: $\rho(e^{i\alpha}b)=\cos^2\alpha+\sin^2\alpha(\beta_1\cos^2\Theta+\beta_2\sin^2\Theta).$

\noindent In all three cases $R(b)=1.$

On the other hand, there exists a unit vector $b^\ast\in\Bbb C^2$ such that
$R(b^\ast)<1,$ and so $p_2\not\in\mathcal T$. To define $b^\ast$, put $\Theta=\pi/4$ $\phi_1=0$
and $\phi_2=\pi/2.$ Then $2\rho(e^{i\alpha}b^\ast)=1+\beta_2+(\beta_1-\beta_2)\sin^2\alpha.$
As $\beta_1<\beta_2<1$, we conclude that
$R(b^\ast)=\frac{1+\beta_2}{2}<1.$\qed

\setcounter{equation}{0}
\section{Estimates in a maximal basis}\label{max}

The aim of this section is, using the estimates for invariants in
terms of a minimal basis, to prove that they remain true in terms
of a maximal basis, in spite of the counterexample from the last
section. A similar approach allows to confirm the correctness of
other results using maximal basis.

Let $D\subset\CC^n$ be a $\C$-convex domain not containing a
complex line (i.e. each nonempty intersection of $D$ with a
complex line is biholomorphic to $\D$). For $z\in D$, let
$e_1,\dots,e_n$ be a minimal basis of $D$ at $z,$ and
$a_1,\dots,a_n$ -- a reordered maximal basis of $D$ at $z,$
meaning that the new $a_1$ is the old $\wdtl a_1,$ but $a_2=\wdtl
a_n, a_3=\wdtl a_{n-1},\dots,a_n=\wdtl a_2.$ Let $d_1\le\dots\le
d_n$ and $s_1\le\dots\le s_n$ be the corresponding numbers (recall
that $d_1=s_1=d_D(z)$). Put $p_D(z)=\prod_{j=1}^n d_j$ and
$s_D(z)=\prod_{j=1}^ns_j.$ As before, $K_D(z)$ denotes the Bergman
kernel (on the diagonal). Let $F_D(z;X)$ be any of the metrics of
Carath\'eodory, Kobayashi or Bergman. For $X\in\C^n$ put
$$E_D(z;X)=\sum_{j=1}^n\frac{|\langle X,e_j\rangle|}{d_j},\quad
A_D(z;X)=\sum_{j=1}^n\frac{|\langle X,a_j\rangle|}{s_j}.$$ We will write
$f(z)\lesssim g(z)$ if $f(z)\le c_n g(z)$ for some
constant $c_n>0$ depending only on $n;$ $f(z)\sim g(z)$ means that
$f(z)\lesssim g(z)\lesssim f(z).$ By Proposition \ref{bound1},
Theorems \ref{ber2}, \ref{ber3} and Proposition \ref{ber4} we know that
$$K_D(z)\sim 1/p^2_D(z),\quad F_D(z;X)\sim E_D(z;X)\sim 1/d_D(z;X)$$
(as noted, under the much stronger requirements that the domain be $\C$-convex,
smooth, bounded and of finite type these estimates
follow also from \cite{Blu,Lie}). For brevity we will sometimes omit the
arguments $z$ and $X$. Lemma \ref{bound11} easily implies that
$$K_D\lesssim 1/s^2_D,\quad F_D\lesssim A_D.$$ In particular,
$$1/d_D(z;X)\sim E_D(z;X)\lesssim A_D(z;X)$$
As mentioned, the main corollary from the incorrect property $(\ast)$
for the maximal bases (for a bounded smooth
$\C$-convex domain of finite type) is the fact that $$A_D(z;X)\sim_D 1/d_D(z;X),$$
where the constant in $\sim_D$ depends on $D.$ Based on this fact,
in \cite{Chen,McN2,McN3} it is shown that
$$K_D\sim_D 1/s_D^2,\quad F_D\sim_D A_D.$$
The next two propositions show that anyway these estimates are correct.

The first estimate can be also obtained from \cite{Hef2} in the case of a
bounded smooth $\C$-convex domain of finite type.
The proof there uses the incorrectly proven estimate $1/d_D(z;X)\sim_D A_D(z;X),$
but a closer look shows that
one can use only the correct part of that estimate, $1/d_D(z;X)\lesssim_D A_D(z;X).$

\begin{proposition}\label{max1} Let $D\subset\C^n$ be a $\C$-convex domain,
not containing complex lines. Then for each point $z\in D$ we have $d_j\sim s_j$, $j=1,\dots,n.$
\end{proposition}

Once again observe that the constant in $\sim $ depends only on the dimension $n$ of $D.$
\smallskip

\beginproof We first prove that $s_j\lesssim d_j.$
As $E_D\lesssim A_D,$ it suffices to understand that if $E_D\le cA_D,$
then $s_j\le c' d_j,$ where $c'=n!c.$

The formula for the determinant of the unitary transformation
between the two bases implies that $\prod_{j=1}^n|\langle
a_j,e_{\sigma(j)}\rangle|\ge 1/n!$ for some permutation of
$\sigma$ of $\{1,\dots,n\}$ In particular, $|\langle
a_j,e_{\sigma(j)}\rangle|\ge 1/n!.$ Then $E_D(z;a_j)\le c
A_D(z;a_j)$ implies $s_j\le c'd_{\sigma(j)}.$

Suppose now that $c'd_k<s_k$ for some $k.$ Then
$$c'd_k<s_k\le s_j\le c'd_{\sigma(j)},\quad j\ge k.$$
Consequently $\sigma(j)>k$ for each $j\ge k,$
a contradiction, as $\sigma$ is a permutation.

These arguments show that $\wdtl d_j\sim d_j,$ where $\wdtl d_j$
are the corresponding numbers for another minimal basis of $D$ at
$z.$ Thus we can assume that $e_1=a_1.$ We know that $s_1=d_1.$ It
remains to show that $s_k\gtrsim d_k$ for $k\ge 2.$ Choose a unit
vector in $\span(e_k,\dots,e_n)$ that is orthogonal to
$a_{k+1},\dots, a_n$ ($a_n'=e_n$ if $k=n$). Then $a_k'$ is also
orthogonal to $a_1=e_1.$ Consequently $s_k\ge d_D(z;a_k')$ ( by
the construction of a maximal basis). On the other hand, as $a_k'$
is orthogonal to $e_1,\dots,e_{k-1},$ then
$$\frac{1}{d_D(z;a_k')}\sim E_D(z;a_k')=
\sum_{j=k}^n\frac{|\langle
a_k',e_j\rangle|}{d_j}\lesssim\frac{1}{d_k}.$$ So $s_k\ge
d_D(z;a_k')\gtrsim d_k.$\qed
\smallskip

\begin{proposition}\label{max2} Let $D$ be as in Proposition \ref{max2}.
Then $A_D\sim E_D.$
\end{proposition}

\beginproof In view of the inequality $E_D\lesssim A_D$ and
Proposition \ref{max1} ($s_k\sim d_k$), it suffices
to prove that $$\frac{|\langle X,a_k\rangle|}{d_k}\lesssim E_D (z;X)$$
for each $k.$

Put $b_{jk}=\langle a_j,e_k\rangle.$ As
$$\frac{1}{d_j}\sim\frac{1}{d_D(z;a_j)}\sim
E_D(z;a_j)\ge\frac{|b_{jk}|}{d_k},$$ it follows that
$|b_{jk}|\lesssim d_k/d_j.$ The unitary transformation with matrix
$B=(b_{jk})$ transforms the basis $e_1,\dots,e_n$ into the basis
$a_1,\dots,a_n$. For the inverse matrix $C=(c_{jk})$ we have
$$
|c_{jk}|\le\sum_{\sigma}|b_{1\sigma(1)}\dots
b_{k-1,\sigma(k-1)}b_{k+1,\sigma(k+1)}\dots b_{n,\sigma(n)}|
$$
$$
\lesssim\sum_{\sigma}\frac{d_{\sigma(1)}}{d_1}\dots\frac{d_{\sigma(k-1)}}{d_{k-1}}
\frac{d_{\sigma(k+1)}}{d_{k+1}}\dots\frac{d_{\sigma(n)}}{d_n}=
\sum_{\sigma}\frac{d_k}{d_j}=(n-1)!\frac{d_k}{d_j},
$$
where $\sigma$ varies over all bijections between the sets $\{
1,\dots, k-1, k+1, \dots, n\}$ and $\{ 1,\dots, j-1, j+1, \dots,
n\}$.

Consequently
$$
\frac{|\langle X,a_k\rangle|}{d_k}\le \sum_{j=1}^n|\langle
X,e_j\rangle|\frac{|b_{kj}|}{d_k} =\sum_{j=1}^n|\langle
X,e_j\rangle|\frac{|\overline c_{jk}|}{d_k}\lesssim E_D.\qed$$

\noindent{\bf Remark.} The constructions of minimal and maximal
bases can be generalized in the following way: we choose
"minimal"\ discs at steps $1,\dots,k$ and "maximal"\ discs at
steps $k+1,\dots, n-1$ (the $n$-th choice is canonical); $k=n-1$
gives a minimal basis, $k=1$ -- a maximal one, and $k=0$ -- a
basis without "minimal"\ discs. Note that Propositions \ref{max1}
and \ref{max2} remain true when $s_j$ are substituted by the
numbers in the new basis and we express $A_D$ in this basis. (This
construction has an obvious real analogue).
\smallskip

\setcounter{equation}{0}
\section{Localizations}\label{loc}
It is natural to ask whether the results from the preceding
sections have local character, i.e. whether $\C$-convexity is a
local notion (like convexity and pseudoconvexity) and whether the
behavior of the considered invariant metrics near a boundary point
of a given domain is similar to that of an intersection of the
domain with a neighborhood of the points.

It is hard to get localization results for the Carath\'eodory
metric and here we are not going to deal with them. Anyway such
results can be found in \cite{Nik0}.

First we will discuss the local character of $\C$-convexity. As
noted in Section \ref{con}, each bounded $\C$-convex domain is
(weakly) linearly convex, and the converse is true under the
additional assumption of a $\mathcal C^1$-smooth boundary.

The next proposition shows that this fact has local character.

\begin{proposition}\label{loc-con} Let $a$ be a $\mathcal C^k$-smooth
boundary point ($2\le k\le\infty$) of a domain $D\subset\C^n$ that is
locally weakly linearly convex near $a,$*\footnote{*cf.
Section \ref{con}} i.e. for each $b\in\partial D$ near $a$
there exists a neighborhood $U_b$ so that $D\cap U_b\cap T^\CC_b(\partial
D)=\varnothing.$ Then there exists a neighborhood $U$ of $a,$ for which
$D\cap U$ is a $\mathcal C^k$-smooth $\Bbb C$-convex domain.
\end{proposition}

Clearly this proposition has "convex"\ and "pseudoconvex"\ analogues,
proven in a similar, but easier, way.
\smallskip

\beginproof We can assume that $a=0.$ Denote by $H_f(z;X)$ the Hessian
of a $\mathcal C^2$-smooth function $f$. Put $B_s=\Bbb
B_n(0,s)$ ($s>0$) and $$r(z)=\left\{\begin{array}{ll}
-d_D(z),&z\in D\\
d_D(z),&z\not\in D.\end{array}\right.$$ The differential
inequality for $r^2$ in the proof of \cite[Proposition 2.5.18
(ii)$\Rightarrow$(iii)]{APS} easily implies that there exists an
$\eps>0$ so that $r$ is a $\mathcal C^k$-smooth defining function
of $D$ at $B_{3\eps}$ and $H_r(z;X)\ge 0,$ if $\langle\partial
r(z),\overline X\rangle=0$ and $z\in D\cap B_{2\eps}.$ Then the
proof of \cite[Lemma 1]{Die-For1} shows that there exists a $c>0$
such that $H_r(z;X)\ge-c||X||\cdot|\langle\partial r(z),\overline
X\rangle|,$ $z\in D\cap B_{2\eps}.$ We can suppose that $2\eps
c\le 1$ and $D\cap B_\eps$ is connected. Choose a smooth function
$\chi$ so that $\chi(x)=0$ for $x\le\eps^2$ and
$\chi'(x),\chi''(x)>0$ for $x>\eps^2.$ Put
$\theta(z)=\chi(||z||^2).$ We can find a $C\ge 1/2$ such that
$$B_{2\eps}\Supset G'=\{z\in B_{2\eps}:0>\rho(z)=r(z)+C\theta(z)\}\subset D.$$

Later, the inequalities $2c\eps\le1$ and $|\langle\partial\theta(z),\overline
X\rangle|\le\chi'(||z||^2)||z||\cdot||X||$ yield
$\chi'(||z||^2)||X||>c|\langle\partial\theta(z),\overline X)|,$
if $z\in B_{2\eps}\setminus\overline{B_\eps}$ and $X\neq 0.$ This,
together with $$H_r(z;X)\ge-c||X||\cdot|\langle\partial r(z),\overline X\rangle|,\quad z\in G',$$
$$
H_{\rho}(z;X)=H_r(z;X)+4C\chi''(||z||^2) \mbox{Re}^2\langle
z,X\rangle+2C\chi'(||z||^2)||X||^2,
$$
$C\ge 1/2,$ and the triangle inequality, implies that
$$H_{\rho}(z;X)\ge-c||X||\cdot|\langle\partial\rho(z),\overline X\rangle|,
\quad z\in\overline{G'}.$$
In addition, the last inequality is strict, if
$z\in\overline{G'}\setminus\overline{B_\eps}$ and $X\neq 0.$ This
shows that $\partial\rho\neq 0$ on $\partial
G'\setminus\overline{B_\eps};$ (otherwise $\rho$ would attain a local minimum at
some point of this set, which is clearly impossible). Thus $\partial\rho\neq 0$ on $\partial G'.$

Let $G$ be a connected component of $G'$ that contains $D\cap
B_\eps.$ Then \cite[Proposition 2.5.18]{APS} (see also
\cite[Proposition 4.6.4]{Hor2}) implies that $G$ is a $\mathcal
C^k$-smooth $\Bbb C$-convex domain. It remains to put
$U=B_n(0,\eps)\cup G.$\qed
\smallskip

Now we will discuss the localization of the Kobayashi metric. First
recall that if $D$ is a hyperbolic domain (i.e. the Kobayashi pseudodistance $k_D$ is a distance),
then the following weak localization takes place (see e.g. \cite{Nik0}):

\begin{proposition}\label{locw-kob} If $V\Subset U$ are neighborhoods of a boundary point
of a hyperbolic domain $D\subset\C^n,$ then there exists a constant
$C\ge 1$ such that for each $z\in D\cap V$ and each $X\in\Bbb C^n$ we have
$$\kappa_D(z;X)\le\kappa_{D\cap U}(z;X)\le C\kappa_D(z;X).$$
\end{proposition}

Propositions \ref{loc-con} and \ref{locw-kob} show that all the
above results for the Kobayashi metric, as well as those connected
with types and multitypes, have local character in the case of
bounded domains.

To this aim note the following. If $a$ is a boundary point of a
bounded domain $D\subset\C^n,$ then it is easily seen that for
each neighborhood $U$ of $a$ we have
$$p_D(z)\sim_\ast p_{D\cap U}(z),\ s_D(z)\sim_\ast
s_D(z),\ d_D(z;X)\sim_\ast d_{D\cap U}(z;X),$$
$$E_D(z;X)\sim_\ast E_{D\cap U}(z;X),\ A_D(z;X)\sim_\ast A_{D\cap
U}(z;X)$$ near $a;$ here the constant in $\sim_\ast$ depends on $D$ and $U.$

Thus we get the following corollary of Propositions
\ref{bound1}, \ref{loc-con} and \ref{locw-kob}.

\begin{corollary}\label{loc-cor}
Let $a$ be a boundary point of a bounded domain $D\subset\C^n$,
as in Proposition \ref{loc-con}. Then
$$\kappa_D(z;X)\sim_D d_D(z;X)\sim_D E_D(z;X)$$
near $a.$ (the constant $\sim_D$ depends on $D$).
\end{corollary}

Recall that the constant $\sim_D$ depends on $D.$

Now we will sharpen the last corollary, if $\partial D$ does not
contain analytic discs through $a$ (by Proposition \ref{type7}
this is equivalent to $\partial D$ not containing linear discs
through $a$).

\begin{proposition}\label{locs-kob}
Let $a$ be a boundary point of a bounded domain $D\subset\C^n,$ as
in Proposition \ref{loc-con}. Also assume that $\partial D$ does
not contain analytic discs through $a.$ Then
$$\frac{1}{4}\le\liminf_{z\to a}\kappa_D(z;X)d_D(z;X)\le\limsup_{z\to a}\kappa_D(z;X)d_D(z;X)\le
1$$ uniformly on $X\in(\CC^n)_\ast.$
\end{proposition}

As in $\partial D$ there are no analytic discs through $a$ by
Propositions \ref{bound2}, \ref{bound3} and \ref{loc-con} imply
that for each sufficiently small neighborhood $U$ of $a$ we have
$\lim_{z\to a}d_D(z;X)=\infty$ uniformly on all unit vectors in
$\C^n.$ Then by shrinking $U$ (if necessary), $d_D(z;X)=d_{D\cap
U}(z;X)$ for each $z$ near $a$ (also $E_D(z;X)=E_{D\cap U}(z;X)$
and $A_D(z;X)=A_{D\cap U}(z;X)$).

After these remarks, Proposition \ref{locs-kob} follows from the following
strict localization for the Kobayashi metric (cf. Proposition \ref{c.pr0}).

\begin{proposition}\label{loc-kob}
Let $D\subset\C^n$ be a bounded domain, whose boundary does not
contain nontrivial analytic discs through a point $a\in\partial
D.$ Suppose that there exists a neighborhood $U$ of $a$ and a
function $f\in\mathcal O(D\cap U)$ such that $\lim_{z\to
a}|f(z)|=\infty.$ Then for each neighborhood $V$ of $a$ we have
$$\lim_{z\to a}\frac{\kappa_{D\cap V}(z;X)}{\kappa_D(z;X)}=1$$
uniformly on $X\in(\CC^n)_\ast.$
\end{proposition}

\beginproof Using the condition for the discs, as in
the proof of Proposition \ref{c.pr0}, it follows that each sequence
of analytic discs $\varphi_j$ , for which $\varphi_j(0)\to a,$
converges to $a$ uniformly (on the compact subsets of $\D$).
Then the proposition is contained in \cite[Corollary 2.3.4]{Nik0}.\qed
\smallskip

As a planar domain, having at least two points in its boundary, is
hyperbolic (see e.g. \cite{Jar-Pfl1}), the proof of the above
proposition shows that it can be strengthened for $n=1.$

\begin{proposition}\label{loc-kob-pla} Let $a$ be a boundary point
of a domain $D\subsetneq\C\setminus\{a\}.$ Then for each neighborhood
$V$ of $a$ we have $$\lim_{z\to a}\frac{\kappa_{D\cap V}(z;1)}{\kappa_D(z;1)}=1.$$

In particular, if $a$ is an isolated boundary point of $D,$ then
$$\lim_{z\to a}\kappa_D(z;1)|z|\log|z|=-\frac{1}{2}.$$
\end{proposition}

Note that the last equality follows from the formula
$\kappa_{\D_\ast}(z;1)=-\frac{1}{2|z|\log|z|}$ (see e.g.
\cite{Jar-Pfl1}).

Proposition \ref{loc-kob-pla} generalizes essentially, and in a short way, \cite[Theorem 1]{KLZ}.

\setcounter{equation}{0}
\section{Localization of the Bergman kernel and the Bergman metric}\label{loc-ber}

In this section we provide localization theorems for the Bergman
kernel and the Bergman metric, that together with Proposition
\ref{loc-con} allow us to localize the results from the preceding
sections that deal with the Bergman kernel and the Bergman metric.

In the case when $D\subset\Bbb C^n$ is a bounded pseudoconvex
domain, the corresponding results are well known (see e.g. \cite{DFH}).

\begin{theorem}\label{clas} Let $V\Subset U$
be neighborhoods of a boundary point $z_0$ of a bounded
pseudoconvex domain $D\subset C^n.$ Then there exists a constant
$c\ge 1$ such that for each $z\in D\cap V$ and for each $X\in\Bbb
C^n$ we have
$$c^{-1}K_{D\cap U}\le K_D(z)\le K_{D\cap U}(z),$$
$$c^{-1}B_{D\cap U}(z;X)\le B_D(z;X)\le c B_{D\cap U}(z;X).$$
\end{theorem}

By imitating the proof of Corollary \ref{loc-cor}, we get

\begin{corollary}\label{loc-ber1}
Let $a$ be a boundary point of a bounded domain $D\subset\C^n$,
as in Proposition \ref{loc-con}. Then
$$K_D(z)\sim_D 1/p^2_D(z),\quad B_D(z;X)\sim_D 1/d_D(z;X)\sim_D E_D(z;X)$$
near $a.$
\end{corollary}

Note that for the Bergman kernel, the localization is strict if
$z_0\in\partial D$ is a holomorphic peak point in the most general
sense, i.e. there exists a function $p\in\mathcal O(D,\D)$ such
that
$$\lim_{z\to z_0}p(z)=1>\sup_{D\setminus U}|p|$$
for each neighborhood $U$ of $z_0.$ This is proved in the
fundamental work \cite{Hor1} of L. H\"ormander as an application
of the $L^2$-estimates for the $\overline\partial$-problem. More
general results can be found in \cite{Her2}.

One of the goals of this section is to carry over this result in
the case of an arbitrary pseudoconvex domain (not necessarily
bounded). We will say that $z_0\in\partial D$ is a locally
holomorphic peak point, if $z_0$ is a holomorphic peak point of
$D\cap U$ for some neighborhood $U$ of $z_0.$

\begin{theorem}\label{peak} Let $U$ be a neighborhood of a boundary locally
holomorphic peak point $z_0$ of a pseudoconvex domain $D\subset C^n.$ Then
$$\lim_{z\to z_0}\frac{K_{D\cap U}(z)}{K_D(z)}=1,$$
$$\lim_{z\to z_0}\frac{B_{D\cap U}(z;X)}{B_D(z;X)}=1$$
uniformly on $X\in(\CC^n)_\ast.$
\end{theorem}

In particular, $K_{D\cap U}(z)>0$ and so $B_{D\cap U}(z;X)$
exists for $z$ close to $z_0.$
\smallskip

To prove Theorem \ref{peak}, we need a localization
lemma for the pluricomplex function of Green $g_D$
(for the definition see Section \ref{desc}).

\begin{lemma}\label{green} Let $U$ ne a neighborhood of a locally
holomorphic peak point $z_0$ of a pseudoconvex domain $D\subset
C^n.$ Then
$$\liminf_{z\to z_0,w\in D\setminus U}g_D(z,w)=0.^*\footnote{*Here and further we assume $D\setminus U\neq\varnothing.$}$$

In addition, there exists a neighborhood $V\subset U$ of $z_0$ such that
$$\inf\{z\in D\cap V,w\in D\cap U\setminus\{z\}:g_D(z,w)-g_{D\cap U}(z,w)\}=0.$$
In particular, we have strong localization for the Azukawa metric:
$$\lim_{z\to z_0}\frac{A_{D\cap U}(z;X)}{A_D(z;X)}=1$$
uniformly on $X\in(\CC^n)_\ast.$
\end{lemma}

The first equality means that $D$ has the so-called property $(P)$
(see e.g. \cite{Com1}), that is applied in problems about Bergman
invariants, as well in pluripotential theory.
\smallskip

\beginproof We use that each locally holomorphic peak point is a
plurisubharmonic peak point, i.e. there exists a negative function
$\varphi\in\PSH(D)$ such that
$$\lim_{z\to z_0}\psi(z)=0>\inf_{D\setminus U_1}\psi$$
for each neighborhood $U_1$ of $z_0.$ Indeed, one can assume that
$p$ is a holomorphic peak function on $D\cap U_1$ at $z_0.$ Then
it suffices to choose a neighborhood $U_2\Subset U_1$ of $z_0$
such that $G=D\cap U_1\setminus U_2\neq\varnothing$ and to put
$\delta=\sup_G|p|,$ $$\varphi=-1+\left\{\begin{array}{ll}
\max(\delta,|p|)&D\cap U_2\\
\delta&D\setminus U_2\end{array}\right.$$ (This argument shows
that the notion "plurisubharmonic peak point"\ has local
character.)

On the other hand, each plurisubharmonic peak point is a
plurisubharmonic antipeak point (see e.g. \cite{Gau2}), i.e. there
exists a negative function $\psi\in\PSH(D)$ such that
$$\lim_{z\to z_0}\psi(z)=-\infty<\inf_{D\setminus U_1}\psi$$
for each neighborhood $U_1$ of $z_0;$ for this purpose we put
$\psi=-\log(-\varphi).$

Now we pass essentially to the proof. Let us first suppose that
$z_0$ is a holomorphic peak point of $D\cap U.$ Let
$W\Subset U$ be a neighborhood of $z_0.$ We can choose another neighborhood
$V\Subset W$ of $z_0$ so that $\inf_{D\setminus W}\psi\ge
c:=1+\sup_{D\cap\partial V}\psi.$ Fix $z\in D\cap V$ and put
$d(z)=\inf_{w\in D\cap\partial V}g_{D\cap U}(z,w),\
u(z,w)=(c-\psi(w))d(z)$ for $w\in D.$ As $u(z,w)\le g_{D\cap
U}(z,w)$ for $w\in D\cap\partial V$ and $u(z,w)\ge 0>g_{D\cap
U}(z,w)$ for $w\in D\cap U\setminus W,$ the function
$$v(z,w)=\left\{\begin{array}{ll}
g_{D\cap U}(z,w),&w\in D\cap V\\
\max\{g_{D\cap U}(z,w),u(z,w)\},&w\in D\cap W\setminus V\\
u(z,w),&w\in D\setminus W
\end{array}\right.$$
is plurisubharmonic on the second variable and has a logarithmic
singularity at $z.$ Also, $v(z,w)<cd(z)$ and so $g_D(z,w)\ge
v(z,w)-cd(z).$ As $v(z,w)=u(z,w)\ge 0$ for $w\in D\setminus W,$ we
get that $g_D(z,w)\ge -cd(z)$ for $w\in D\setminus W.$ Since
$$g_{D\cap U}(z,w)\ge\left|\frac{p(w)-p(z)}{1-\overline{p(z)}p(w)}\right|,$$
$\lim_{z\to z_0}d(z)=0$ and then $\lim_{z\to z_0}\inf_{w\in D\setminus W}g_D(z,w)=0,$
which proves the first equality in the lemma.

Let now $U$ be an arbitrary neighborhood of $z_0.$ We repeat the above
considerations. Using the first equality in the lemma for $V$ instead of
$U$ and the inequality $g_{D\cap U}\ge g_D,$ we get
$\lim_{z\to z_0}d(z)=0.$ Then the equality $v(z,w)=g_{D\cap
U}(z,w)$ for $w\in D\cap V$ implies the second equality in the lemma.\qed
\smallskip

\noindent{\bf Remark.} The above proof shows that for each
neighborhood $U$ of a plurisubharmonic antipeak point $z_0$, there
exists neighborhood $V\subset U$ of $z_0$ so that
$$\lim_{z\to z_0}\inf\{w\in D\cap V\setminus\{z\}:g_D(z,w)-g_{D\cap U}(z,w)\}>-\infty.$$
In particular, we have the following weak localization for the Azukawa metric:
for each neighborhood $U$ of $z_0$, there exist a constant
$C>0$ and a neighborhood $V\subset U$ of $z_0$ so that for each $z\in
D\cap V$ and for each $X\in\Bbb C^n$ we have
$$C^{-1}A_{D\cap U}(z;X)\le A_D(z;X)\le A_{D\cap U}(z;X).$$
Note that each boundary point of a bounded domain is a
plurisubharmonic antipeak point, as shown by the function
$\frac{\log(||z-z_0||}{\diam(D)}$).
\smallskip

A key role in the proof of Theorem \ref{peak} will be played by
the following lemma (replacing the existence of a bounded strictly
plurisubharmonic functions on bounded domains).

\begin{lemma}\label{strong} For each plurisubharmonic antipeak point
$z_0$ of open set $D\subset\CC^n$, there exists a neighborhood $V$
containing it, a number $c>0$ and a bounded function $s\in PSH(D)$
such that $-1<s\le 0$ and the function $s(z)-c||z||^2$ is
plurisubharmonic in $D\cap V.$
\end{lemma}

\beginproof Let $\varphi$ be a plurisubharmonic antipeak function for $z_0,$
and $W$ be a bounded neighborhood of $z_0$ such that
$D\cap\partial W\neq\varnothing.$ Then
$$m=\inf_{D\cap\partial W}(\varphi-||\cdot-z_0||^2)>-\infty$$
and consequently
$$\tilde s=\left\{\begin{array}{ll}
\max\{\varphi,||\cdot-z_0||^2+m\},&D\cap W\\
\varphi,&D\setminus W
\end{array}\right.$$
is a bounded plurisubharmonic function on $D$ that coincides with
$||\cdot-z_0||^2+m$ in some neighborhood $V$ of $z_0.$ It remains
to put $$s=\frac{\sup_D\tilde s-\tilde s}{\sup_D\tilde s-\inf_D
\tilde s}.\qed$$

\noindent{\it Proof of Theorem \ref{peak}}. Recall that
$$B_D(z;X)=\frac{M_D(z;X)}{\sqrt{K_D(z)}},
$$
where $M_D(z;X)=\sup\{|f'_z(X)|:f\in
L_h^2(D),\,\|f\|_D=1,\;f(z)=0\}$.

We will only prove that
$$\lim_{z\to z_0}\frac{M_{D\cap U}(z;X)}{M_D(z;X)}=1$$
uniformly on $X\in(\CC^n)_\ast.$ The proof of the equality
$$\lim_{z\to z_0}\frac{K_{D\cap U}(z)}{K_D(z)}=1$$ is analogical
(even simpler) and we omit it. These two equalities imply the theorem.

By shrinking (if necessary) the neighborhood $V$ in Lemma
\ref{strong}, we can assume that $V\subset U$ and that there exists a
locally holomorphic peak function $p$ for $z_0,$ defined on
$D\cap V.$ Let $\chi$ be a smooth function with support in $V$ such that
$0\le\chi\le 1$ and $\chi\equiv 1$ in a neighborhood $V_1\Subset V$ of $z_0.$
By Lemma \ref{green}, there exists a neighborhood $V_2\Subset
V_1$ of $z_0$ so that
$$m=\inf\{g_D(z,w):z\in D\cap V_2, w\in D\setminus V_1\}>-\infty.$$
For $k\in\NN,$ $z\in D\cap V_2$ and $f\in L_h^2(D\cap U)$ such that
$f(z)=0$, put $\alpha={\overline\partial}(\chi fp^k)$ and continue $\alpha$
trivially as a $\overline\partial$-closed
$(0,1)$-form on $D.$ Let
$$\beta=\exp(-2(n+j)g_D(z,\cdot)-s),$$ where $s$ is the function from Lemma
\ref{strong}. As $-\log\beta-c||\cdot||^2$ is a plurisubharmonic function on the open set $\{w\in
D:\alpha(w)\neq 0\},$ from the proof of \cite[Theorem
2.2.1']{Hor1} it follows that there exists a smooth function $h$ on $D$ such that
$\overline\partial h=\alpha$ and $$\int_D|h|^2\beta\le c^{-1}\int_D|\alpha|^2\beta.$$
Then $g=\chi fp^k-h$ is a holomorphic function on $D.$ As the right-hand side
of the above inequality is bounded, so is the left-hand side. Then $h(z)=0$ and hence
$$g'_z(X)=(p(z))^kf'_z(X).$$

In addition, from $g_D<0$ and $s<0$ it follows that
$$||h||_D^2\le\int_D|h|^2\beta.$$ On the other hand, if
$C=\exp(-2(n+j)m+1)\sup|{\overline\partial\chi}|^2$ and $q=\sup_{D\cap V\setminus V_1}|p|,$ then
$$\int_D|\alpha|^2\beta\le Cq^{2k}.$$
Putting $C_1=\sqrt{C/c},$ the last three inequalities imply
$$||g||_D\le 1+C_1q^k.$$

Now the definition of $M_D$ implies that
$$M_{D\cap U}(z;X)\ge M_D(z;X)\ge{M_{D\cap U}(z;X)|p(z)|^{2k}\over(1+C_1 q^k)^2}.$$
Leaving $z\to z_0,$ then $k\to\infty$ and using that $\lim_{z\to z_0}p(z)=1$ and $q<1,$
we get the desired equality.\qed
\smallskip

From the above proof (for $k=0$) we get

\begin{corollary}\label{antipeak}
If $U$ is a neighborhood of a plurisubharmonic antipeak point $z_0$ of an (arbitrary)
domain $D\subset\CC^n,$
then there exist a constant $c>0$ and a neighborhood $V\subset U$ of $z_0$ so that
$$c^{-1}K_{D\cap U}\le K_D(z)\le K_{D\cap U}(z),$$
$$c^{-1}B_{D\cap U}(z;X)\le B_D(z;X)\le c B_{D\cap U}(z;X)$$
for each $z\in D\cap V$ and for each $X\in\Bbb C^n.$

In particular, such a localization holds for an arbitrary boundary
point $z_0$ of a domain $D\subset\CC$ whose complement is not a
polar set.
\end{corollary}

Recall that a set $E\subset C$ is called polar if $E\subset u^{-1}(-\infty)$ for some $-\infty\not\equiv
u\in\SH(D).$ If the complement of a domain $D\subset\CC$ is not polar, then $K_D>0;$ otherwise $K_D\equiv 0.$

To see that $z_0$ is a subharmonic antipeak point of $D,$
it suffices to note that for each sufficiently small neighborhood
$V$ of $z_0$, the complement of $G=D\cup V$ is not polar. Then
$g_G(z_0,\cdot)$ is a bounded function on $G$ outside an arbitrary
neighborhood of $z_0$ and so it is a subharmonic antipeak function for $z_0.$

Corollary \ref{antipeak} can be applied for proving that the completeness of the Bergman distance of a
planar domain with a non-polar complement has local character (see the work of the author \cite{Nik2}).
\smallskip

To apply Theorem \ref{peak}, note that if $a$ is a boundary point
of a domain $D\subset\C^n$ as in Proposition \ref{loc-con}, and in
addition $a$ is of finite type, then it is a locally holomorphic
peak point. Indeed, as noted at the end of Section \ref{type}, $a$
is a semiregular point and then it suffices to use the main result
in \cite{Yu2}. A more general result in the smooth $\C$-convex
case can be found in \cite{Die-For2}. So we get

\begin{corollary}\label{loc-ber0}
Let $a$ be a smooth boundary point of finite type of a (not necessarily bounded)
 domain $D\subset\C^n$ as in Proposition \ref{loc-con}.
Then $$K_D(z)\sim 1/p^2_D(z),\quad B_D(z;X)\sim 1/d_D(z;X)\sim E_D(z;X).$$
near $a.$
\end{corollary}

Recall that the constants in $\sim$ depend only on $n.$ This
corollary essentially strengthens some of the main results in \cite{Chen,McN2,McN3,Blu}.

The next proposition allows us to sharpen this result, as well as
Corollary \ref{locs-kob} in the case of convex domains. In less
generality this proposition is formulated in \cite{Sib} with only a hint for a proof.

\begin{proposition}\label{con-peak} Let $D\subset\CC^n$ be a convex
domain. Then $a\in\partial D$ is a holomorphic peak point exactly when $L_a=\{0\}.$
\end{proposition}

\beginproof The necessity of the condition $L_a=\{0\}$ is almost
obvious. Indeed, suppose that there exists a holomorphic peak
function $f$ for $D$ at $a,$ but $L_a\neq\{0\}.$ By the remark after the end
of the proof of Proposition \ref{bound3}, one can find a
vector $X\neq 0,$ number $\varepsilon>0$ and a sequence of points
$z_j\to a$ so that $\Delta_X(z_j,\eps)\subset D.$ Then,
considering the restriction of $f$ on the complex line through
$z_j$ in the direction of $X,$ we get a contradiction with the maximum principle.

Let now $L_a=\{0\}.$ Then $D$ does not contain complex lines;
otherwise $D$ would be biholomorphic to $\C\times D',$ and the
corresponding biholomorphism would be continued in a neighborhood
of $a$ (see the proof of Proposition \ref{c.pr4}) and so $D$ would
contain analytic discs, and so also linear discs through $a$ (see
the remark after the proof of Proposition \ref{type7}); a
contradiction. Consequently $D$ is biholomorphic to a bounded
domain and the corresponding biholomorphism is continued in a
neighborhood of $a$ (see the proof of Proposition \ref{c.pr4}).
Thus we can suppose that $D$ is bounded domain. Note that if $c$
is a positive number such that $\ds c\inf_{z\in D}\Ree(z_1)>-1$ (D
is bounded), then the function $f_1(z)=\exp(z_1+cz_1^2)$ belongs
to $A(D)=\O(D)\cap C(\overline D)$ and $|f_1(z)|<1$ for
$z\in\overline D\setminus\{z_1=0\}$. This easily implies (cf.
\cite{Gam}) that $\supp\ \mu\subset D_1:=\partial D\cap\{z_1=0\}$.
Since $L(0)=0$, the origin is a boundary point of the compact
convex set $D_1$. As above, we may assume that
$D_1\subset\{z\in\CC^n:\Ree(z_2)\le 0\}$ ($z_2$ is independent of
$z_1$) and then construct a function $f_2\in A(D)$ such that
$|f_2(z)|<1$ for $z\in D_1\setminus\{z_2=0\}$. This implies that
$\supp\ \mu\subset D_1\cap\{z_2=0\}$. Repeating this argument we
conclude that $\supp\ \mu=\{0\}$, i.e. $0$ is peak point for the
algebra $A(D)$ (see e.g. \cite{Gam}), which even means that there
exists a function $f\in A(D)$ such that $f(a)=1$ and $|f(b)|<1$
for each point $\overline D\ni b\neq a.$ \qed

\begin{corollary}\label{loc-ber2}
Let the pseudoconvex domain $D\subset\C^n$ be locally convex
near its boundary point $a$. If $\partial D$ does not contain analytic
discs through $a$, then $$K_D(z)\sim 1/p^2_D(z),\quad B_D(z;X)\sim\kappa_D(z;X)\sim 1/d_D(z;X)\sim E_D(z;X)$$
near $a.$
\end{corollary}

The estimate for $\kappa_D(z;X)$ follows from the strong localization for the Kobayashi metric near a
locally holomorphic peak point of an arbitrary (not necessarily bounded) domain (see e.g.
\cite[Theorem 2.3.9]{Nik0}).
\smallskip

\noindent{\bf Remark.} Corollary \ref{loc-ber2} immediately implies that
under these assumptions we get
\begin{equation}\label{her}
\lim_{z\to a}B_D(z;X)=\infty
\end{equation}
locally uniformly on $X\in(\C^n)_\ast.$ Thus we carry over (in an easy way) the main result from \cite{Her1}
even for unbounded domains (the proof there is based on the $\overline\partial$-technique of Ozawa--Takegoshi).
In the case of bounded pseudoconvex domains that are locally $\C$-convex near
$a,$ the equality (\ref{her}) also remains true due to Theorem
\ref{clas}. This is another strengthening of the mentioned result. On the other hand, using
the inequality $\gamma_D\le B_D$ and Proposition \ref{bound3} we can "reverse"\ the above considerations,
i.e. from (\ref{her}) to get $L_a=\{0\}.$

\setcounter{equation}{0}
\section{Boundary behavior of invariant metrics of planar domains}\label{plane}

After discussing the boundary behavior of the invariant metrics of domains in $\C^n,$
it is natural to see whether these results can be more precise for planar domains.
In this quite short section we will prove the following

\begin{proposition}\label{pla}
If $a_0$ is a $\mathcal C^1$-smooth boundary point of a domain
$D\subset\C$, then
$$\lim_{a\to a_0}\gamma_D(a;1)d_D(a)=\lim_{a\to
a_0}\kappa_D(a;1)d_D(a)=\frac{1}{2},$$
$$\lim_{a\to a_0}K_D(a)d^2_D(a)
=\frac{1}{4\pi}\hbox{ and  }\lim_{a\to
a_0}B_D(a;1)d_D(a)=\frac{\sqrt2}{2}.
$$
\end{proposition}

The condition for smoothness is essential,
as shown e.g. by the first quadrant.
\smallskip

\beginproof The proposition for the Carath\'eodory and Kobayashi metrics is
equivalent to the inequalities
\begin{equation}\label{plak}
\limsup_{a\to a_0}\kappa_D(a;1)d_D(a)\le\frac{1}{2},
\end{equation}
\begin{equation}\label{plac}
\limsup_{a\to a_0}\gamma_D(a;1)d_D(a)\ge\frac{1}{2}.
\end{equation}

Inequality (\ref{plak}) is contained in \cite[p. 60]{Nik0} in a more general situation
(we are not including its proof here): Let $a_0$ be
a $\mathcal C^1$-smooth boundary point of a domain $D\subset\C,$
$X_a\to X$ for $a\to a_0$. If $X_N$ is the projection of $X$ onto
the complex normal to $\partial D$ at $a_0,$ then
\begin{equation}\label{plag}
\limsup_{a\to a_0}\kappa_D(a;X_a)d_D(a)\le\frac{||X_N||}{2}.
\end{equation}

Now we will prove the less trivial inequality (\ref{plac}) (via
the Pinchuk scaling method).

We can assume that $a_0=0.$ For each point $a\in D$ close to
$0,$ there exists a point $\wdht a\in\partial D$ such that $\|a-\wdht
a\|=d_D(a)$ and $a$ lies on the internal normal to $\partial D$ at
$\wdht a$. Let $r$ be a $\mathcal C^1$-smooth defining function for $D$ near $0.$
Put $\Phi_a(z)=\frac{\partial r}{\partial
z}(\wdht a)(\wdht a-z)$. Let also
$$
E_\eps=\{z\in\CC:\Re z>-\eps|z|\},\quad
F_\eps=\{z\in\CC:|z|>\eps\}.
$$
For each sufficiently small $\eps>0$ we have $\Phi_a(D)\subset
E_\eps\cup F_\eps$ for $|a|<\eps$. As $\wdtl a=\Phi_a(a)>0$,
\begin{equation}\label{pla1}
\gamma_D(a;1)\ge\gamma_{E_\eps\cup F_\eps}(\wdtl a;X(a))=
\gamma_{G_{\eps,a}}(1;1)\frac{|X(a)|}{\wdtl
a}=\frac{\gamma_{G_{\eps,a}}(1;1)}{d_D(a)},
\end{equation}
where $X(a)=-\frac{\partial r}{\partial z}(\wdht a)$ and
$G_{\eps,a}=E_\eps\cup F_{\frac{\eps}{\wdtl a}}$.
Note that
\begin{equation}\label{pla2}
\lim_{a\to a_0}\gamma_{G_{\eps,a}}(1;1)=\gamma_{E_\eps}(1;1)
\end{equation} and \begin{equation}\label{pla3}
\lim_{\eps\to
0+}\gamma_{E_\eps}(1;1)=\gamma_{E_0}(1;1)=\frac{1}{2}.
\end{equation}
Then (\ref{plac}) follows from (\ref{pla1}), (\ref{pla2}) and (\ref{pla3}).

To prove (\ref{pla2}), we denote by $H_\eps$ and $H_{\eps,a}$ the
images of $E_\eps$ and $G_{\eps,a}$, respectively, for the mapping
$z\to2/(z+1),$ if $\wdtl a<\eps<1$. Then $H_\eps$ and $\wdtl
H_{\eps,a}=H_{\eps,a}\cup\{0\}$ are bounded simply connected
domains and consequently $C_{H_\eps}=K_{H_\eps}$ and
$C_{H_{\eps,a}}=C_{\wdtl H_{\eps,a}}=K_{\wdtl H_{\eps,a}}$. Now,
using the Montel theorem, it is easily seen that
$$\lim_{a\to a_0}K_{\wdtl H_{\eps,a}}(1;1)=K_{H_\eps}(1;1),$$
which implies (\ref{pla2}).

The equality (\ref{pla3}) is proven in a similar way (or by using
that $E_{\eps}$ and $E_0$ are conformally equivalent).

The proposition for the Bergman kernel and Bergman metric is obtained
analogically (having in mind that $B_D(z;1)=\frac{M_D(z;1)}{\sqrt{K_D(z)}}$)
and we omit the proof.\qed
\smallskip

\noindent{\bf Remark.} Under the somewhat stronger requirement that the boundary
be Dini-smooth near $a_0,$ the proposition for the Bergman kernel, as well as for
the metrics of Bergman and Kobayashi, can be also proven using that:

$\bullet$ each $\mathcal C^1$-smooth boundary point $a$ of a domain
$D\subset C$ is a locally holomorphic peak point and so we have strong
localization for these invariants;

$\bullet$ there exists a neighborhood $U$ of $a$ so that $G=D\cap
U$ is a bounded Dini-smooth simply connected domain and so the
Riemann mapping between $G$ and $\D$ is continued to a $\mathcal
C^1$-diffeomorphism between $\overline G$ and $\overline\D.$

\end{document}